\documentclass[twocolumn,secnumarabic,amssymb,superscriptaddress,nobibnotes,aps,prx,longbibliography]{revtex4-2}

\usepackage{amsmath,amsfonts,bm,amssymb,mathtools}
\usepackage{graphicx}
\usepackage{float}
\usepackage{subfigure}
\usepackage{color}
\usepackage{algorithm}
\usepackage{algpseudocode}
\usepackage[colorlinks=true,citecolor=blue,linkcolor=blue,urlcolor=blue]{hyperref}
\usepackage{ntheorem}
\usepackage{array}
\usepackage{multirow}
\usepackage{makecell}
\usepackage{booktabs}

\newtheorem{definition}{Definition}[section]
\usepackage{bibunits}

\newcommand{\rvx}{\mathbf{x}}
\newcommand{\rvy}{\mathbf{y}}
\newcommand{\rve}{\mathbf{e}}
\newcommand{\rvs}{\mathbf{s}}
\newcommand{\rvt}{\mathbf{t}}
\newcommand{\rvv}{\mathbf{v}}
\newcommand{\rvz}{\mathbf{z}}
\newcommand{\rvb}{\mathbf{b}}
\newcommand{\rvf}{\mathbf{f}}
\newcommand{\rvn}{\mathbf{n}}
\newcommand{\rvh}{\mathbf{h}}
\newcommand{\vtheta}{\boldsymbol{\theta}}
\newcommand{\vzero}{\boldsymbol{0}}
\newcommand{\rvg}{\mathbf{g}}
\newcommand{\rmF}{\mathbf{F}}
\newcommand{\rmP}{\mathbf{P}}
\newcommand{\rmA}{\mathbf{A}}
\newcommand{\rmD}{\mathbf{D}}
\newcommand{\rmB}{\mathbf{B}}
\newcommand{\rmZ}{\mathbf{Z}}
\newcommand{\rmY}{\mathbf{Y}}
\newcommand{\rmG}{\mathbf{G}}
\newcommand{\rmW}{\mathbf{W}}
\newcommand{\rmS}{\mathbf{S}}
\newcommand{\rmH}{\mathbf{H}}
\newcommand{\rmC}{\mathbf{C}}
\newcommand{\rmR}{\mathbf{R}}
\newcommand{\rmI}{\mathbf{I}}
\usepackage{titletoc}
\titlecontents{suppsection}
[0pt]                              
{\addvspace{0.5em}\bfseries}       
{\contentslabel{2.2em}}            
{}                                 
{\titlerule*[0.7pc]{.}\contentspage}

\titlecontents{suppsubsection}
[2em]                              
{}
{\contentslabel{2.2em}}            
{}
{\titlerule*[0.7pc]{.}\contentspage}
\begin{document}


	\title{{\color{black}Ptolemy’s Equant Equates to a Universal Dynamical Clock via Machine Learning}}
	\author{Jingdong Zhang}
	\affiliation{School of Mathematical Sciences, Fudan University, Shanghai 200433, China.}
	\affiliation{Research Institute of Intelligent Complex Systems, Fudan University, Shanghai 200433, China}
	\affiliation{I-X
		and Department of Mathematics, Imperial College London, London, SW7 2AZ, United Kingdom.}
	\affiliation{Institute for Complex Systems and Mathematical Biology, University of Aberdeen, Aberdeen AB24 3UE, United Kingdom.}
	
	\author{Luan Yang}
	\affiliation{Research Institute of Intelligent Complex Systems, Fudan University, Shanghai 200433, China}
	\affiliation{Department of Psychiatry, University of Cambridge, Cambridge CB2 1TN, United Kingdom}

	\author{Murilo~S.~Baptista}\email{murilo.baptista@abdn.ac.uk}
\affiliation{Institute for Complex Systems and Mathematical Biology, University of Aberdeen, Aberdeen AB24 3UE, United Kingdom.}
	
	\author{Zefeng Zhang}
	\affiliation{Research Institute of Intelligent Complex Systems, Fudan University, Shanghai 200433, China}

	\author{Qunxi Zhu}\email{qxzhu@fudan.edu.cn}
	\affiliation{Research Institute of Intelligent Complex Systems, Fudan University, Shanghai 200433, China}
	
	\author{Wei Lin}\email{wlin@fudan.edu.cn}
	\affiliation{School of Mathematical Sciences, Fudan University, Shanghai 200433, China.}
	\affiliation{Research Institute of Intelligent Complex Systems, Fudan University, Shanghai 200433, China}
	
	\author{Celso Grebogi}
	\affiliation{Institute for Complex Systems and Mathematical Biology, University of Aberdeen, Aberdeen AB24 3UE, United Kingdom.}

	\date{\today}
	\begin{abstract}
		{\color{black}
			Oscillatory dynamics arise ubiquitously in nonlinear systems, yet identifying a physically interpretable phase and phase dynamics in nonlinear, high-dimensional oscillations remains a central unresolved problem. Here we establish the principle of a universal dynamical clock, a physical perspective in which oscillations of arbitrary dimensionality and geometry are equivalently {\color{black} represented as uniform rotation through an equant-induced nonlinear viewing coordinate}, {\color{black}inspired by Ptolemy’s equant and formalised through an areal-uniformity principle reminiscent of Kepler’s second law}. By proposing and implementing a machine-learning framework, we demonstrate the existence of such an equant for a broad class of oscillatory dynamics and construct the associated dynamical clock and phase dynamics under a wide range of additive forces, including noise, periodic perturbations, and coupling effects. The critical role of the dynamical clock in uncovering new physical rules and phenomena across diverse systems is demonstrated in various significant findings: (i)  collective oscillations in Escherichia coli populations obey a previously unexplained superlinear scaling law, resolving a long-standing open problem posed in 2004; (ii) the response mechanisms of engineered genetic circuits to changes in gene expression and environmental conditions; (iii) a classical-mechanics counterpart of the Berry geometric phase emerges naturally from the phase of the dynamical clock; and (iv) {\color{black}the optimal equant non-uniformity provides a geometric early-warning signal for critical transitions and enables prediction of critical parameters}.  By providing an operational and system-agnostic phase dynamics that can be constructed directly from data, the dynamical clock enables principled classification, comparison, and control of oscillatory systems, and offers a new route to understanding how specific dynamical regimes support distinct functional behaviours in networked systems.
		}

	\end{abstract}
	\maketitle
	
	\section{Introduction}
	
	{\color{black}Oscillatory dynamics are a fundamental mode of motion in nonlinear systems, arising across a wide range of natural and engineered contexts, from the firing of single neurons~\cite{hodgkin1952quantitative,hodgkin1952propagation} and collective behaviour of cellular populations~\cite{chen2017weak} to the dynamics of synthetic genetic circuits~\cite{brophy2014principles}. Despite their ubiquity, a basic theoretical question remains unresolved: how to define a phase and phase dynamics that are physically interpretable for nonlinear, high-dimensional oscillations. While phase is routinely invoked to characterise rhythmic behaviour and collective dynamics, its definition often relies on model-specific constructions or low-dimensional assumptions, such as classic Kuramoto models~\cite{kuramoto1984chemical,yeung1999time}. These assumptions obscure the physical meaning of phase when applied to complex oscillatory systems, leaving the relationship between observed oscillatory trajectories and an underlying, physically grounded phase description poorly understood. As a consequence, such models frequently fail to capture intrinsic frequencies and phase response mechanisms observed in real-world oscillations, where multidimensional irregularities and dynamical complexity give rise to experimental phenomena that remain unexplained within existing theoretical frameworks.


		To study the evolution of phase, phase reduction has been applied across diverse contexts, from experimental measurements of phase response curves in cyanobacterial circadian rhythms~\cite{kiyohara2005novel} to computer-aided phase response equations~\cite{ermentrout2003simulating}. Most existing approaches focus on how phase responds to perturbations, implicitly assuming the existence of a well-defined phase while leaving its physical interpretation largely unexamined. Common identification techniques, such as the adjoint and direct methods~\cite{glass1988clocks}, rely on local approximations or numerical procedures, limiting their reliability in high-dimensional and irregular oscillatory systems. As a result, phase reduction lacks a unifying structural principle, constraining its predictive power in complex settings including neural networks, cardiac tissue, and genetic circuits. This fundamental limitation has been clearly articulated by~\cite{pietras2019network}: “There is no unique method for a rigorous reduction to a phase model of networks of coupled oscillators, nor a straightforward recipe along which the phase dynamics should be reduced.”

	}

	\begin{figure*}[htp]
		\centering
		\includegraphics[width=\textwidth]{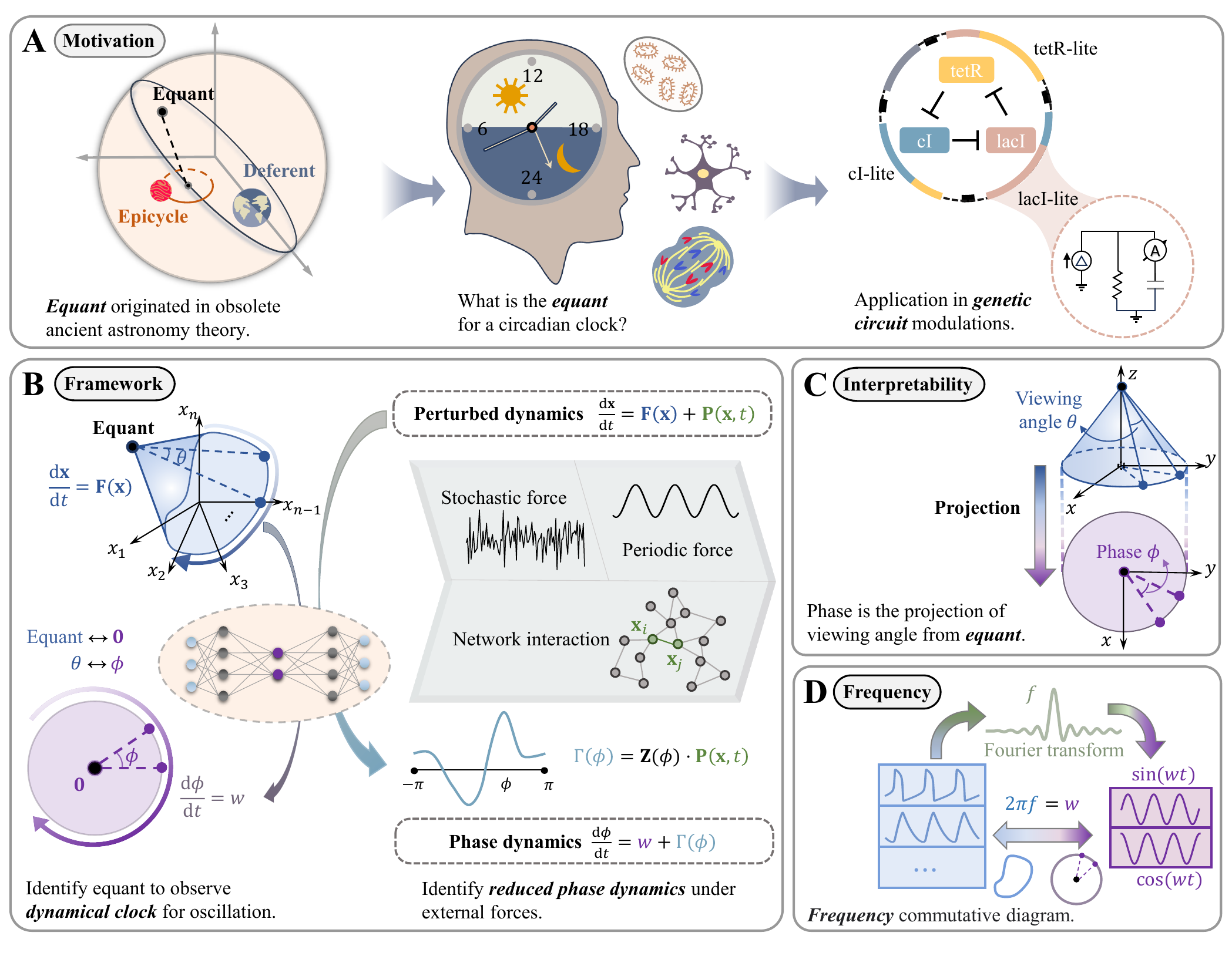}
		\caption{\textbf{Decoding Rhythmic Phenomena with Ptolemy’s Equant}. (\textbf{A}), 
			The deferent-and-epicycle model~\cite{evans2013mechanical} explained retrograde motion by placing each planet on a small epicycle revolving around a larger circle, the deferent (Fig.~\ref{fig1}A). Building on this model, Ptolemy introduced the \textbf{\emph{equant}}, an imaginary point from which an observer perceives the centre of the epicycle moving along the deferent at a constant angular speed. We identify the existence of equant for circadian clock, e.g., genetic clock, neuronal rhythm and mitosis, and apply it in genetic circuit design. (\textbf{B}),  The machine learning framework specifies the dynamical clock over the unit circle and equate it with a general oscillatory system, and identifies phase dynamics under various perturbations. (\textbf{C}), The phase can be physically interpreted as the projection of the viewing angle  observed from the equant. (\textbf{D}), The dominant Fourier frequency of oscillatory trajectories rescales the natural frequency of the dynamical clock. 
		}
		
		\label{fig1}
	\end{figure*} 
	
	The quest for simple principles to describe oscillatory motion dates back to ancient astronomy, where Greek scholars sought to reconcile the ideal of uniform circular motion with the observed retrograde trajectories of planets. To achieve this, Ptolemy introduced the concept of the equant, an imaginary point from which planetary orbits appear as constant-rate rotations~\cite{toomer1998ptolemy} (Fig.~\ref{fig1}A). This viewpoint anchored the “phase” of planetary motion to the angle of constant-rate rotation. Although Ptolemy’s geocentric view was later replaced by the heliocentric model, the underlying geometric insight remains instructive for understanding general oscillations. Here, we propose the principle of \textbf{\emph{dynamical clock}} to equate a large class of non-linear high-dimensional periodic oscillators with a two-dimensional uniform circular motion. Viewed from an equant-like perspective {\color{black}equipped with a  nonlinear viewing coordinate}, this clock has a pointer advancing at a constant angular speed set by the frequency of the original oscillation. Under this principle, phase is naturally the angular speed multiplied by elapsed time.

	With the rapid growth of artificial intelligence (AI), physics-informed machine learning has paved a new way for uncovering meaningful physical concepts and mathematical insights of natural phenomena~\cite{iten2020discovering,wang2023scientific,davies2021advancing,zhang2024learning,zhang2024machine}. Building on this progress, we develop a machine-learning framework to validate the dynamical clock. This framework provides a unifying perspective to decipher diverse oscillations and derives explicit reduced phase dynamics under realistic perturbations (Fig.~\ref{fig1}B).
	The key idea is to reinterpret a complex oscillation as a simplified two-dimensional clock, whose phase, termed the interpretable phase, is defined from the viewpoint of the equant. To build intuition, we consider a three-dimensional oscillator rotating uniformly in the $(x,y)$-plane (Fig.\ref{fig1}C). From an equant placed on the $z$-axis, the interpretable phase corresponds to the viewing angle of this rotation. This relation directly links the reduced phase dynamics to the dominant frequency of the original time series (Fig.~\ref{fig1}D).
	

	In this paper, we introduce the principle of universal dynamical clock, which provides an operational and geometrically grounded framework for defining phase dynamics in nonlinear oscillatory systems of arbitrary dimension. In Sec.~\ref{sec dynamical clock}, we formulate the dynamical clock as an equivalence between a general periodic oscillation and a two-dimensional uniform circular motion, thereby establishing a principled notion of phase with constant angular velocity (\ref{subsec step1}). Furthermore, we develop a data-driven construction that realises the dynamical clock directly from observational time series by combining manifold learning with invertible neural networks, (\ref{subsec step2}) yielding an interpretable and dynamically meaningful phase coordinate (\ref{subsec interpret}). 
	In Sec.~\ref{Methods phase reduction}, we identify the analytical formulation of phase dynamics for complex oscillatory systems perturbed by external forces, including  the stochastic force (\ref{subsec stochastic}), the
	periodic force (\ref{subsec period}), and the network interaction (\ref{subsec network}). In Sec.~\ref{sec results}, we demonstrate the power of this framework by  resolving a long-standing open problem in collective oscillations of Escherichia coli populations under quorum sensing (\ref{sec Ecoli}),  enabling precise modulation of synthetic genetic circuits based on the experimental data  (\ref{sec genetic}); {\color{black}establishing the classical-mechanics counterpart of the Berry geometric phase based on the dynamical clock~(\ref{appendix berry})}, {\color{black}revealing the role of early warning of critical parameters for different rhythms based on the  equant (\ref{sec dimension})}; and uncovering phase synchronisation mechanisms under complex perturbations  (\ref{sec FHN}). We conclude in Sec.~\ref{sec discussion} by discussing the broader implications of the dynamical clock for phase dynamics and collective behaviour in complex oscillatory systems.
	
	\section{Universal Dynamical Clock}\label{sec dynamical clock}
	\subsection{Problem setup}
	We consider a globally (or locally) attractive limit cycle oscillator described by a general non-linear dynamical system 
	\begin{equation}\label{eq1}
		\dfrac{\mathrm{d}\rvx}{\mathrm{d}t}=\rmF(\rvx),~\rvx\in\mathbb{R}^n,
	\end{equation} 
	wherein we denote by $\rmC$ the stable limit cycle.
	We assume there exists a periodic scalar function $\phi(\rvx)$ on $\rmC$ with period $2\pi$ such that its time derivative is a constant $w$, i.e.,
	\begin{equation*}
		\dfrac{\mathrm{d}\phi(\rvx)}{\mathrm{d}t}=w,~\phi\in[0,2\pi),
	\end{equation*}
	where $\phi$ and $w$ are known, respectively, as the phase function and the natural frequency~\cite{monga2019phase}, respectively. Such a phase function establishes a homeomorphism between the limit cycle and the unit circle, and
	we denote by $\boldsymbol{\chi}(\phi)\in\rmC$ the inverse map of phase. According to chain derivation, we have,
	\begin{equation*}
		\frac{\mathrm{d}\phi(\rvx)}{\mathrm{d}t}=\frac{\partial \phi}{\partial\rvx}\cdot\rmF(\rvx)=w,
	\end{equation*}
	which is used to construct loss function in our machine learning framework. In general systems the limit cycle $\rmC$ is a complex $1$-dimensional manifold without analytical formulation, which causes challenge for identifying $\phi$ and $\boldsymbol{\chi}$. For brevity, we denote $\frac{\mathrm{d}\rvx}{\mathrm{d}t}$ by $\dot{\rvx}$  below.

	{\color{black}
		
	To physically interpret the phase, otherwise a mathematical coordinate, we introduce the equant as an observational point $\rvx^\ast$ from which the viewing geometry of the limit cycle can be mapped to a uniform dynamical clock. The uniformity is therefore understood in the nonlinear coordinate induced by the equant and phase map, rather than as uniform Euclidean angular motion around $\rvx^\ast$. Motivated by Ptolemy’s equant and to avoid the non-uniqueness of purely topological viewing-angle conditions, we select this point through an areal-uniformity criterion inspired by Kepler’s second law: the swept area along the limit cycle should be as uniform as possible.

	\begin{definition}\label{def equant}
		We define the equant for the system in Eq.~\eqref{eq1} as a point $\rvx^\ast$ such that $\mathrm{(i)}$ for the viewing line $l_{\rvx^\ast\rvx(t)}$ connecting $\rvx^\ast$ and $\rvx(t)\in\rmC$, there
		exists a homeomorphism between the phase $\phi$ and the viewing angle $\theta$ induced by the viewing lines,
		$
		\theta(t, t')
		:=
		\angle\left(l_{\rvx^\ast\rvx(t)},l_{\rvx^\ast\rvx(t')}\right)
		$; and  $\mathrm{(ii)}$ it minimises the normalised areal variance, i.e.,
		\begin{equation}
			\rvx^\ast \in \arg\min_{\rvx} \hat{K}(\rvx),
			~
			\hat{K}(\rvx)=\frac{K(\rvx)}{\bar{A}^2(\rvx)} .
		\end{equation}
		Here, $\hat{K}(\rvx)$ is a dimensionless measure of the non-uniformity of the swept areas observed from $\rvx$. Specifically,
		\begin{equation*}
			\bar{A}(\rvx)
			=
			\frac{1}{T}\int_0^{T}
			\frac{1}{2}
			\left\|
			\big(\rvx-\rvx(t)\big)\wedge \rmF(\rvx(t))
			\right\|
			\mathrm{d}t
		\end{equation*}
		denotes the mean areal velocity over one period, and
		\begin{equation*}
			K(\rvx)
			=
			\frac{1}{T}\int_0^{T}
			\left|
			\frac{1}{2}
			\left\|
			\big(\rvx-\rvx(t)\big)\wedge \rmF(\rvx(t))
			\right\|
			-
			\bar{A}(\rvx)
			\right|^2
			\mathrm{d}t
		\end{equation*}
		denotes the variance of the swept areal velocity over one period.   $\rvx(t)\subset \rmC$ is the periodic solution with period $T$.
		We define the optimal equant non-uniformity as
		$
		\hat{K}^\ast=\min_{\rvx}\hat{K}(\rvx)$.
	\end{definition}

	According to Definition~\ref{def equant}, $\hat{K}$ is non-negative and continuous. Hence, provided that the admissible search domain of $\rvx$ is bounded, there exists an equant that minimises $\hat{K}$. Moreover, if the minimiser set $\arg\min_{\rvx}\hat{K}(\rvx)$ is a singleton, the equant is unique.

}
	
	In order to take perturbations ensuing from the real-world scenarios into consideration, we extend the domain of $\phi$ to the vicinity of $\rmC$. Given any phase value $\psi\in[0,2\pi)$, we define the isochron $\rmI_{\psi}$ as the manifold that intersects $\rmC$ with point $\boldsymbol{\chi}(\psi)$ such that any trajectory initiated from $\rmI_{\psi}$ converges to the same limit, i.e.,
	\begin{equation*}
		\lim_{t\to\infty}\Vert\rvx(t;0,\rvx_1)-\rvx(t;0,\rvx_2)\Vert=0,~\forall \rvx_1,\rvx_2\in \rmI_\psi.
	\end{equation*}
	The isochron can be regarded as the orthogonal space of the limit cycle, and thus it is a $(n-1)$-dimensional manifold.
	We can naturally extend $\phi$ to the isochrons $\{\rmI_\psi:\psi\in[0,2\pi)\}$ as $\phi(\rvx)=\psi$ for any $\rvx\in\rmI_\psi$~\cite{winfree1980geometry,winfree1974patterns,guckenheimer1975isochrons}. Then we consider the perturbed limit cycle oscillator described as,
	\begin{equation*}	\dot{\rvx}=\rmF(\rvx)+\rmP(\rvx,t),~\rvx\in\mathbb{R}^n,t\ge0,
	\end{equation*}
	where $\rmP(\rvx,t)$ is the external force representing the perturbation. The corresponding perturbed phase dynamics becomes
	\begin{equation*}
		\dot{\phi}=\dfrac{\partial \phi}{\partial\rvx}\cdot\left(\rmF(\rvx)+\rmP(\rvx,t)\right)=w+\dfrac{\partial \phi}{\partial\rvx}\cdot\rmP(\rvx,t).
	\end{equation*}
	By denoting the phase response curve as
	\begin{equation}\label{eq PRC}
		\rmZ(\phi)=\frac{\partial \phi}{\partial \rvx}\bigg|_{\rvx=\boldsymbol{\chi}(\phi)},
	\end{equation}
	we obtain the closed-form of phase dynamics on $\rmC$ as
	\begin{equation}\label{eq phase}
		\dot{\phi}=w+\Gamma(\phi,t),~\Gamma(\phi,t)=\rmZ(\phi)\cdot\rmP(\boldsymbol{\chi}(\phi),t).
	\end{equation}
	We require that the strength of perturbation is restricted to a certain range such that the perturbed dynamics fluctuate in the vicinity of $\rmC$, which is a fundamental assumption in phase reduction theory~\cite{brown2004phase,goldobin2010dynamics,NOVICENKO20121090}.
	We consider three major forms of perturbation:  $\text{(i)}$ the stochastic force driven by noise; $\text{(ii)}$ the periodic external force; and $\text{(iii)}$ and the interplay among oscillators $\{\rvx_i\}_{i=1}^N$ connected in a network. The detailed formulations of these perturbations are provided in Section~\ref{Methods phase reduction}.
	
	Given the general oscillatory dynamics in Eq.~\eqref{eq1}, the goal of our framework is {\color{black}threefold},
	\begin{itemize}
		\item [(i)] finding analytical expression of the phase $\phi(\rvx)$ and its inverse $\boldsymbol{\chi}(\phi)$ to obtain the homeomorphism between the general oscillation and the universal dynamical clock,
		\item [(ii)] discovering the equant $\rvx^\ast$ to interpret the physical meaning of the phase, and
		\item [(iii)] investigating the influence of perturbations to oscillation by obtaining the analytical formulation of perturbed phase dynamics in Eq.~\eqref{eq phase}.
	\end{itemize}
	
	\subsection{First step: Equating the limit cycle with the unit circle}\label{subsec step1}
	In order to obtain the phase function, we note that the limit cycle $\rmC$ is a $1$-dimensional manifold, and the phase $\phi(\rvx)\in[0,2\pi)$ belongs to the quotient space $\mathbb{R}/_\sim$ with an equivalent relation as 
	\begin{equation*}
		x\sim y\Leftrightarrow x-y=2k\pi,~k\in\mathbb{Z}.
	\end{equation*}
	We seek the phase function $\phi:\rmC\to \mathbb{R}/_\sim$ and its inverse $\boldsymbol{\chi}:\mathbb{R}/_\sim\to\rmC$ via using an auto-encoder neural network.  To integrate the quotient manifold structure into the machine learning framework, we note that $\mathbb{R}/_\sim$ is equivalently related to  the unit circle $\rmS^1=\{(x,y):x^2+y^2=1\}$ via the exponent map $e^{\text{i}\phi}$.  Therefore, we seek the phase function $\phi$ and its inverse $\boldsymbol{\chi}$ via an auto-encoder neural network with the latent space $\rmS^1$ as: 
	\begin{equation*}
		\begin{aligned}
			\text{(Encoder)}~&\rvx\in\boldsymbol{C}\to \rvy\in \rmS^1: \text{Enc}(\rvx)=\dfrac{\text{NN}_{\boldsymbol{\theta}_\phi}(\rvx)}{\|\text{NN}_{\boldsymbol{\theta}_\phi}(\rvx)\|},\\
			\text{(Decoder)}~&\rvy\in \rmS^1\to \rvx\in\rmC: \text{Dec}(\rvy)=\text{NN}_{\boldsymbol{\theta}_{\boldsymbol{\chi}}}(\rvy),\\
			&\phi_{\boldsymbol{\theta}_\phi}(\rvx)= \rvh \left(\text{Enc}(\rvx)\right),\\
			&\boldsymbol{\chi}_{\boldsymbol{\theta}_{\boldsymbol{\chi}}}(\phi)= \text{Dec}\left(\rvh^{-1}(\phi)\right),
		\end{aligned}
	\end{equation*}
	where Enc: $\rmC\to\rmS^1$ is the neural network (NN) parametrised by $\vtheta_\phi$, {\color{black}$\|\cdot\|$ is the $\mathrm{L}^2$ norm,} Dec: $\rmC\to\rmS^1$ is the neural network parametrised by $\vtheta_{\boldsymbol\chi}$, and $\rvh$ is the homeomorphism between $\rmS^1$ to $\mathbb{R}/_\sim$ as
	\begin{equation*}
		\begin{aligned}
			\rvh(x,y)&=\left\{\begin{aligned}
				&\arcsin y,~x\ge0,~y\ge0\\
				&\arccos x,~x<0,~y\ge0,\\
				&\pi-\arcsin y,~x\le0,~y<0\\
				&2\pi-\arccos x,~x>0,y<0.\\
			\end{aligned}\right.   \\
			\rvh^{-1}(\phi)&=(\cos(\phi),\sin(\phi)).
		\end{aligned}
	\end{equation*}
	We also set the natural frequency as a learnable quantity $w$.  Then, we train the parameters $\boldsymbol{\theta}_{\phi}$, $\boldsymbol{\theta}_{\boldsymbol{\chi}}$ and $w$ using the loss function $\mathcal{L}_{\text{phase}}$ as follows, 
	\begin{equation*}
		\mathcal{L}_{\text{phase}}(\boldsymbol{\theta}_{\phi},\boldsymbol{\theta}_{\boldsymbol{\chi}},w)=\lambda_1\mathcal{L}_1(\boldsymbol{\theta}_{\phi},w)+\lambda_2\mathcal{L}_2(\boldsymbol{\theta}_{\boldsymbol{\chi}}),
	\end{equation*}
	with
		\begin{equation}\label{eq mc loss1}
		\begin{aligned}
			\mathcal{L}_1(\boldsymbol{\theta}_{\phi},w)&=\dfrac{1}{K}\sum_{i=1}^K\bigg\Vert\dfrac{\partial\phi_{\boldsymbol{\theta}_{\phi}}}{\partial\rvx}\cdot \rmF(\rvx_i)-w\bigg\Vert_2^2,\\
			\mathcal{L}_2(\boldsymbol{\theta}_{\boldsymbol{\chi}})&=\dfrac{1}{K}\sum_{i=1}^K\|\rvx_i-\boldsymbol{\chi}_{\boldsymbol{\theta}_{\boldsymbol{\chi}}}\left(\phi_{\boldsymbol{\theta}_{\phi}}(\rvx_i)\right)\|_2^2,\\
		\end{aligned}
	\end{equation}
{\color{black}	Here $\lambda_{1,2}$ are the predefined regularisation parameters, and we collect the trajectory data $\{\rvx_i\}_{i=1}^K$ over the limit cycle $\rmC$ as training data. 	The first term  $\mathcal{L}_1$ in the loss function measures whether the auto-encoder maps the original dynamics near the limit cycle to the universal dynamical clock, and the second term establishes the analytical homeomorphism between the limit cycle and $\rmS^1$. }

	{\color{black}
		Recent machine-learning approaches have provided powerful tools for data-driven phase computation. In particular, phase-autoencoder methods learn circular latent variables whose evolution over a prescribed time interval $\tau$ follows a rigid rotation with frequency $w$~\cite{yawata2024phase,hiruta2025autoencoder}, while Koopman-based methods compute global isochrons and phase coordinates through spectral properties of the flow~\cite{mauroy2012use,mauroy2014global}. These approaches primarily focus on computing phase coordinates, isochrons, or phase sensitivity functions.
		
		In contrast, our framework learns the phase through the differential constraint in Eq.~\eqref{eq mc loss1}
		which can be imposed pointwise on limit-cycle samples with velocity information. Hence, the first step of our method does not require regularly sampled time-series data or off-cycle samples. This differs from previous approaches and classical adjoint method~\cite{ermentrout1996type}, which requires both $\rmF(\rvx)$ and $\nabla\rmF(\rvx)$ along the limit cycle and usually relies on time-parametrised orbit data. Although our main training procedure does not explicitly constrain transversal gradients away from the limit cycle, the learned phase map, inverse map, and PRC are sufficient for constructing the reduced phase dynamics considered in this work. In the Supplementary Material, we compare our ML-based PRCs with adjoint-computed PRCs and show comparable performance while using only limit-cycle data.

 
 For completeness, the Supplementary Material also presents an extended version of our framework in which transversal gradient information is incorporated into the training process. This extension covers both model-driven settings, where $\nabla \mathbf{F}$ is known, and data-driven settings, where only trajectory data $\{\mathbf{x}_i\}$ are available.

}

	\subsection{Second step: Discovering the  equant with the invertible neural network}\label{subsec step2}
	Given the equant $\rvx^\ast$ and the limit cycle $\rmC$, the line segments $l_{\rvx^\ast\rvx(t)}$ between the equant and the oscillator state on the limit cycle form a manifold  $\mathbb{M}=\{\tilde{\rvx}:\tilde{\rvx}\in l_{\rvx^\ast\rvx},\rvx\in\rmC\}$. 
	{\color{black}
	The raw Euclidean viewing angle $\theta$ observed from the equant is not generally identical to the phase $\phi$: its period may differ from $2\pi$, and its relation to $\phi$ can be nonlinear for irregular limit cycles.
}
	For universal dynamical clock on the unit circle, the centre of the unit disk $\mathbb{D}=\{(x,y):x^2+y^2\le1\}$ is the equant and naturally equates the phase as its viewing angle.  		
		{\color{black}
			To capture the relation between the phase and the equant for general oscillations, we extend the learned homeomorphism between the limit cycle $\rmC$ and the unit circle $\mathbb S^1$ to a map between the line-of-sight surface $\mathbb M$ and the unit disk $\mathbb D$. The equant in Definition~\ref{def equant} is selected by an areal-uniformity criterion inspired by Kepler's second law, which requires the swept areal velocity observed from the equant to be as uniform as possible. This criterion links the observer point to the intrinsic temporal evolution of the limit cycle and provides a dynamically meaningful measure of oscillatory non-uniformity.
			
	
	As a dimensionless measure, $\hat{K}^\ast$ enables comparison across different limit-cycle systems and quantifies the geometric non-uniformity of their oscillatory motion. It can also provide an indicator of how the limit-cycle geometry changes under parameter variation. In Results~\ref{sec dimension}, we show that, for the FitzHugh--Nagumo system undergoing a supercritical Hopf bifurcation, $\hat{K}^\ast$ decreases towards zero as the bifurcation parameter approaches the tipping point, thereby offering a pre-critical geometric signal of the transition. In Appendix~\ref{appendix proof}, we further prove that $\hat{K}^\ast>0$ for the FitzHugh--Nagumo system under the considered parameter regime, indicating that exact areal uniformity cannot be achieved before the transition.

}
		
	 We proceed to find the equant by extending the learned homeomorphism $\boldsymbol{\chi_{\boldsymbol{\theta}_{\boldsymbol{\chi}}}}$ between the limit cycle and $\rmS^1$ to the invertible map between $\mathbb{M}$ and $\mathbb{D}$. Such an extension $\rvg$ should satisfy 
	 {\color{black}
	\begin{equation*}
		\begin{aligned}
			&\rvg\big|_{\rmC=\partial\mathbb{M}} = \text{Enc},~\rvg^{-1}\big|_{\rmS^1=\partial\mathbb{D}} = \text{Dec},\\
			&\rvg(\rvx^\ast)=\boldsymbol{0},~\rvx^\ast\in\arg\min K(\rvx),\\
			&\rvg\big|_{\mathbb{M}}: r\boldsymbol{\chi}(\phi)+(1-r)\rvx^\ast\mapsto r(\cos\phi,\sin\phi), \\
		\end{aligned}
	\end{equation*}
}
	We employ  a special class of generative model, the invertible NN~\cite{dinh2016density}, to parametrise the extension as $\rvg_{\vtheta_{\rvg}}$. The structure of our NN is comprised of the following three modules.
	\begin{itemize}
		\item [(i)]
		\textbf{Invertible mapping.}~We respectively split the input $\rvx=(x_1,...,x_n)^\top$ and the output $\rvy=(y_1,...,y_n)^\top=\rvf_I(\rvx)$ of the invertible non-linear mapping into two halves, respectively, as $\rvx=(\rvx_{1:d}^\top,\rvx_{d+1:n}^\top)^\top$ and $\rvy=(\rvy_{1:d}^\top,\rmY_{d+1:n}^\top)^\top$ , in which $d<n$ is a hyperparameter. The {\color{black}split} variables obey the following equations,
		\begin{equation*}
			\begin{aligned}
				&\quad\quad\rvy_{1:d} = \rvx_{1:d},\\
				y_{j} &= x_j\exp(s_j(\rvx_{1:d}))+t_j(\rvx_{1:d}),~d+1\le j\le n,
			\end{aligned}
		\end{equation*}
		where $\rvs=(s_{d+1},...,s_n)$ and $\rvt=(t_{d+1},...,t_n)$ are  parametrised NN functions. The corresponding formulation of inverse mapping is trivially calculated as follows,
		\begin{equation*}
			\begin{aligned}
				\rvx_{1:d} &= \rvy_{1:d},\\
				x_{j} &= \exp\left(-s_j(\rvy_{1:d})\right)\left[y_j-t_j(\rvy_{1:d})\right],~d+1\le j\le n.
			\end{aligned}
		\end{equation*}
		\item [(ii)]
		\textbf{Permutation mapping.}~Since the invertible mapping keeps half of the state variables unchanged, we introduce the permutation mapping to transfer the unchanged variables to the next layer with a non-linear activation function. The permutation mapping and its inverse mapping are given by 
		\begin{equation}\label{eq permutation}
			\begin{aligned}
				&\rvf_P((x_1,\cdots,x_n)^\top)=(x_{\alpha(1)},\cdots,x_{\alpha(n)})^\top,~\alpha\in \rmS_n\\
				&\rvf_P^{-1}((x_1,\cdots,x_n)^\top)=(x_{\beta(1)},\cdots,x_{\beta(n)})^\top,~\beta=\alpha^{-1},
			\end{aligned}
		\end{equation}
		where $\rmS_n$ is the permutation group on $\{1,\cdots,n\}$. In implementation, we randomly sample the permutation element and obtain its inverse element via the efficient \texttt{randperm} and \texttt{argsort} methods in Pytorch package~\cite{paszke2019pytorch} as $	\alpha=\texttt{randperm}(n),~\beta=\texttt{argsort}(\alpha)$. 
		
		By alternatively combining the invertible mapping and the permutation mapping, we obtain the parametrised invertible function, 
		\begin{equation}\label{eq INN}
			\rvg_{\boldsymbol{\theta}_{\rvg}}=\rvf_I\circ\rvf_P\cdots\rvf_P\circ\rvf_I.
		\end{equation}
		\item [(iii)]
		\textbf{Dimension alignment.} Although both $\mathbb{M}\subset\mathbb{R}^n$ and $\mathbb{D}\subset\mathbb{R}^2$ are $2$-dimensional manifolds, the above parametrized function cannot be directly employed to find the target bijection since the input and the output of $\rvg_{\boldsymbol{\theta}_{\rvg}}$ belong to $\mathbb{R}^n$. To align the dimension between $\rvg_{\boldsymbol{\theta}_{\rvg}}$ and $\rvg$, we augment the unit disk $\mathbb{D}$ to an $n$-dimensional space as $\mathbb{D}_n=\text{Aug}(\mathbb{D},n-2)$, where $\text{Aug}(\cdot,\cdot)$ is an augmenting operation defined as
		\begin{equation*}	
			\text{Aug}(\mathbb{D},k)=\mathbb{D}\times\{0\}^{k}.	
		\end{equation*}  
		In this way, the parametrised invertible function relates to the target bijection as
		\begin{equation*}
			\rvg_{\boldsymbol{\theta}_{\rvg}}(\rvx)=\text{Aug}(\rmG(\rvx),n-2)\in\mathbb{D}_n,~\text{for all}~\rvx\in\mathbb{M}.
		\end{equation*} 
		
	\end{itemize}
	Finally, we train the NNs according to the following loss function $\mathcal{L}_{\text{equant}}$ as
	\begin{equation*}
		\mathcal{L}_{\text{equant}}(\boldsymbol{\theta}_{\rmG})=\lambda_3\mathcal{L}_3(\vtheta_{\rvg})+\lambda_4\mathcal{L}_4(\vtheta_{\rvg})+\lambda_5\mathcal{L}_5(\vtheta_{\rvg}),
	\end{equation*}
	
	with
	\begin{widetext}
		\begin{equation}\label{eq mc loss2}
			\begin{aligned}
				\mathcal{L}_{3}(\vtheta_{\rvg})&=\dfrac{1}{K}\sum_{i=1}^K\|\rvg_{\boldsymbol{\theta}_{\rvg}}(\rvx_i)-\text{Enc}(\rvx_i)\|_2^2,\\
\mathcal{L}_4(\vtheta_{\rvg})&=	{\color{black}\frac{\sum_{i=1}^{K}\Delta_i\left(A_i-\sum_{k=1}^{K}\Delta_kA_k\right)^2}{\left(\sum_{i=1}^K\Delta_iA_i\right)^2},~A_i=\frac{1}{2}\|\left(\rvg_{\boldsymbol{\theta}_{\rvg}}^{-1}(\boldsymbol{0})-\rvx_i\right)\wedge\rmF(\rvx_i)\|,~\Delta_i=\frac{(t_{i+1}-t_{i})\omega}{2\pi}},\\
				\mathcal{L}_5(\vtheta_{\rvg})&=\dfrac{1}{K}\sum_{i=1}^K\sum_{j=1}^M\bigg\|\rvg_{\boldsymbol{\theta}_{\rvg}}\left(\dfrac{j}{M}\rvg_{\boldsymbol{\theta}_{\rvg}}^{-1}(\boldsymbol{0})+\dfrac{M-j}{M}\rvx_i\right)-\dfrac{M-j}{M}\rvg_{\vtheta_{\rvg}}(\rvx_i)\bigg\|_2^2.
			\end{aligned}
		\end{equation}
	\end{widetext}

		
	{\color{black}Here we require the dataset $\{\rvx_i=\rvx(t_i)\}_{i=1}^K\subset\rmC$ contain temporally ordered limit-cycle samples with known sampling intervals, and $\lambda_{3,4,5}$ to be predefined parameters.} The loss term $\mathcal{L}_3$ guarantees that $\rvg_{\vtheta_{\rvg}}$ is an extension of the homeomorphism between the limit cycle and $\rmS^1$, {\color{black} $\mathcal{L}_4$ is used to find the  equant that has minimal $\hat{K}$ value}, and $\mathcal{L}_5$ measures whether the line of sight $l_{\rvg_{\boldsymbol{\theta}_{\rvg}}^{-1}(\boldsymbol{0})\rvx_i}$ in the original space corresponds to the line of sight $l_{\boldsymbol{0}\rvg_{\vtheta_{\rvg}}(\rvx_i)}$ of universal dynamical clock, i.e., 
	\begin{equation*}
		\begin{aligned}
			\rvg_{\boldsymbol{\theta}_{\rvg}}\left(l_{\rvg_{\boldsymbol{\theta}_{\rmG}}^{-1}(\boldsymbol{0})\rvx_i}\right)&=l_{\boldsymbol{0}\rvg_{\vtheta_{\rvg}}(\rvx_i)},~\rvg_{\boldsymbol{\theta}_{\rvg}}\left(l_{\rvg_{\boldsymbol{\theta}_{\rmG}}^{-1}(\boldsymbol{0})\rvx_i}\right)\\
			&=\left\{\rvg_{\boldsymbol{\theta}_{\rvg}}(\rvx):\rvx\in l_{\rvg_{\boldsymbol{\theta}_{\rvg}}^{-1}(\boldsymbol{0})\rvx_i}\right\}.
		\end{aligned}
	\end{equation*}
	{\color{black}If the time series data is sampled equidistantly over time, we simply set $\Delta_i\equiv 1/K$ in Eq.~\eqref{eq mc loss2}.

We emphasize that the above two steps play distinct roles in our framework and are trained sequentially rather than jointly. In the first step, the autoencoder is trained on the uncoupled limit-cycle dynamics to learn the phase map $\phi:\rmC\to S^1$, its inverse map $\boldsymbol{\chi}:S^1\to\rmC$, and the associated phase sensitivity function $\rmZ(\phi)=\dfrac{\partial\phi}{\partial\rvx}\big|_{\rvx=\boldsymbol{\chi}(\phi)}$. These quantities are sufficient for deriving the reduced phase dynamics under weak perturbations in Sec.~\ref{Methods phase reduction}. After the autoencoder has been trained, its parameters are fixed. In the second step, the invertible network is trained separately to extend the learned boundary map from the limit cycle to a bijection between the surface $\mathbb{M}$ and the unit disk $\mathbb{D}$, such that the corresponding equant achieves the optimal areal non-uniformity. The role of the equant identified by this areal-uniformity criterion is geometric and interpretative: $\hat K^\ast$ captures the areal non-uniformity of the limit cycle and can be used to assess whether the system is approaching a critical transition.

}
	\subsection{Definition of the interpretable phase}\label{subsec interpret}
	
	To comprehend the phase function from a physical view-point with the equant, we consider the following system for example,
	\begin{equation}\label{eq circle}
		\begin{aligned}
			&\dot{x}=y,~\dot{y}=-x,~\dot{z}=-z,\\
			&x(0)=\frac{\sqrt{2}}{2},~y(0)=\frac{\sqrt{2}}{2}~z(0)=z_0.\\
		\end{aligned}
	\end{equation}
	The limit cycle of this system is the unit circle $\rmS^1$ in the $(x,y)$ plane, and the oscillator rotates uniformly about $\rmS^1$. Therefore, any point on $z$-axis is an equant denoted by $\rvx^\ast(z)$, and the viewing angle observed from the origin  $\rvx^\ast(0)$ is the phase.
	Since there exists a linear projection from line $l_{\rvx^\ast(z)\rvx}$ to line {\color{black}$l_{\rvx^\ast(0)\rvx}$} for any $\rvx\in\rmS^1$, we can rescale the viewing angle of any equant to the phase by multiplying by a constant. However, such linear operation does not exist for a general oscillator with irregular limit cycle. In our framework, we address this issue by identifying an invertible map between the surface $\mathbb{M}$ and the domain $\mathbb{D}$. This map acts similarly to a projection, as demonstrated in the example. Finally we come up with the definition of interpretable phase with the equant in our framework as follows.
	\begin{definition}\label{def 1}
		For a phase function $\phi$, if there exists a bijection $\rvg$ between $\mathbb{M}$ and $\mathbb{D}$ such that $\mathrm{(i)}$ $\rvg$ maps the lines of sight emanating from the equant $\rvx^\ast$ to the lines of sight emanating from the centre of $\mathbb{D}$, i.e., $	\rvg(l_{\rvx^\ast\rvx})=l_{\boldsymbol{0}\rvg(\rvx)}$ for any $\rvx\in\mathbb{M}$, and $\mathrm{(ii)}$ for the temporal trajectory $\rvx(t)$ over the limit cycle, the intersection angle between the mapped lines $\rvg(l_{\rvx^\ast\rvx(t)})$ exactly is equal to the phase, i.e., $\phi(\rvx(t))=\angle(\rvg(l_{\rvx^\ast\rvx(t)}),l_{\vzero\rve_1})$ with $\rve_1=(1,0)^\top$. Then we define $\phi$ as the interpretable phase.
	\end{definition}

	\section{Analytical phase reduction under perturbations }\label{Methods phase reduction}
	
	After the training stage, we obtain  $\phi_{\vtheta_\phi}(\rvx)$ and $\boldsymbol{\chi}_{\vtheta_{\boldsymbol\chi}}(\phi)$. Hence, we obtain the phase response curve as 
	\begin{equation*}
		\rmZ(\phi)=\dfrac{\partial\phi_{\boldsymbol{\theta}_{\phi}}}{\partial\rvx}\bigg|_{\rvx=\boldsymbol{\chi}_{\boldsymbol{\theta}_{\boldsymbol{\chi}}}(\phi)},
	\end{equation*} 
	which can be efficiently  computed by \textsc{autograd} method in \textsc{PyTorch} \cite{paszke2017automatic}.
	We utilise these functions to deduce the analytical form of phase response function $\Gamma(\phi,t)$ against three different perturbations. For simplicity, we omit the subscripts $\vtheta_\phi$ and $\vtheta_{\boldsymbol{\chi}}$ in the following derivations. To proceed, we systematically identify the analytical formulation of reduced phase dynamics for general oscillatory systems under external perturbations, including  the stochastic force driven by intrinsic or extrinsic noise~\cite{gonze2002robustness}, the
	periodic force from existing rhythms in nature~\cite{yan2019robust}, and the interplay among oscillators connected in a network~\cite{pecora1998master}.

	\subsection{Stochastic force}\label{subsec stochastic}

	Firstly, we focus on the stochastic oscillation governed by a stochastic differential equation (SDE)
	\begin{equation*}
		\begin{aligned}
			\mathrm{d}{\rvx}&=\rmF(\rvx)\mathrm{d}t+\rmG(\rvx)\mathrm{d}\rmW(t),\\
		\end{aligned}
	\end{equation*}
	where $\rmG(\rvx)\in\mathbb{R}^{n\times m}$ is the diffusion term, and $\rmW(t)\in\mathbb{R}^m$ is the standard Wiener process. By applying I$\hat{\text{t}}$o's formula~\cite{oksendal2013stochastic} and the averaging method, we obtain the phase dynamics on the limit cycle as
		\begin{equation*}
			\begin{aligned}
				\mathrm{d}{\phi}(\rvx) &=\rmZ(\phi)\cdot\mathrm{d}\rvx+\frac{1}{2}\mathrm{d}\rvx^\top\underbrace{\nabla^2\phi(\rvx)\big|_{\rvx=\boldsymbol{\chi}(\phi)}}_{\text{denoted by}~\rmY(\phi)}\mathrm{d}\rvx\\
				&=\rmZ(\phi)\cdot\left[\rmF(\boldsymbol{\chi}(\theta))\mathrm{d}t+\rmG(\boldsymbol{\chi}(\theta))\mathrm{d}\rmW(t)\right]\\
				&\quad+\frac{1}{2}\text{Tr}\left[\rmG(\boldsymbol{\chi}(\phi))^\top \rmY(\phi)\rmG(\boldsymbol{\chi}(\phi))\right]\mathrm{d}t\\
				&=\left\{w+\frac{1}{2}\text{Tr}\left[\rmG(\boldsymbol{\chi}(\phi))^\top \rmY(\phi)\rmG(\boldsymbol{\chi}(\phi))\right]\right\}\mathrm{d}t\\
				&\quad+\boldsymbol\Gamma(\phi)\cdot\mathrm{d}\rmW(t),\\
				\boldsymbol\Gamma(\phi)&=(\rmZ(\phi)\cdot\rmG_1(\boldsymbol{\chi}(\phi)),...,\rmZ(\phi)\cdot\rmG_m(\boldsymbol{\chi}(\phi))^\top,
			\end{aligned}
		\end{equation*}
	here $\text{Tr}[\cdot]$ is the trace operator. 	The corresponding Fokker-Planck equation of the phase distribution $p(\phi,t)$ is,
	\begin{widetext}
		\begin{equation}\label{eq PR noise FK}
			\begin{aligned}
				\partial_t p(\phi,t)&=-\partial_\phi\left\{\left(w+\frac{1}{2}\text{Tr}\left[\rmG(\boldsymbol{\chi}(\phi))^\top \rmY(\phi)\rmG(\boldsymbol{\chi}(\phi))\right]p(\phi,t)\right)\right\} +\frac{1}{2} \sum_{i=1}^m\partial_{\phi}^2\left[\rmZ(\phi)\cdot\rmG_i(\boldsymbol{\chi}(\phi))p(\phi,t)\right].
			\end{aligned}
		\end{equation}
	\end{widetext}
	The influence of noise to oscillation is reflected in the stationary distribution of the above Fokker-Planck equation. 
	
	To further simplify the phase dynamics, we consider the phase decomposition $\phi=wt+\psi$, with $\psi$ representing the effect of noise to the phase. Then we have,
	
	\begin{equation}\label{eq PR noise slow}
		\begin{aligned}
			\mathrm{d}{\phi} &=w+\mathrm{d}{\psi},\\
			\mathrm{d}{\psi}&=\frac{1}{2}\text{Tr}\left[\rmG(\boldsymbol{\chi}(wt+\psi))^\top \rmY(wt+\psi)\rmG(\boldsymbol{\chi}(wt+\psi))\right]\mathrm{d}t\\
			&\quad+\boldsymbol\Gamma(wt+\psi)\cdot\mathrm{d}\rmW(t),\\
		\end{aligned}
	\end{equation}
	When the noise perturbation is weak, $\psi$ is a slow component and approximately stays
	constant over one cycle of oscillation with frequency $w$, the drift term on the right hand side of Eq.~\eqref{eq PR noise slow} changes the $\psi$ over one period $T=\frac{2\pi}{w}$ as
	\begin{equation}\label{eq PR noise delta psi}
		\begin{aligned}
			\delta\psi&=\int_0^{T}\frac{1}{2}\text{Tr}\left[\rmG(\boldsymbol{\chi}(wt))^\top \rmY(wt)\rmG(\boldsymbol{\chi}(wt))\right]\mathrm{d}t\\
			&=\frac{T}{2\pi}\int_0^{2\pi}\frac{1}{2}\text{Tr}\left[\rmG(\boldsymbol{\chi}(\theta))^\top \rmY(\theta)\rmG(\boldsymbol{\chi}(\theta))\right]\mathrm{d}\theta
		\end{aligned}
	\end{equation}
	Then the dynamics of $\delta\psi$ during one period is equivalent to a constant,
	\begin{equation}\label{eq PR noise delta w}
		\begin{aligned}
			\delta w&=\frac{\delta\psi}{T}\\
			&=\frac{1}{2\pi}\int_0^{2\pi}\frac{1}{2}\text{Tr}\left[\rmG(\boldsymbol{\chi}(\theta))^\top \rmY(\theta)\rmG(\boldsymbol{\chi}(\theta))\right]\mathrm{d}\theta.
		\end{aligned}
	\end{equation}
	Putting Eqs.~\eqref{eq PR noise delta psi},\eqref{eq PR noise delta w} into Eq.~\eqref{eq PR noise slow} leads to the following dynamics,
	\begin{equation*}
		\begin{aligned}
			\mathrm{d}{\psi} &= \delta w\mathrm{d}t+\boldsymbol\Gamma(wt+\psi)\cdot\mathrm{d}\rmW(t).
		\end{aligned}
	\end{equation*}
	Similarly, we perform the time-average operation to the diffusion term and obtain the simplified phase dynamics as,
	\begin{equation*}
		\begin{aligned}
			\mathrm{d}{\psi} &= \delta w\mathrm{d}t+\sum_{i=1}^m\Gamma_i\mathrm{d}W_i(t),\\
			\Gamma_i&=\frac{1}{2\pi}\int_0^{2\pi}\rmZ(\theta)\cdot\rmG_i(\boldsymbol{\chi}(\theta))\mathrm{d}\theta.
		\end{aligned}
	\end{equation*}
	Therefore, we obtain the simplified phase dynamics as,
	\begin{equation}\label{eq SDE phase}
		\begin{aligned}
			\mathrm{d}{\phi} = \left(w+\delta w\right)\mathrm{d}t+\sum_{i=1}^m\Gamma_i\mathrm{d}W_i(t).
		\end{aligned}
	\end{equation}
	%
	Although we use the standard Wiener process in our derivation, the results presented above can be easily extended to the more general case involving correlated Wiener processes.
	
	To proceed, we investigate the synchronisation behaviour of collective dynamics driven by scalar common noise and independent noise as,
	\begin{equation}\label{eq common+independ}
		\mathrm{d}{\rvx}_i=\rmF(\rvx)\mathrm{d}t+\rmG(\rvx)\mathrm{d}\xi(t)+\rmH(\rvx)\mathrm{d}\eta_i(t),~i=1,\cdots,N.
	\end{equation}
	We denote the corresponding phase dynamics as,
	\begin{equation*}
		\begin{aligned}
			\mathrm{d}{\phi}_i &= \left(w+\delta w\right)\mathrm{d}t+ \Gamma_1(\phi_i)\mathrm{d}\xi+\Gamma_2(\phi_i)\mathrm{d}\eta_i,\\
			\delta w &= \dfrac{1}{2\pi}\int_0^{2\pi}\frac{1}{2}\text{Tr}\left[(\rmG+\rmH)^\top \rmY(\theta)(\rmG+\rmH)\right]\mathrm{d}\theta.\\
		\end{aligned}
	\end{equation*}
	Here $\Gamma_1(\phi)=\rmZ(\phi)\cdot \rmG(\boldsymbol{\chi}(\phi))$, $\Gamma_2(\phi)=\rmZ(\phi)\cdot \rmH(\boldsymbol{\chi}(\phi))$, and the Wiener processes $\xi$ and $\{\eta_i\}_{i=1}^N$ are common noise and independent noise, respectively.
	The synchronisation behaviour of the noise driven dynamics is reflected in the phase difference  {\color{black}$\delta\phi_i=\phi_i-\phi_{i+1},~1\le i\le N-1$}. We follow the approach in~\cite{nakao2007noise} to deduce the stationary distribution of the phase difference as
	{\color{black}
		\begin{equation}\label{eq FK common+independ}
		P^\ast(\delta\phi)=\dfrac{p_0}{\gamma_1(0)-\gamma_1(\delta\phi)+\gamma_2(0)},
	\end{equation}
	where $p_0$ is the normalisation constant, $\gamma_i(\psi), i=1,2$ are time-average functions defined as
	\begin{equation*}
	\gamma_i(\psi)=\frac{1}{2\pi}\int_0^{2\pi}\Gamma_i(\theta)\Gamma_i(\psi+\theta)\mathrm{d}\theta,~i=1,2.
	\end{equation*}
	Therefore, the long term collective behaviour of the noise perturbed oscillators can be readily determined by calculating $P^\ast$ using our framework. For completeness, we provide a full derivation of Eq.~\eqref{eq FK common+independ} in Appendix~\ref{appendix deri}.
}

	
	\subsection{Periodic force}\label{subsec period}
	We are concerned whether an oscillator driven by the external periodic force
	$\rvf(t)$ with the frequency $\Omega$ shows phase lock phenomenon, i.e., if the frequency of driven
	oscillator converges to $\Omega$. {\color{black}We require the periodic forcing to be weak so that the perturbed dynamics can be  typically analysed in the vicinity of a small neighbourhood of the original limit cycle.} The perturbed dynamics is described as
	\begin{equation*}
		\dot{\rvx}=\rmF(\rvx)+ \rvf(t),~\rvf(t)=\rvf(t+T).
	\end{equation*}
	The corresponding phase dynamics is,
	\begin{equation*}
		\begin{aligned}
			\dot{\phi}&=\rmZ(\phi)\cdot\left[\rmF(\rvx)+\rvf(t)\right]\\
			&=w+\rmZ(\phi)\cdot\rvf(t).
		\end{aligned}
	\end{equation*}
	To investigate the influence of the periodic force to the phase frequency, we consider the phase difference $\delta \phi=\phi-\Omega t$ whose dynamics is governed by, 
	\begin{equation}\label{eq phase period}
		\begin{aligned}
			\dot{\delta \phi}&=\dot{\phi}-\Omega\\
			&=\Delta+\rmZ(\Omega t+\delta\phi)\cdot\rvf(t),\\
			\Delta&=w-\Omega.~
		\end{aligned}
	\end{equation}
	By time-averaging the weak term over one cycle of oscillation of frequency $\Omega$, we obtain the closed-form dynamics of phase difference as,
	\begin{equation*}
		\begin{aligned}
			\dot{\delta \phi}&=\Delta+\Gamma(\delta \phi),\\
			\Gamma(\delta \phi) &= \dfrac{1}{2\pi}\int_0^{2\pi}\rmZ(\delta \phi+\theta)\rvf\left(\frac{\theta}{\Omega}\right)\mathrm{d}\theta.
		\end{aligned}
	\end{equation*} 
	The phase lock occurs once the vector field $\Delta+\Gamma(\delta\phi)$ has a stable zero solution. Therefore, by viewing the intersection of constant $\Delta$ and function $\Gamma(\delta \phi)$ in the diagram, our framework efficiently discerns whether the periodic force leads to phase lock or divergence. 
	
	\subsection{Network interaction}\label{subsec network}
	We consider the collective dynamics of coupled oscillators, expressed in a general form as
	\begin{equation}
		\begin{aligned}	\dot{\rvx}_i&=\rmF_i(\rvx_i)+\sum_{j=1}^NA_{ij}\rmH(\rvx_i,\rvx_j),\\
		\end{aligned}
	\end{equation}
	where $\rmA = (A_{ij})\in\mathbb{R}^{N\times N}$ captures the interacting structure among the oscillators, and  $\rmH(\rvx_i,\rvx_j)$ represents the influence of $j$th oscillator to $i$th oscillator. {\color{black}We first consider the identical oscillators with $\rmF_i\equiv\rmF$.} Then we have the phase dynamics in the vicinity of $\rmC$ as
	\begin{equation}\label{eq PR network}
		\begin{aligned}
			\dot{\phi}_i &= \rmZ(\phi_i)\cdot\rmF(\rvx_i)+\sum_{j=1}^NA_{ij}\rmZ(\phi_i)\cdot\rmH(\rvx_i,\rvx_j)\\
			&=w+\sum_{j=1}^NA_{ij}\rmZ(\phi_i)\cdot\rmH(\boldsymbol{\chi}(\phi_i),\boldsymbol{\chi}(\phi_j))\\
		\end{aligned}
	\end{equation}
	from which the phase response function of network interaction is identified by a bivariate function $\tilde{\Gamma}(\phi_i,\phi_j)=\rmZ(\phi_i)\cdot\rmH({\color{black}\boldsymbol{\chi}(\phi_i)},\boldsymbol{\chi}(\phi_j))$. Even though the form $\tilde{\Gamma}(\phi_i,\phi_j)$ is clear and precise, a more predominant way is to simplify the phase response function as a univariate function $\Gamma(\phi_j-\phi_i)$ by averaging $\Gamma$ over one period~\cite{kuramoto1984chemical,nakao2016phase}. In this way, the phase response only  relies on the phase difference and the reduced phase dynamics is described as,
	\begin{equation}\label{eq PR network final}
		\begin{aligned}
			\dot{\phi}_i &=w+
			\sum_{j=1}^{N}A_{ij}\Gamma(\phi_j-\phi_i),\\
			\Gamma(\phi_j-\phi_i)&=\dfrac{1}{2\pi}\int_0^{2\pi}\rmZ(\theta+\phi_i)\rmH\left(\boldsymbol{\chi}(\theta+\phi_i),\boldsymbol{\chi}(\theta+\phi_j)\right)\mathrm{d}\theta.
		\end{aligned}
	\end{equation}
	The proof is mainly based on the method in~\cite{kuramoto2019concept}, {\color{black}we provide the detailed derivation in Appendix~\ref{appendix deri} for completeness.}

	Next, we extend our framework to heterogeneous oscillators as follows,
	\begin{equation*}
		\begin{aligned}				 
			\rmF_i(\rvx_i)=\rmF(\rvx_i)+\delta\rmF_i(\rvx_i),~\dfrac{\|\delta\rmF_i(\rvx_i)\|}{\|\rmF_i(\rvx_i)\|}\ll1,
		\end{aligned}
	\end{equation*}
	{\color{black}here, $\delta \rmF_i(\rvx_i)$ denotes the deviation from the original dynamics caused by heterogeneous perturbations; for example, when a parameter $\mu$ is perturbed to $\mu_i$ for the $i$-th oscillator, then we have $\delta \rmF_i(\rvx_i)=\rmF_i(\rvx_i,\mu_i)-\rmF(\rvx_i,\mu)$. The corresponding} phase dynamics is governed by
	\begin{equation}\label{eq PR network heter final}
		\begin{aligned}
			\dot{\phi}_i&=(w+\delta w_i)+\sum_{i=1}^{N}A_{ij}\Gamma(\phi_j-\phi_i),\\
			\delta w_i&=\rmZ(\phi_i)\cdot\delta\rmF(\rvx_i)=\dfrac{1}{2\pi}\int_0^{2\pi}\rmZ(\theta)\cdot\delta\rmF_i(\boldsymbol{\chi}(\theta))\mathrm{d}\theta.
		\end{aligned}
	\end{equation} 
	{\color{black}The detailed proof is also provided in Appendix~\ref{appendix deri}.}

	Although we derive the phase reduction for each of the three perturbations separately, the phase dynamics resulting from the combination of these perturbations can  easily be obtained by summing the corresponding phase response functions.
	
	For the analytical expression of the phase dynamics, we note that the NN functions are analytical when the parameters are fixed after training. To  enhance the interpretability of the established phase response functions, we approximate the NN-based function by a linear combination of known basis functions as
	\begin{equation*}
		\Gamma(\phi)=\sum_{i=1}^{M}a_if_i(\phi),
	\end{equation*}
	where $\{f_i:\mathbb{R}\to\mathbb{R}\}_{i=1}^{M}$ is the predefined dictionary of basis functions. Since the phase is a periodic function, the basis function can be further required as $f_i:\mathbb{R}/_\sim\to\mathbb{R}$, e.g., the Trigonometric basis.\\

	\section{Results and Applications}\label{sec results}
	The previous sections primarily establish the universal dynamical clock and the associated phase dynamics. In this section, we highlight the merits of this principle in understanding and regulating rhythmic oscillations from distinct perspectives: (i) resolving a long-standing open problem on synchronisation in \textit{Escherichia coli}~\cite{garcia2004modeling}; (ii) enabling precise modulation of synthetic genetic circuits based on the experimental data; (iii) {\color{black}establishing a classic-mechanics counterpart of the Berry geometric phase;} (iv) revealing the role of spatial dimensionality in human perception of rhythms; and (v) uncovering phase synchronisation mechanisms under complex perturbations. By rendering phase an explicit and operational variable, the dynamical clock enables the systematic classification of oscillatory regimes and provides a principled basis for analysing, comparing, and modulating network dynamics underlying distinct functional and computational behaviours.

	\begin{figure*}[htp]
		\centering
		\includegraphics[width=\textwidth]{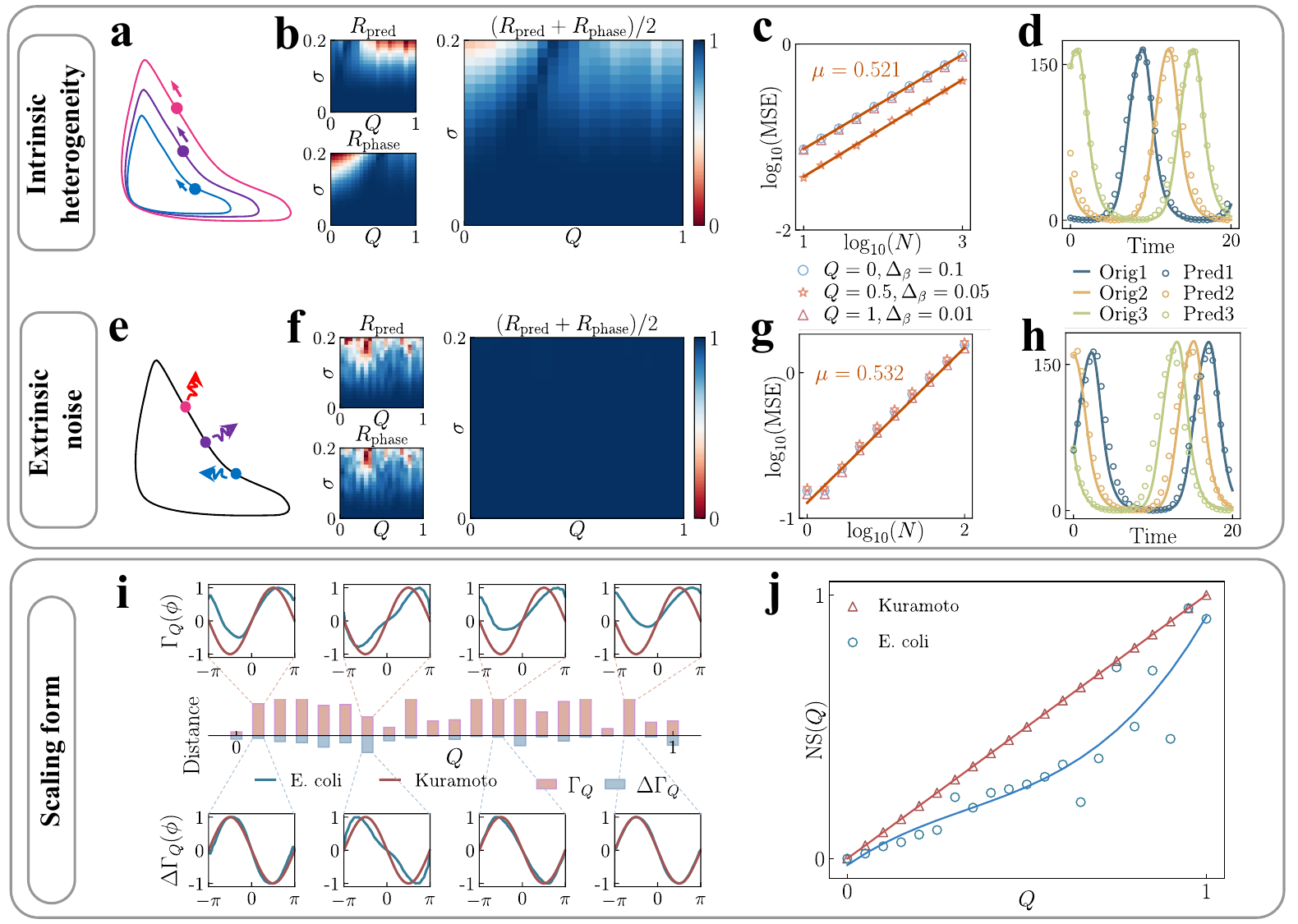}
		\caption{\textbf{Scaling form of E. coli cell populations under quorum sensing}. (\textbf{A}) Illustration of the influence of intrinsic {\color{black}heterogeneity} to the limit cycle of E. coli cell.   (\textbf{B}) Similarity between different synchronisation indices with the original synchronisation index $R_{\text{orig}}$ over different cell densities $Q$ and noise strengths $\sigma$.
			(\textbf{C}) Mean square error (MSE) between 
			the original trajectories and the predicted trajectories under different cell numbers $N$.
			(\textbf{D}) Trajectories of three randomly selected cells. (\textbf{E}-\textbf{H}) show the similar results under extrinsic noise. (\textbf{I}) Response functions of phase and phase difference  under different cell densities $Q$. The bar diagram shows the distance between the response functions with the Kuramoto model. (\textbf{J}) Dynamics of phase difference between any two E. coli cells under different cell densities $Q$. (\textbf{K}) Natural strength as a function of the cell density $Q$, solid lines represent the fitting results to the data points.}
		\label{fig3}
	\end{figure*}

	\subsection{Scaling form of Escherichia coli population}\label{sec Ecoli}
	{\color{black}
		The dynamical clock derived from the equant reveals a unifying geometric structure underlying spontaneous rhythms in nature. Beyond offering a new perspective, it enables the resolution of scientific problems governed by rhythmic dynamics. We first address an open problem posed by Strogatz in 2004: does synchronisation in Escherichia coli populations under quorum sensing follow Kuramoto-like scaling~\cite{garcia2004modeling}? Quorum sensing, a global mean-field coupling mechanism facilitating cell-to-cell communication, is known to induce collective genetic oscillations~\cite{mcmillen2002synchronizing}. While its influence on collective behaviour has been quantified quantitatively~\cite{taylor2009dynamical}, its effect on phase dynamics remains unresolved.

		The main challenge in identifying the reduced phase dynamics of \textit{E. coli} population under quorum sensing arises from the persistent mean-field coupling at synchrony, which causes each cell’s dynamics to deviate from the dynamical clock observed from the equant in the uncoupled case.} To address this issue, we reorganise the original mathematical model such that the coupling term vanishes when the cellular states are equal to the mean field by incorporating a portion of the coupled terms into the self-dynamics of each cell (see \textbf{Appendix~\ref{appendix Ecoli}} for a detailed derivation). We notice that the reorganised self-dynamics varies with the coupling strength, which represents the cell density, and we identify the specialised equant with its corresponding analytical reduced phase dynamics for each cell density separately. 
	
	To test the accuracy of the established phase dynamics, we focus on two major perturbations affecting population diversity: (i) intrinsic {\color{black}heterogeneity}, which contributes to heterogeneity in oscillation frequency and amplitude across cells (Fig.~\ref{fig3}A), and (ii) stochastic extrinsic noise, which continuously perturbs cellular dynamics (Fig.~\ref{fig3}E).
	Specifically, we collect state trajectories, predicted state trajectories, and phase trajectories under varying noise intensities and cell densities. {\color{black}For measuring the statistics of the trajectories, we} compute the synchronization degree of the original trajectories ($R_{\text{orig}}$), the reconstructed trajectories ($R_{\text{pred}}$), and the phase trajectories ($R_{\text{phase}}$).
	 {\color{black}Here, $R_{\text{orig}}$ and $R_{\text{pred}}$ are computed from the original and reconstructed observable trajectories, respectively, using the observable-based variance-ratio synchrony measure proposed in~\cite{garcia2004modeling},
	 \begin{equation*}
	 	\begin{aligned}
	 		R&=\dfrac{\langle M^2\rangle-\langle M\rangle^2}
	 		{\frac{1}{N}\sum_{i=1}^N\left(\langle b_i^2\rangle-\langle b_i\rangle^2\right)},\\
	 		M&=\dfrac{1}{N}\sum_{i=1}^N b_i(t),
	 	\end{aligned}
	 \end{equation*}
	 where $\langle\cdot\rangle$ denotes the time average. By contrast, $R_{\text{phase}}$ is the Kuramoto order parameter computed directly from the learned phase trajectories.}
	By comparing $R_{\text{orig}}$ and $R_{\text{pred}}$ in Figs.~ \ref{fig3}B \& F, we find that the identified phase dynamics remain accurate under most noise intensity and cell density combinations, except in regions of strong noise, where perturbed cellular rhythms deviate from the original limit cycle. A similar trend holds for $R_{\text{orig}}$ and $R_{\text{phase}}$. Moreover, we observe that the average of $R_{\text{pred}}$ and $R_{\text{phase}}$ aligns more closely with $R_{\text{orig}}$.
	Figures~\ref{fig3}C \& \ref{fig3}G show that the prediction error between the original and predicted trajectories follows a power-law scaling with respect to the number of cells, i.e.,~$\text{MSE}\approx N^\mu$. Additionally, Figures~\ref{fig3}D \& \ref{fig3}H further verify the accuracy of the established phase dynamics, as the real and predicted trajectories remain perfectly consistent over one complete oscillation period.

	Given the accurate reduced phase dynamics, we can determine the scaling form of the \textit{E. coli} population under quorum sensing by analysing the phase response function $\Gamma_Q(\phi)$ across different coupling strengths $Q$.
	While direct plots of $\Gamma_Q(\phi)$ for various $Q$ values show significant deviations from a sinusoidal function, the phase difference dynamics between any two cells resemble those of the Kuramoto model. Specifically, the phase difference $\delta\phi$ evolves according to $\Delta_Q\Gamma(\delta\phi)=\Gamma_Q(-\delta\phi)-\Gamma_Q(\delta\phi)$, whose graph closely resembles the $-\sin(\delta\phi)$ function in the Kuramoto model. Figure~\ref{fig3}I provides statistical comparisons between the phase response functions and the Kuramoto model, further supporting the observation that phase differences evolve according to Kuramoto-type dynamics. To further examine this relationship, we define the natural strength $\text{NS}(Q)$ as the amplitude of the periodic function $\Delta\Gamma_Q$. Thus, the scaling law of $\text{NS}(Q)$ with respect to $Q$ quantifies the influence of coupling on phase dynamics. For reference, in the Kuramoto model, $\text{NS}(Q)$ scales linearly. As shown in Fig.~\ref{fig3}J, $\text{NS}(Q)$ approximates the linear function $y=x$ at the endpoints of $[0,1]$ but remains below it within the interval, exhibiting second-order scaling. In conclusion, {\color{black}our results provide a negative solution to the open problem:} the quorum sensing dynamics of \textit{E. coli} differ from the standard Kuramoto model, but the phase difference dynamics exhibit a scaling behaviour similar to the Kuramoto model, with a super-linear increase in phase response intensity as the coupling strength increases.
	
	{\color{black}We further demonstrate the predictive power of the dynamical clock in data-driven settings and systematically assess the effectiveness of the unifying geometric structure across diverse real-world oscillatory systems, as shown in the following results.}\\

	\begin{figure*}[htp]
		\centering
		\includegraphics[width=\textwidth]{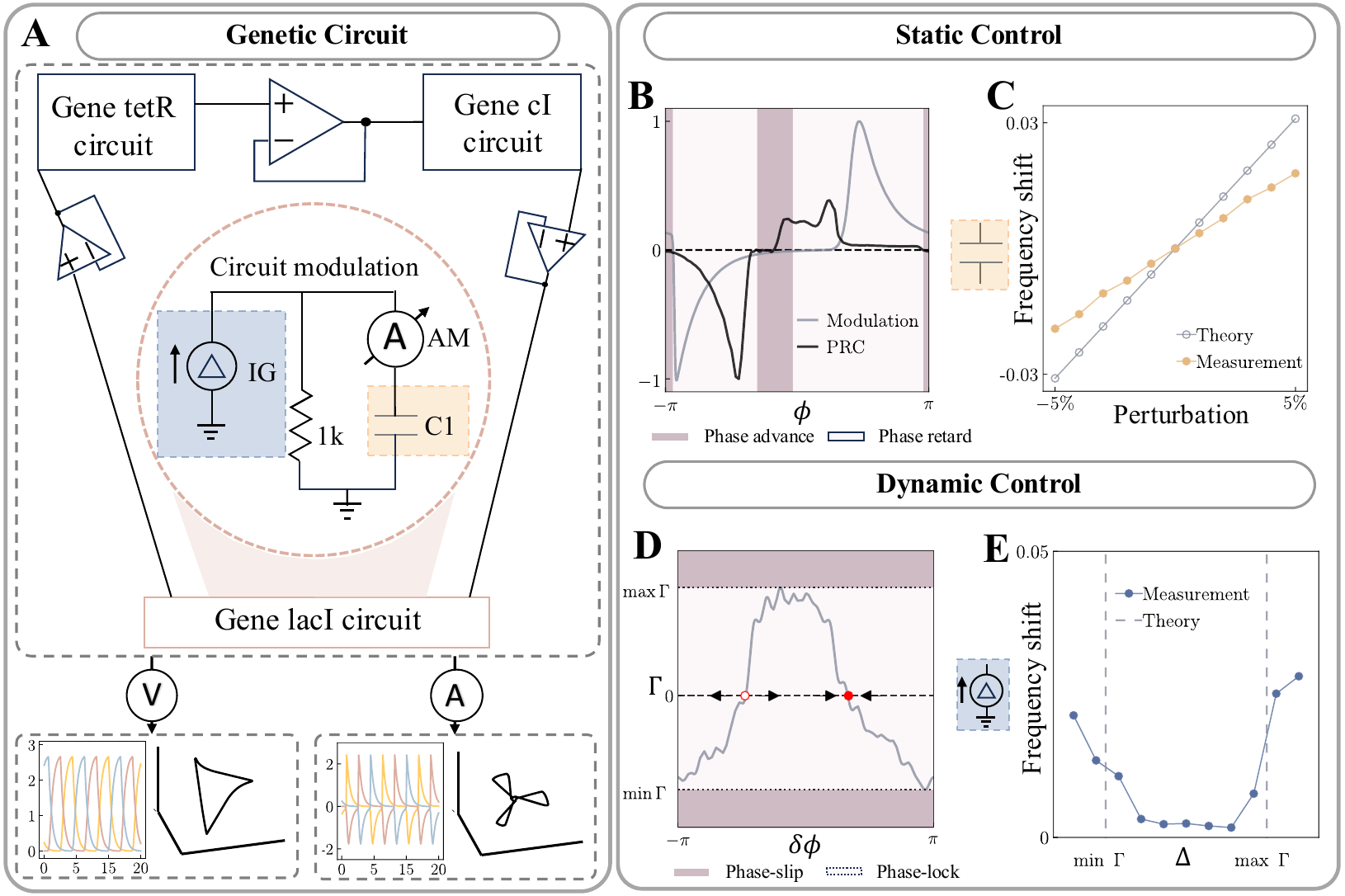}
		\caption{\textbf{Modulation of genetic circuit: A data-driven realization}. (\textbf{A}) Illustration of modulated genetic circuits and the data from the circuits including voltages and currents. (\textbf{B}) The PRC and modulation together advance the phase when they have the same signs, otherwise retard the phase. (\textbf{C}) The numerical validation of the theoretical prediction. (\textbf{D})  The phase response functions against periodic force. (\textbf{E}) The frequency difference between the modulated genetic circuit and the external force versus $\Delta=w-\Omega$, the difference between the natural frequency and external frequency.  }
		\label{fig_genetic}
	\end{figure*}

	\subsection{Modulating genetic circuit: A data-driven realization}\label{sec genetic}       
	Cells naturally navigate environments, communicate, and construct complex patterns by initiating gene expression in response to specific signals, enabling the programming of cells to perform sophisticated tasks or to synthesize complex chemicals and materials. {\color{black}We illustrate the power of the dynamical clock, grounded in the equant, to modulate genetic circuits by selecting appropriate regulators and tuning expression levels, thereby enabling advanced cellular programming.} Specifically, we apply our approach to the Repressilator, a three-gene inhibitory ring oscillator previously implemented in E. coli~\cite{hellen2017electronic}.  {\color{black}Instead of using the data produced by simulating the model of differential equations, we use an electronic circuit to mimic the genetic oscillator and to produce the data.  The dynamical process of three RNA concentrations and the feedback coupling effect among the three components are modelled by the dynamic resistor-capacitor  circuits and the repressor circuits, respectively.} The detailed experimental configurations are provided in Appendix~\ref{appendix gene}.
	
	Notably, we focus on two critical aspects of genetic circuit design, as emphasized in~\cite{brophy2014principles,holtz2010engineering}: static control, which involves tuning expression parameters to balance circuit behaviour, and dynamic control, which underscores adaptation to external environmental changes, such as peripheral clocks.  In the Repressilator circuit, the protein concentration and the rate of gene expression are, respectively, represented by the voltage and the current recordings experimentally collected from the real circuit, as illustrated in Fig.~\ref{fig_genetic}A. We analyse these experimental recordings to extract the phase response curve (PRC) of the dynamical clock, as observed from the machine-learning-coined equant.
	
	{\color{black}For static control, we modulate the transcription rate of the lacI gene, corresponding to tuning the capacitance in the circuit.} The resulted phase shifts in RNA concentration cycle are determined by the alignment between the PRC and the modulation function derived from our framework: when the two share the same sign, the phase advances; when they differ, it delays (see Fig.~\ref{fig_genetic}B). This modulation strategy yields predicted oscillation frequencies that closely match experimental measurements across a broad spectrum of perturbation amplitudes, as shown in Fig.~\ref{fig_genetic}C.

	To achieve dynamic control, we apply a periodic sinusoidal current via an external source to the RC circuit of the lacI oscillator, entraining or perturbing repressor dynamics in analogy to synthetic approaches for engineering biological clocks. We investigate phase-lock phenomena by conducting bifurcation analysis based on the phase response function. Figure~\ref{fig_genetic}D delineates the parameter space of frequency difference $\Delta$, identifying a phase-lock region (light area) where the circuit frequency converges to that of the periodic forcing. The results shown in Fig.~\ref{fig_genetic}E validate the accuracy of the predicted boundary between phase-lock and phase-slip regions.

		{\color{black}
		\subsection{Dynamical Clock and Berry Geometrical Phase}\label{appendix berry}
		In quantum mechanics, the Berry geometrical phase arises in a quantum system when it undergoes cyclic adiabatic evolution~\cite{berry1984quantal}. The quantum adiabatic theorem~\cite{simon1983holonomy} applies to such a case when the system's Hamiltonian $\rmH(\rmR(t))$ varies periodically with time $t$, where $\rmR(t)$ is the Hamiltonian parameter. When the parameter varies along a closed loop $\rmC_\rmR$ with period $T_\rmR$, the state of the quantum system over one cycle is described by
		\begin{equation}\label{eq berry}
			\resizebox{\linewidth}{!}{$ |\Psi(T_\rmR)\rangle=\exp\left[\text{i}\gamma_n(\rmC_\rmR)\right]\exp\left\{\frac{-\text{i}}{\hbar}\int_0^{T_\rmR}E_n(\rmR(t))\mathrm{d}t\right\}|\Psi(0)\rangle.
				$}
		\end{equation}
		Here,  the first exponential term is the geometric term, with $\gamma_n(\rmC_\rmR)$ being the Berry geometrical phase, $E_n$ is the $n$th eigenvalue of the Hamiltonian, and the second exponential term is the ``dynamic phase factor". The Berry phase is often analysed  in quantum mechanics, and plays a significant role in phenomena such as the Aharonov-Bohm effect~\cite{aharonov1959significance}. 
		
		Akin to the geometric phase in quantum systems, previous studies have extended the concept to classical dynamics, encompassing the Pancharatnam phase in classical optics~\cite{pancharatnam1956generalized}, Hannay's angle in integrable Hamiltonian systems~\cite{hannay1985angle}, and holonomy in the Foucault pendulum~\cite{wilczek1989geometric}. The mathematical formulation of the geometric phase for dissipative systems with cyclic attractors was established under adiabatic conditions by~\cite{kepler1991geometric}.
		
		However, existing approaches lack a unified framework for computing the geometric phase in general limit cycle systems with adiabatic parameter evolution. In our work, we address this gap by constructing a homeomorphism between the limit cycle and the unit circle, thereby capturing the impact of cyclic parameter variation on phase dynamics through our phase reduction method. This allows us to present a unified method for determining the geometric phase in general limit cycle systems.
		
		To formalize this approach, we first define the geometric phase as follows.
		
		\noindent\textbf{Definition 3~(Geometrical Phase for Classical Dynamical System)} \textit{
			We consider the limit cycle dynamical system as
			$\dot{\rvx}=\rmF(\rvx,\rmR(t))$, with $\rmR(t)\equiv0$ degenerating the system to the Eq.~\eqref{eq1}. We denote by $\rmC$ the limit cycle of the system $\dot{\rvx}=\rmF(\rvx,\vzero)$, and $\mathcal{T}:~\rmC\to \rmS^1$ the homeomorphism between the limit cycle and the unit circle. When the parameter $\rmR( t)$ varies along a closed loop $\rmC_\rmR$ with period $T_\rmR\to\infty$ (adiabatic condition), the state over one period of parameter is described by}
		\begin{equation}\label{eq berry classic}
			\resizebox{\linewidth}{!}{$ 
				\mathcal{T}(\rvx(T_\rmR)) = \exp(\mathrm{i}\gamma(\rmC_\rmR))\exp\left(\mathrm{i}\int_0^{T_\rmR}w(\rmR(t))\mathrm{d}t\right)\mathcal{T}(\rvx(0)).
				$}
		\end{equation}
		\textit{Here $w(\rmR(t))$ is the instantaneous frequency of the dynamics $\rmF(\rvx,\rmR(t))$ along the limit cycle,  and $\gamma(\rmC_\rmR)$ is the geometrical phase of the limit cycle system. }

		\noindent\textbf{Remark 1}~\textit{
			In the above Definition, we are considering the dynamics with irregular limit cycle. For an explanation to the Eq.~\eqref{eq berry classic}, we take the simple case with $\rmC=\rmS^1$ for an example. In this case, the initial state $\rvx(0)$ and the final state $\rvx(T_\rmR)$ are on the unit circle. Thus, the two exponent terms in Eq.~\eqref{eq berry classic} contribute to the accumulative angle variation during one period $T_\rmR$, that is, 
			\begin{equation}
				\Delta \theta_{T_\rmR}=\gamma+\int_0^{T_\rmR}w(\rmR(t))\mathrm{d}t.
			\end{equation}
			Here,~$\Delta \theta_{T_\rmR}$ represents the total variation quantity of angle between $\rvx(0)$ and $\rvx(T_\rmR)$ during one period $T_\rmR$.}
		
		We provide a unifying perspective to better understand the Eqs.~\eqref{eq berry} and \eqref{eq berry classic} in the following remark.
		
		\noindent\textbf{Remark 2}~\textit{
			In quantum mechanics, the  dynamical phase factor $\exp\left\{\frac{-\text{i}}{\hbar}\int_0^{T_\rmR}E_n(\rmR(t))\mathrm{d}t\right\}$ in Eq.~\eqref{eq berry} represents the accumulation of energy over time. Here, let $\rvn(\rmR)$ be the $n$th eigenstate of the Hamiltonian $\rmH(\rmR)$. The corresponding eigenvalue equation is: $\frac{-i}{\hbar}\rmH(\rmR)\rvn(\rmR)=\frac{-i}{\hbar}E_n(\rmR)\rvn(\rmR)$. This implies that energy serves as a frequency along the eigenstate direction $\rvn(\rmR)$ in the Schrödinger equation. In our classical dynamical system, an analogous term $\exp\left(\mathrm{i}\int_0^{T_\rmR}w(\rmR(t))\mathrm{d}t\right)$ characterizes the effect of frequency along the limit cycle. From the perspective of Koopman operator theory~\cite{brunton2016discovering}, lifting the finite-dimensional limit cycle dynamics $\dot{\rvx}=\rmF(\rvx,\rmR(t))$ into an infinite-dimensional function space results in a dissipative operator $\mathcal{K}(\rmR(t))$. The lifted dynamics is given by: $\dot{g}(\rvx)=\mathcal{K}(\rmR(t))g(\rvx)$ for any function $g(\rvx)$. Since the system possesses a limit cycle, the Koopman operator $\mathcal{K}$ admits a unitary eigenpair $(\lambda_{\pm}(\rmR(t)),g_{\pm}(\rmR(t)))$. Here, $g_{\pm}(\rmR(t))$ corresponds to a function defined on the limit cycle, while the eigenvalue $\lambda_{\pm}(\rmR(t))$ is related to the phase frequency $w(\rmR(t))$ as $\lambda_{\pm}=\cos(w)\pm\mathrm{i}\sin(w)$. Thus, we establish a unified perspective on geometric phase in both quantum and classical mechanics by considering the eigendecomposition of their respective infinite-dimensional operators.}

		To proceed, we show how to calculate the Berry geometrical phase using our framework.  In this paper, we consider a simple case wherein the cyclic parameter is additive, that is, $\rmF(\rvx,\rmR(t))=\rmF(\rvx)+\rmR(t)$. Here, we regard the parameter as perturbation. This kind of perturbation is common in the neuronal dynamics such as the injected currents to the neuron. To maintain the limit cycle, we require the parameter $\rmR$ to be in the vicinity of zero. Then the evolution of the perturbed dynamics is approximated by the corresponding perturbed phase dynamics over the limit cycle as
		\begin{equation}\label{eq berry 1}
			\dot{\phi}=w+\rmZ(\phi)\cdot\rmR(t).
		\end{equation}
		The phase shift between Eq.~\eqref{eq berry 1} and the unperturbed phase dynamics over one cycle of the perturbation is exactly the Berry geometrical phase as defined in~\cite{kepler1991geometric}, 
		\begin{equation}\label{eq berry 2}
			\gamma(\rmR)=\int_0^{T_\rmR}\rmZ(\phi(t))\cdot\rmR(t)\mathrm{d}t.
		\end{equation}
		The integral term in Eq.~\eqref{eq berry 2} is a functional depending on the perturbation $\rmR(t)$ and the phase trajectory $\phi(t)$. Nevertheless, in the adiabatic condition, the period $T_\rmR$ tends to infinity and largely {\color{black}exceeds} the phase period $2\pi/w$. The resultant integral eliminates the local effect and hence becomes an global invariant, as observed in~\cite{kepler1991geometric}.
		Notice that Eq.~\eqref{eq berry 2} is the Berry phase for systems whose cyclic parameter is additive. For general systems, where $\rmF(\rvx,\rmR)$ is not additive, our proposed Eq.~\eqref{eq berry classic} is the extension of the Berry phase for classical systems.
	
	}
	{\color{black}
To numerically validate the above formulation, we further consider the FitzHugh--Nagumo system under weak cyclic parameter perturbations. Specifically, we introduce periodic perturbations into the recovery-variable vector field
$
\dot{y}=x+a-by
$ by varying either the additive parameter \(a\) or the multiplicative parameter \(b\) along a closed cycle as, 
$
a(t)=a+\sigma\sin(\pi t/T_R)$,~$
b(t)=b+\sigma\sin(\pi t/T_R),$ 
where \(T_R\) is the perturbation period and \(\sigma\) is the perturbation strength. For the perturbation of \(a\), the induced phase shift is directly given by Eq.~\eqref{eq berry 2}. For the perturbation of \(b\), the perturbation enters the vector field through the term \(-by\), and hence the corresponding phase response is proportional to \(Z_y(\phi)\chi_y(\phi)\), where \(Z_y\) denotes the \(y\)-component of the phase response curve and \(\chi_y\) denotes the \(y\)-component of the inverse phase map. Therefore, the accumulated geometric phase over one perturbation cycle is predicted by
\begin{equation}\label{eq geo b}
	\gamma_b(R)
	=
	-\sigma\int_0^{T_R}
	Z_y(\phi(t))\chi_y(\phi(t))
	\sin(\pi t/T_R)\,\mathrm{d}t .
\end{equation}
For each perturbation cycle, we compute the geometric phase in two independent ways. First, the theoretical geometric phase \(\gamma_{\mathrm{theory}}\) is obtained from the reduced phase dynamics derived from the Eqs.~\eqref{eq berry 2},\eqref{eq geo b}. Second, the numerical geometric phase \(\gamma_{\mathrm{num}}\) is measured directly from simulations of the full perturbed FitzHugh--Nagumo system by applying the learned phase map to the state after one perturbation period. As shown in Fig.~\ref{fig_geo}, the theoretical predictions closely match the geometric phases obtained from the full numerical simulations across different perturbation periods \(T_R\) and perturbation strengths \(\sigma\). This agreement is observed for cyclic perturbations of both \(a\) and \(b\), indicating that the learned dynamical-clock phase provides an operational and quantitatively accurate way to compute the accumulated geometric phase in the weak-perturbation regime. The corresponding phase-space trajectories remain close to the original limit cycle, further supporting the validity of the phase-reduction approximation used in this calculation.

	\begin{figure}[htp]
	\centering
	\includegraphics[width=0.5\textwidth]{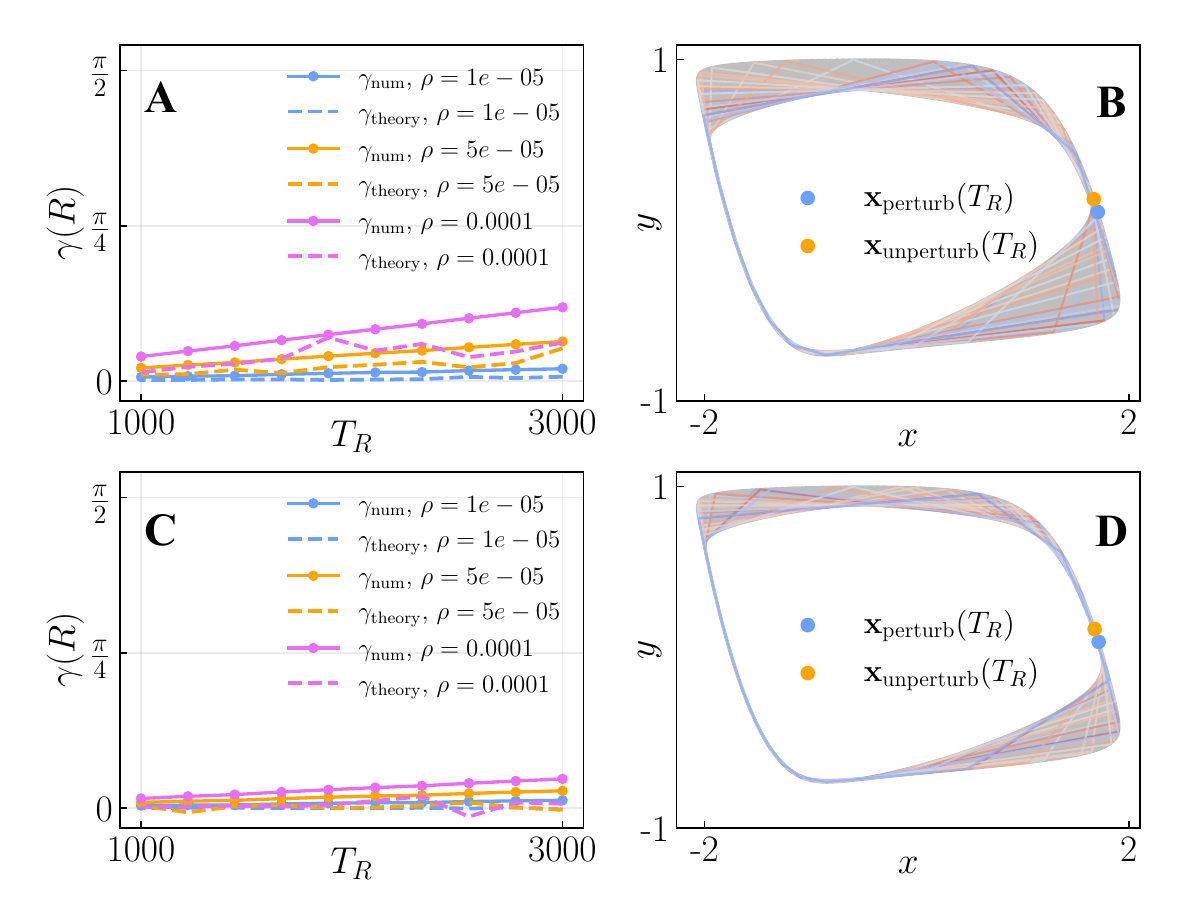}
	\caption{{\color{black}\textbf{Geometric phase under weak cyclic perturbations}. (\textbf{A}, \textbf{C}) Geometric phase of the FitzHugh--Nagumo system under weak cyclic perturbations of the additive parameter $a$ and the multiplicative parameter $b$, respectively. Solid curves denote $\gamma_{\mathrm{num}}$ measured from full numerical simulations, whereas dashed curves denote $\gamma_{\mathrm{theory}}$ predicted by the learned phase dynamics. Different colours indicate different perturbation strengths $\rho$.
		(\textbf{B}, \textbf{D}) Corresponding phase-space trajectories under cyclic perturbations of $a$ and $b$, respectively. The displacement between $x_{\mathrm{perturb}}(T_R)$ and $x_{\mathrm{unperturb}}(T_R)$ reflects the accumulated geometric phase after one perturbation cycle.}}
	\label{fig_geo}
\end{figure}  

}

	\begin{figure*}[htp]
		\centering
		\includegraphics[width=\textwidth]{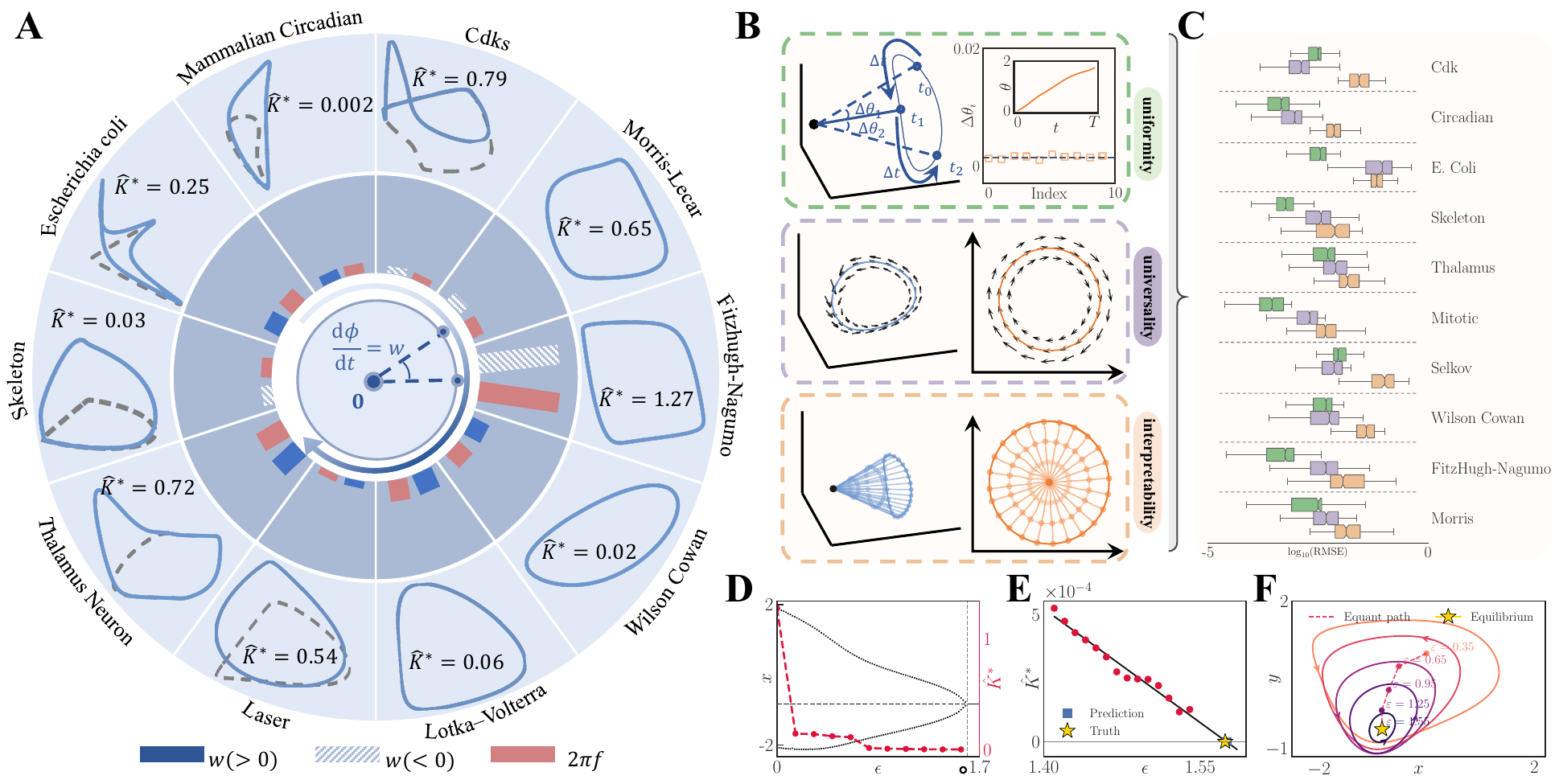}
		\caption{\textbf{Validating the equant across oscillatory systems and predicting critical transitions}. (\textbf{A}) Summary of $10$ different, representative oscillatory systems includes the shape of limit cycle, dimension of ambient space and the equant, Fourier frequency $2\pi f$ and natural frequency $w$.   
			(\textbf{B}) Illustration of the considered performance metrics. {\color{black}(\textbf{C})} Results of accuracy, uniformity and interpretability for the obtained equants in terms of root mean square error (RMSE). {\color{black}(\textbf{D}) Bifurcation diagram of the FitzHugh--Nagumo system with respect to parameter $\varepsilon$. A stable limit cycle surrounds an unstable equilibrium before the tipping point, whereas only a stable equilibrium remains after the transition. Red dots indicate the optimal equant non-uniformity $\hat{K}^\ast$ for each parameter value. 
				(\textbf{E}) Pre-critical linear scaling of $\hat{K}^\ast$ over $\varepsilon\in[1.40,1.55]$, whose extrapolation to $\hat{K}^\ast=0$ accurately predicts the tipping point. 
				(\textbf{F}) As $\varepsilon$ approaches the tipping point, both the limit cycle and the  equant converge towards the equilibrium.}}
		\label{fig2}
	\end{figure*}     
	\subsection{Validating the equant and predicting critical transitions}\label{sec dimension}
	
	\subsubsection{Performance metrics}
	
	We first provide the detailed way to assess the performance of our framework in finding the phase and the equant.
	\begin{itemize}
		\item \textbf{Frequency.} We compare the dominant Fourier frequency $f$ calculated from the time series data with the learned natural frequency $w$ of the phase, and we expect $f-2\pi |w| $ to be close to zero. We note that the sign of $w$ describes the essential rotation direction of the oscillation along the limit cycle.
		\item \textbf{Uniformity.} For regular time series data $\{\rvx_i\}_{i=1}^K$ with time interval $\delta t$, we use the following term as the Accuracy index to test whether the learned equant observes uniform dynamical clock, 
		\begin{equation*}
			\begin{aligned}
				\theta_i&=\angle\left(l_{\boldsymbol{0}\rvg_{\boldsymbol{\theta}_{\rvg}}(\rvx_i)},l_{\boldsymbol{0}\rvg_{\boldsymbol{\theta}_{\rvg}}^{-1}(\rvx_{i+1})}\right),\\
				\delta\theta_i&=\theta_i-\frac{1}{K-1}\sum_{k=1}^{K-1}\theta_k,~i=1,...,K-1.
			\end{aligned}
		\end{equation*}  
		Then we draw the box-plot of $\{\delta\theta_i\}_{i=1}^{K-1}$ in {\color{black}Fig.~\ref{fig2}}.
		
		\item \textbf{Universality.} To ensure the existence of the universal dynamical clock, we not only require the auto-encoder to establish the homeomorphism between the limit cycle and the unit circle, but also expect the auto-encoder to transform the vector field on the limit cycle to the vector field of the dynamical clock on the unit circle. Thus, we assess the universality of the dynamical clock of each system as,  
		\begin{equation*}
			e_i = \left|w-\frac{\partial{\phi}_{\vtheta_\phi}}{\partial \rvx_i} \cdot \rmF(\rvx_i)\right|,~i=1,...,K.
		\end{equation*} 
		Here we require that the data $\{\rvx_i\}_{i=1}^K$ is sampled from the uniform distribution over the limit cycle, and we summarise the statistics of the error $\{e_i\}_{i=1}^K$ by box-plotting  in {\color{black}Fig.~\ref{fig2}}.

		\item \textbf{Interpretability.} According to the Definition~\ref{def 1}, we interpret the phase with the equant by transforming the line segments emanating from the equant to line segments emanating from the centre of the unit circle. To validate whether the learned invertible function $\rvg_{\vtheta_{\rvg}}$ satisfies the conditions in Definition~\ref{def 1}, we consider the direction discrepancy between two unit vectors $\rvv_1$, $\rvv_2$ as 
		\begin{equation*}
			d(\rvv_1,\rvv_2)=\sqrt{\dfrac{1-\rvv_1\cdot\rvv_2|\rvv_1\cdot\rvv_2|}{2}}.
		\end{equation*}  
		The discrepancy takes the minimum $0$ (resp. maximum $1$) when the two vectors have the same (resp. opposite) direction. Then, we test whether $\rvg_{\vtheta_{\rvg}}(l_{\rvx^\ast\rvx})=l_{\vzero\rvg_{\vtheta_{\rvg}}(\rvx)}$ by comparing the discrepancy between $l_{\vzero\rvg_{\vtheta_{\rvg}}(\rvx)}$ and $l_{\vzero\rvg_{\vtheta_{\rvg}}(\rvy)}$ for any point $\rvy\in l_{\rvx^\ast\rvx}$. Specifically, for dataset $\{\rvx_i\}_{i=1}^K$ over the limit cycle, we calculate the discrepancy error as 
		\begin{equation*}
			d_{kj}=d\left(\rvg_{\vtheta_{\rvg}}(\rvy_{kj}),\rvg_{\vtheta_{\rvg}}(\rvx_{k})\right),
		\end{equation*}
		where {\color{black}$\{\rvy_{kj}\}_{j=1}^M$ are regularly arranged points on the line segment $l_{\rvx^\ast\rvx_k}$. Here $M$ denotes the number of interpolation points sampled on each line segment from the equant to a limit-cycle point. We set $M=6$ over all the tested systems, and show their statistics in Fig.~\ref{fig2}}.
		
	\end{itemize}

	\subsubsection{Results on diverse systems}
	From the intricate dance of microscopic cells to the majestic rhythms of celestial bodies, spontaneous oscillations emerge as the most striking phenomenon, often oscillating about an irregular and asymmetrical limit cycle that challenges human perception of complex systems. Using our machine learning framework, we identify the equant and derive the functions that describe the phase evolution for a range of oscillatory systems from electrical circuits to mammalian cell clocks, as illustrated in Fig.~\ref{fig2}A. The detailed formulations of the considered dynamics are provided in Supplementary Materials (SM).
	
	We validate the concept of the dynamical clock by testing it from four perspectives: Frequency, Uniformity, Universality and Interpretability. For \textbf{Frequency}, we consider the discrepancy between the phase frequency and the main Fourier frequency of the oscillatory trajectories in terms of both absolute value and sign. The observer at the equant sees the oscillator moving the same  {\color{black}interpretable phase increment} during the same time interval. We verify the \textbf{Uniformity} of the learned equant by comparing the growth of {\color{black}interpretable phases} observed from equant at different time with constant time intervals. {\color{black}To verify the \textbf{Universality} of the dynamical clock, we require that our method transforms the vector field of an oscillatory system near its limit cycle into a uniform circular flow near the unit circle. We assess the accuracy of this phase reduction by quantifying how closely the transformed dynamics align with ideal uniform rotation, and define universality as the average deviation from this ideal behaviour along the limit cycle. A small deviation indicates a high degree of universality.} For the \textbf{Interpretability}, we measure it by how accurate the extended function transforms line segments emanating from the equant to line segments emanating from the centre of the unit circle. If the extended function is perfectly accurate, line segments emanating from the equant are mapped to straight lines in the unit disk, and the index is close to zero. We summarise the detailed mathematical definitions of these indices in ~\cite{methods}. 
	In Fig.~\ref{fig2}A, we find low frequency (shown as the colorbars) error in terms of absolute value across all the systems. There are however four systems that have opposite signs of phase frequency and Fourier frequency: Skeleton, FitzHugh-Nagumo, Morris-Lecar, and Cdks cell  cycles, representing an essential clockwise rotation property. Figure~\ref{fig2}B illustrates that the performance metrics uniformity universality and interpretability mentioned above, and Figure~\ref{fig2}C shows that our framework performs well in all the tested systems, verifying the accuracy of the located equants and the existence of the dynamical clock.  The results demonstrate the effectiveness of our proposed method for identifying equant for any spontaneous oscillations.

	{\color{black}
In addition to the performance results, we report the values of $\hat{K}^\ast$ for all considered systems, allowing us to compare geometric regularity across different limit cycles. The results reveal two groups with markedly different orders of magnitude. The Wilson-Cowan, Lotka-Volterra, Skeleton, and Mammalian Circadian systems exhibit relatively small values of $\hat{K}^\ast$ ($\sim O(10^{-2})$), whereas the remaining systems show larger values ($\sim O(1)$). 

To further investigate the dynamical implication of $\hat{K}^\ast$, we examine its dependence on the time-scale parameter $\varepsilon$ in the FitzHugh--Nagumo system, as shown in Fig.~\ref{fig2}D--F. Along this parameter route, the system approaches a Hopf transition as $\varepsilon$ increases. We find that the  equant non-uniformity $\hat K^\ast$ provides a geometric indicator of this transition: the optimal value $\hat{K}^\ast$ decreases towards zero as the system approaches the critical point. This suggests that, within this family of oscillators, smaller values of $\hat{K}^\ast$ correspond to a limit cycle that is closer to the impending transition, whereas larger values reflect a finite-amplitude relaxation oscillation that remains farther from the Hopf point.

To quantify this behaviour, we track both the optimal equant and the corresponding $\hat{K}^\ast$ learned by our framework. As $\varepsilon$ increases, the FitzHugh--Nagumo limit cycle gradually loses its relaxation-oscillation geometry and shrinks towards the equilibrium associated with the Hopf transition. Consistently, the optimal equant moves towards the equilibrium (Fig.~\ref{fig2}F), while $\hat{K}^\ast$ decreases monotonically over the tested parameter range (Fig.~\ref{fig2}D). This behaviour indicates that $\hat{K}^\ast$ captures the progressive loss of slow-fast geometric distortion and can therefore serve as an early-warning signal for the approaching transition along this parameter route.

More importantly, the pre-critical decay of $\hat{K}^\ast$ exhibits an approximately linear scaling with respect to $\varepsilon$. By fitting the relation between $\hat{K}^\ast$ and $\varepsilon$ over the range $\varepsilon\in[1.40,1.55]$ and extrapolating to $\hat{K}^\ast=0$, we obtain a predicted critical value $\varepsilon_c=1.575$, in close agreement with the true Hopf point $\varepsilon_c=1.574$ (Fig.~\ref{fig2}E). This result suggests that the equant not only provides a physical interpretation for the dynamical clock, but also offers a practical route for tipping-point prediction from pre-critical oscillatory trajectories.

}

	\begin{figure*}[t]
		\centering
		\subfigure{\label{fig_threea}}
		\includegraphics[width=\textwidth]{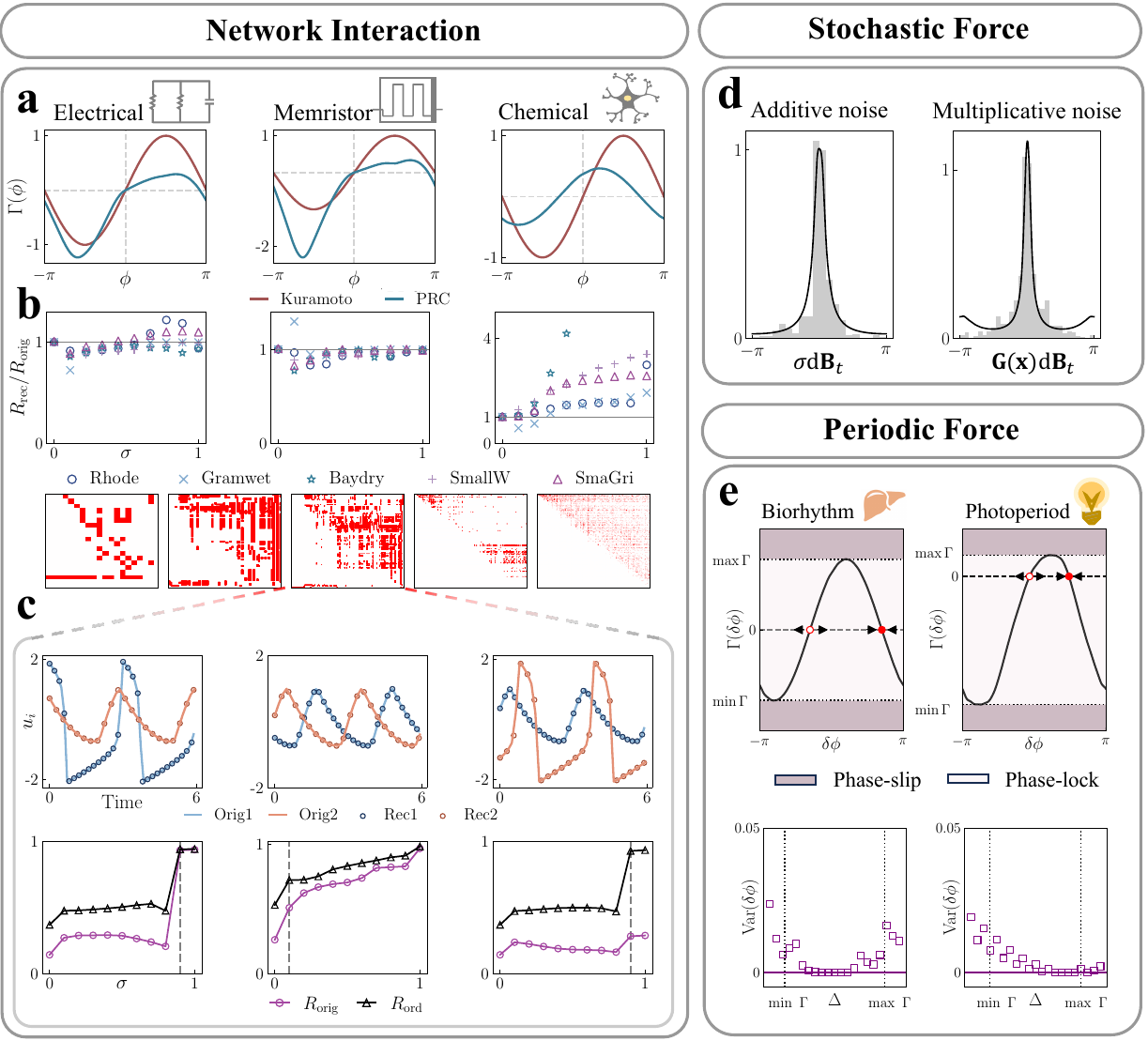}
		\caption{\textbf{Phase synchronisation mechanism against various external forces}. (\textbf{A}) The phase response functions against different coupling functions. (\textbf{B}) The statistics under each coupling function over five realistic networks. (\textbf{C}) We fix the network Baydry and compare the original and the reconstructed trajectories (top), and the scaling form of the synchronisation indicators (bottom). (\textbf{D}) The theoretical prediction (solid line) and the numerical results (shaded areas) of the stationary distribution.  (\textbf{E}) The phase response functions against different periodic forces (top), and the corresponding Statistical validation (bottom).}
		\label{figFHN}
	\end{figure*}

	\subsection{Synchronisation mechanism of neuronal activity}\label{sec FHN}
	
	Phase synchronisation phenomenon of rhythmic oscillations has been extensively studied across scales and disciplines in the context of collective dynamics  \cite{rosenblum2001phase,strogatz2018nonlinear,varela2001brainweb,sauseng2008does,palva2005phase}. Among the plentiful investigations, few have elucidated the mechanisms of phase synchronisation under different types of perturbations. To address this, here we illustrate how our machine learning-based phase reduction framework reveals the influence of external perturbations at the level of phase dynamics, using the FitzHugh-Nagumo system.  The detailed experimental settings are summarised in Appendix~\ref{appendix FHN}.

	We first consider the \textbf{network interactions} arising from couplings between nodes distributed over complex networks. Regarding the network structure, we select five realistic networks with sizes ranging from $20$ to $1059$. For the coupling function, we examine three practical types: electrical coupling~\cite{li2007effects}, memristor-based coupling~\cite{volos2015memristor}, and chemical synaptic coupling~\cite{korotkov2019dynamics}. Figure~\ref{figFHN}A illustrates that the phase response functions $\Gamma$ for these couplings differ markedly from that of the classical Kuramoto model. As in the \textit{E. coli} case, we employ the synchronisation degrees $R_{\text{orig}}$ and $R_{\text{pred}}$ to quantitatively assess the effects of couplings under varying coupling strengths $\sigma$ and network structures. As shown in Fig.~\ref{figFHN}B, $R_{\text{pred}}$ closely matches $R_{\text{orig}}$ for electrical and memristor-based couplings across most coupling strengths and networks, whereas notable deviations occur for chemical coupling as $\sigma$ increases. This deviation results from the coupling term does not vanish even when oscillators are synchronised, thereby distorting the original oscillation shape. This property is also evident in the phase response functions, where $\Gamma(0) = 0$ for electrical and memristor couplings, but $\Gamma(0) \neq 0$ for chemical coupling. Consequently, for vanishing couplings, our framework performs reliably across coupling strengths. In contrast, for non-vanishing couplings, a weak coupling strength is necessary to constrain perturbations near the original limit cycle. The accuracy of the reduced phase dynamics is further verified and shown in Fig.~\ref{figFHN}C (top), where the original and the reconstructed trajectories perfectly coincide.  Moreover, Fig.~\ref{figFHN}C (bottom) shows that the order parameter $R_{\text{ord}}$ and the synchronization indicator $R_{\text{orig}}$ exhibit the same tipping point as a function of $\sigma$, demonstrating the capability of our framework to detect critical transitions.

	Next, we consider a cluster of oscillators driven by a \textbf{stochastic force}, incorporating both additive and multiplicative noise. To capture the impact of noise on the collective behaviour of the oscillators, we theoretically derive the stationary distribution $p^\ast$ of their phase differences by solving the Fokker-Planck equation associated with the stochastic phase dynamics obtained via our phase reduction framework. To validate the theoretical predictions, we numerically compute the empirical distribution of 5000 oscillators over a long time span, repeated across 100 independent simulations. As shown in Fig.~\ref{figFHN}D, the theoretical predictions closely match the numerical results.
	
	Finally, we investigate whether the \textbf{periodically forced} oscillator exhibits a phase-lock phenomenon~\cite{hsieh1996phase} in response to external forces originating from natural rhythms such as biorhythms and photoperiods. Specifically, we examine whether the frequency of the forced oscillator converges from its natural frequency $w$ to the external frequency $\Omega$. The asymptotic behaviour of the forced phase is governed by the relationship between $\Delta = w - \Omega$ and $\Gamma$, as illustrated in Fig.~\ref{figFHN}E (top). As $\Delta$ increases, the system first undergoes a saddle-node bifurcation, then the bifurcation separates into a stable equilibrium (solid dot) and an unstable equilibrium (hollow dot), before eventually merging back into a saddle-node as the equilibria vanish. Phase-lock occurs when a stable equilibrium exists (light region), whereas in the absence of a stable equilibrium, the phase frequency diverges, leading to a phase-slip (dark region).  Otherwise, the phase frequency diverges and is characterised as a phase-slip (dark region). To validate this, we uniformly sample different $\Delta$ values in each region and numerically calculate the frequency shift between the forced oscillator and the external force. As shown in Fig.~\ref{figFHN}E (bottom), the results confirm the division into phase-lock and phase-slip regimes.

	\section{Discussion}\label{sec discussion}
	
	{\color{black} Drawing on Ptolemy’s equant and an areal-uniformity principle inspired by Kepler’s second law, we establish the geometric structure behind spontaneous rhythms that pervade nature.} This yields a unifying framework that maps nonlinear, high-dimensional oscillations onto a two-dimensional uniform circular motion, termed the `dynamical clock', which provides a physically grounded interpretation of phase, enabling analytical phase response functions under external perturbations. While Ptolemy's original equant offered an elegant yet speculative explanation of celestial motion, its broader relevance to general nonlinear oscillations has remained largely unexplored. Here, we establish that the equant provides a powerful lens for interpreting the phase dynamics of general nonlinear oscillators.
	{\color{black}
	We further introduce an areal-uniformity criterion for selecting the equant, which provides a quantitative geometric index
	for comparing limit-cycle regularity and detecting proximity to critical transitions. }
	 The key challenge -- {\color{black}identifying the location of the equant and the associated nonlinear viewing map that render the motion on an irregular limit cycle uniformly representable on the universal dynamical clock} -- is overcome by integrating manifold learning and generative modeling techniques from modern machine learning. Our results not only advance phase reduction theory but also illustrate how AI can assist in uncovering physically meaningful representations of complex dynamical phenomena~\cite{iten2020discovering,wang2023scientific}.
	
	As a direct application of our framework, we address the open question regarding the collective behaviour in \textit{E. coli} cells under quorum sensing,  uncovering the mechanism of how cell-to-cell communication determines the onset of synchronisation in genetic oscillators. Next, we precisely characterise the response of rhythmic gene expression against varying external perturbations solely using the synthetic data from a canonical Repressilator genetic circuit.
	{\color{black}
	Although several derivations are presented using additive perturbations for clarity, this is not an intrinsic restriction of the framework. Nonlinear couplings, multiplicative parameter perturbations, or heterogeneous intrinsic dynamics can be treated as effective weak dynamical perturbations around the reference limit cycle, provided that the perturbed trajectories remain in its neighbourhood. Their influence on the phase is then obtained by projecting the induced vector-field difference onto the learned phase sensitivity. This interpretation is supported by the neuronal examples with memristor and chemical synaptic couplings in Results~\ref{sec FHN}, where the coupling terms are nonlinear but the reduced phase dynamics remain accurate within the weak-coupling regime.}

	We further verify the existence of the dynamical clock among diverse oscillatory scenarios.
	{\color{black}
	Beyond validating the existence of the equant, the equant-based areal-uniformity criterion further provides a quantitative way to compare the geometric regularity of different limit cycles. In the FitzHugh--Nagumo example, the optimal equant non-uniformity $\hat{K}^\ast$ decreases as the system approaches the super-critical Hopf bifurcation, suggesting that the equant can also serve as a geometric diagnostic for detecting proximity to critical transitions.}
	
	
	We believe that there are ample opportunities for future research in applying this framework to investigate emergent behaviour in various fields. Potential application areas include phase adjustment of cardiac pacemaker cells~\cite{monga2019optimal}, frequency modulation of biorhythms~\cite{zhong2023modulating}, developing treatments for neurological disorders stemming from abnormal synchrony~\cite{wang2023desynchronizing}, and revealing the mechanism of genetic clocks~\cite{kim2012mechanism}. Furthermore, it's  also interesting to investigate the connection between our interpretable phase and the Berry geometric phase in quantum mechanics. More broadly, with the rapid growth of machine learning techniques, we aspire that our framework could be used to gain deeper insights in interdisciplinary fields related to AI for science across disciplines as varied as biology, physics and mathematics. {\color{black}Another promising direction is to develop
		equant-based areal-uniformity indicators, such as $\hat{K}^\ast$ defined in this work, for early-warning analysis of transitions in biological and neuronal oscillators.}

	{\color{black}
Several limitations also point to natural future directions. First, the present implementation focuses on phase dynamics on the limit cycle. A natural extension is to augment the circular latent phase with transverse amplitude coordinates constrained by Floquet modes or Koopman eigenfunctions, thereby connecting the framework to phase-amplitude reduction~\cite{shirasaka2017phase,wilson2018greater,wilson2019augmented}. Such an extension would require off-cycle data or additional regularisation. Second, for systems without closed-form governing equations and only the trajectory data is available, the dynamics may be reconstructed using data-driven methods such as Neural ODEs developed in the Supplementary Materials, SINDy, or reservoir computing.  For partial observation data, additional state reconstruction, for example through Koopman-based lifting~\cite{zhang2024learning} or delay-coordinate embedding~\cite{botvinick2025invariant}, is required before applying the phase-learning procedure.	
}

	\section*{Acknowledgements}
	L.Y. is supported by the CSC (No. 202406100251). Q. Z. is supported by the China Postdoctoral Science Foundation (No. 2022M720817), by the Shanghai Postdoctoral Excellence Program (No. 2021091), and by the STCSM (Nos. 21511100200, 22ZR1407300, 22dz1200502, and 23YF1402500). W.L. is supported by the
	National Natural Science Foundation of China (Grant nos. 11925103,
	12531018), by the STCSM (Grant nos. 2021SHZDZX0103, 22JC1402500,
	and 22JC1401402), and by the SMEC (Grant no. 2023ZKZD04).

	\appendix

	{\color{black}
\section{Theoretical Derivations}\label{appendix deri}	

\textit{Derivation for Eq.~\eqref{eq common+independ}.}---To derive the analytical formulation of the stationary distribution of phase difference from Eq.~\eqref{eq common+independ}. We recall the phase dynamics is 
\begin{equation*}
	\begin{aligned}
		\mathrm{d}{\phi}_i &= \left(w+\delta w\right)\mathrm{d}t+ \Gamma_1(\phi_i)\mathrm{d}\xi+\Gamma_2(\phi_i)\mathrm{d}\eta_i,\\
		\delta w &= \dfrac{1}{2\pi}\int_0^{2\pi}\frac{1}{2}\text{Tr}\left[(\rmG+\rmH)^\top \rmY(\theta)(\rmG+\rmH)\right]\mathrm{d}\theta.\\
	\end{aligned}
\end{equation*}
Here $\Gamma_1(\phi)=\rmZ(\phi)\cdot \rmG(\boldsymbol{\chi}(\phi))$, $\Gamma_2(\phi)=\rmZ(\phi)\cdot \rmH(\boldsymbol{\chi}(\phi))$, and the Wiener processes $\xi$ and $\{\eta_i\}_{i=1}^N$ are common noise and independent noise, respectively. We first consider the covariance matrix $\rmD=(D_{ij})\in\mathbb{R}^{N\times N}$ of $\{\mathrm{d}\phi_i\}_{i=1}^N$ as,
\begin{equation*}
	\begin{aligned}
		D_{ij}(\boldsymbol{\phi})&=\Gamma_1(\phi_i)\Gamma_1(\phi_j)+\delta_{ij}\Gamma_2(\phi_i)\Gamma_2(\phi_j),\\
		\boldsymbol{\phi}&=(\phi_1,...,\phi_N)^\top.
	\end{aligned}
\end{equation*}
Next, we consider the Fokker-Planck equation of distribution $q$ of slow phase component $\psi_i=\phi_i-(w+\delta w)t$ as, 
\begin{equation*}
	\partial_{t}q(\boldsymbol{\psi},t)=\frac{1}{2}\sum_{i=1}^N\sum_{j=1}^N\partial_{\psi_i}\partial_{\psi_j}\left[D_{ij}(\boldsymbol{\psi}+(w+\delta w)t)q(\boldsymbol{\psi},t)\right].
\end{equation*}
By time-averaging the covariance $D_{ij}$ over one period $T=\frac{2\pi}{w+\delta w}$, we obtain the simplified Fokker-Planck equation,
\begin{equation}\label{eq PR noise diff FK}
	\begin{aligned}
		\partial_{t}q(\boldsymbol{\psi},t)&=\frac{1}{2}\sum_{i=1}^N\sum_{j=1}^N\partial_{\psi_i}\partial_{\psi_j}\big[\big(\gamma_1(\psi_i-\psi_j)\\
		&\quad+\delta_{ij}\gamma_2(\psi_i-\psi_j) \big)q(\boldsymbol{\psi},t)\big].\\
		\gamma_i(\psi)&=\frac{1}{2\pi}\int_0^{2\pi}\Gamma_i(\theta)\Gamma_i(\psi+\theta)\mathrm{d}\theta,~i=1,2.
	\end{aligned}
\end{equation}

Without loss of generality, we focus on the phase difference $\delta\phi=\phi_1-\phi_2=\psi_1-\psi_2$. The evolution equation for distribution of $\psi_1$, $\psi_2$ can be derived
from $q(\boldsymbol{\psi})$ by integrating over all other phase variables as, 
\begin{equation*}
	R(\psi_1,\psi_2,t)=\int q(\boldsymbol{\psi})\mathrm{d}\psi_3\cdots\psi_N.
\end{equation*}
According to Eq.~\eqref{eq PR noise diff FK}, the evolution equation of $R(\psi_1,\psi_2)$ is, 
\begin{equation}\label{eq PR noise 1,2}
	\begin{aligned}
		\partial_tR(\psi_1,\psi_2,t)&=\frac{\gamma_1(0)+\gamma_2(0)}{2}\partial_{\psi_1}^2R+\frac{\gamma_1(0)+\gamma_2(0)}{2}\partial_{\psi_2}^2R\\
		&+\partial_{\psi_1}\partial_{\psi_2}\left[\gamma_1(\delta\phi)R\right].
	\end{aligned}
\end{equation}
By transforming the two variables $\psi_1,~\psi_2$ to the mean $\psi=\frac{\psi_1+\psi_2}{2}$ and the difference $\delta\phi=\psi_1-\psi_2$, we have 
\begin{equation*}
	\begin{aligned}
		R(\psi_1=\psi+\frac{\delta\phi}{2},\psi_2=\psi-\frac{\delta\phi}{2},t)=Q(\psi,t)P(\delta\phi,t).
	\end{aligned}
\end{equation*}

Then we have 
\begin{equation*}
	\begin{aligned}
		\partial_tR&=P\partial_tQ+Q\partial_tP\\
		\partial_{\psi_1}^2R&=\frac{1}{4}P\partial_\psi^2Q+(\partial_\psi Q)(\partial_{\delta\phi} P)+Q\partial_{\delta\phi}^2 P,\\
		\partial_{\psi_2}^2R&=\frac{1}{4}P\partial_\psi^2Q-(\partial_\psi Q)(\partial_{\delta\phi} P)+Q\partial_{\delta\phi}^2 P,\\
		\partial_{\psi_1}\partial_{\psi_2}\left[\gamma_1(\delta\phi)R\right]&=\frac{1}{4}\gamma_1(\delta\phi)P\partial_\psi^2Q-Q\partial_{\delta\phi}^2 \left[\gamma_1(\delta\phi)P\right].
	\end{aligned}
\end{equation*}

Thus, the Eq.~\eqref{eq PR noise 1,2} is decoupled as,

\begin{equation*}
	\begin{aligned}
		\partial_t Q&=\frac{1}{4}\left[\gamma_1(0)+\gamma_1(\delta\phi)+\gamma_2(0)\right]\partial_{\psi}^2Q,\\
		\partial_t P&=\partial_{\delta\phi}^2\left\{\left[\gamma_1(0)-\gamma(\delta\phi)+\gamma_2(0)\right]P\right\}.
	\end{aligned}
\end{equation*}

Then the stationary distribution of phase difference $P(\delta \phi)$ is determined by, 
\begin{equation*}
	P^\ast(\delta\phi)=\dfrac{p_0}{\gamma_1(0)-\gamma_1(\delta\phi)+\gamma_2(0)},
\end{equation*}
where $p_0$ is the normalisation constant. 

\textit{Derivation for Eq.~\eqref{eq PR network final},\eqref{eq PR network heter final}.}---Specifically, for Eq.~\eqref{eq PR network} we split the phase variable as $\phi_i=wt+\psi_i$. Since we require the perturbation, i.e., the network interaction to be weak, the contribution of the perturbation to the phase frequency is weaker than natural frequency $w$. Thus, $wt$ is the fast component of the phase and $\psi_i$ represents the slow component. Then we have,
\begin{equation}\label{eq PR network slow}
	\begin{aligned}
		\dot{\phi}_i&=w+\dot{\psi}_i,\\
		\dot{\psi}_i&=\sum_{j=1}^NA_{ij}\rmZ(wt+\psi_i)\cdot\rmH(\boldsymbol{\chi}(wt+\psi_i),\boldsymbol{\chi}(wt+\psi_j)).\\
	\end{aligned}
\end{equation}
Since the component $\psi_i$ varies slowly during one cycle of oscillation with frequency $w$, we treat the $\psi_i$, $\psi_j$ on the right hand side of Eq.~\eqref{eq PR network slow} as constant over one period of oscillation. Then the vector field of $\psi_i$ on the left hand side of Eq.~\eqref{eq PR network slow} solely varies along $wt$ with period $T=\frac{2\pi}{w}$, the dynamics of which is equal to the time-average over this period under fixed $\psi_i$, $\psi_j$,
\begin{equation}\label{eq PR network average}
	\resizebox{\linewidth}{!}{$
		\begin{aligned}
			\dot{\psi}_i&=\frac{1}{T}\int_0^T\sum_{j=1}^NA_{ij}\rmZ(ws+\psi_i)\cdot\rmH(\boldsymbol{\chi}(ws+\psi_i),\boldsymbol{\chi}(ws+\psi_j))\mathrm{d}s\\
			&=\frac{1}{2\pi}\int_0^{2\pi}\sum_{j=1}^NA_{ij}\rmZ(\theta+\psi_i)\cdot\rmH(\boldsymbol{\chi}(\theta+\psi_i),\boldsymbol{\chi}(\theta+\psi_j))\mathrm{d}\theta\\
			&=\sum_{j=1}^NA_{ij}\left[\frac{1}{2\pi}\int_0^{2\pi}\rmZ(\theta+\psi_i)\cdot\rmH(\boldsymbol{\chi}(\theta+\psi_i),\boldsymbol{\chi}(\theta+\psi_j))\mathrm{d}\theta\right]\\
			&\triangleq \sum_{j=1}^NA_{ij}\Gamma(\psi_j-\psi_i).
		\end{aligned}
		$}
\end{equation}
We note that $\phi_j-\phi_i=\psi_j-\psi_i$, and finally the Eq.~\eqref{eq PR network final} is obtained by putting Eq.~\eqref{eq PR network average} into Eq.~\eqref{eq PR network slow}.

We consider the homogeneous oscillators in the above analysis, next we extend our framework to heterogeneous oscillators as follows,
\begin{equation}\label{eq nonidentical}
	\begin{aligned}				 
		\dot{\rvx}_i&=\rmF_i(\rvx_i)+\sum_{j=1}^NA_{ij}\rmH(\rvx_i,\rvx_j)\\
		\rmF_i(\rvx_i)&=\rmF(\rvx_i)+\delta\rmF_i(\rvx_i),~\dfrac{\|\delta\rmF_i(\rvx_i)\|}{\|\rmF_i(\rvx_i)\|}\ll1.
	\end{aligned}
\end{equation}
The corresponding phase dynamics takes the form as,
\begin{equation*}
	\dot{\phi}_i=w+\rmZ(\phi_i)\cdot\delta \rmF_i(\boldsymbol{\chi}(\phi_i))+\sum_{j=1}^{N}A_{ij}\Gamma(\phi_j-\phi_i).
\end{equation*}
Similarly we consider the dynamics of slow component $\psi_i=\phi_i-wt$ as,

\begin{equation}\label{eq PR network heter phase}
	\dot{\psi}_i=\rmZ(wt+\psi_i)\cdot\delta \rmF_i(\boldsymbol{\chi}(wt+\psi_i))+\sum_{j=1}^{N}A_{ij}\Gamma(\psi_j-\psi_i).
\end{equation}
The first term on the right-hand side of Eq.~\eqref{eq PR network heter phase} is replaced by the time-average over one period of oscillation,
\begin{equation}\label{eq PR network heter average}
	\begin{aligned}
		\dot{\psi}_i&=\frac{1}{T}\int_0^T\rmZ(ws+\psi_i)\cdot\delta \rmF_i(\boldsymbol{\chi}(ws+\psi_i))\mathrm{d}s\\
		&\quad\quad+\sum_{j=1}^{N}A_{ij}\Gamma(\psi_j-\psi_i)\\
		&=\frac{1}{2\pi}\int_0^{2\pi}\rmZ(\theta+\psi_i)\cdot\delta \rmF_i(\boldsymbol{\chi}(\theta+\psi_i))\mathrm{d}s\\
		&\quad\quad+\sum_{j=1}^{N}A_{ij}\Gamma(\psi_j-\psi_i)\\
		&=\underbrace{\frac{1}{2\pi}\int_0^{2\pi}\rmZ(\theta)\cdot\delta \rmF_i(\boldsymbol{\chi}(\theta))\mathrm{d}s}_{(\rmZ,\boldsymbol{\chi}~\text{is periodic by $2\pi$})}+\sum_{j=1}^{N}A_{ij}\Gamma(\psi_j-\psi_i).
	\end{aligned}
\end{equation}
Then we obtain the reduced phase dynamics of Eq.~\eqref{eq nonidentical} as follows,
\begin{equation*}
	\begin{aligned}
		\dot{\phi}_i&=(w+\delta w_i)+\sum_{i=1}^{N}A_{ij}\Gamma(\phi_j-\phi_i),\\
		\delta w_i&=\rmZ(\phi_i)\cdot\delta\rmF(\rvx_i)=\dfrac{1}{2\pi}\int_0^{2\pi}\rmZ(\theta)\cdot\delta\rmF_i(\boldsymbol{\chi}(\theta))\mathrm{d}\theta.
	\end{aligned}
\end{equation*}

}
	
	\section{Addressing the open problem: Quorum sensing scaling law of \textit{E. coli} cell population}\label{appendix Ecoli}
	
	We consider the dynamics of engineering E. coli cell population with $N$ cells as~\cite{garcia2004modeling}, 
	\begin{equation}\label{eq E.coli1}
		\begin{aligned}
			\dot{a}_i&=-a_i+\dfrac{\alpha}{1+C_i^n},\\
			\dot{b}_i&=-b_i+\dfrac{\alpha}{1+A_i^n},\\
			\dot{c}_i&=-c_i+\dfrac{\alpha}{1+B_i^n}+\dfrac{\kappa S_i}{1+S_i},\\
			\dot{A}_i&=\beta_i(a_i-A_i),\\
			\dot{B}_i&=\beta_i(b_i-B_i),\\
			\dot{C}_i&=\beta_i(c_i-C_i),\\
			\dot{S}_i&=-k_{s0}S_i+k_{s1}A_i-\eta(S_i-Q\bar{S}),\\
		\end{aligned}
	\end{equation}
	where $a_i$, $b_i$, and $c_i$ are the concentrations of mRNA transcribed from genes tetR, cI, and lacI in $i$-th cell, respectively, and the concentration of the corresponding proteins are represented by $A_i$, $B_i$, and $C_i$. The cell-to-cell communication in the population, i.e., the quorum sensing, is carried out by a small molecule known as autoinducer. The  concentration of autoinducer inside each cell is denoted by $S_i$, and the quorum sensing is represented as the mean field coupling $\bar{S}=\frac{1}{N}\sum_{j=1}^NS_j$. The coupling intensity $Q\in[0,1]$ relates to the cell density in the population. We set the parameters in Eq.~\eqref{eq E.coli1} as $\alpha=216$, $n=2$, $\kappa=20$, $\beta_i=\beta=1$, $k_{s0}=1$, $k_{s1}=0.01$, and $\eta=2$.
	
	For the stochastic perturbation to the original dynamics, we consider the intrinsic {\color{black}heterogeneity} and the extrinsic noise, respectively, to the parameter $\beta_i$, which is the ratio between the mRNA and the protein lifetimes. The intrinsic {\color{black}heterogeneity} is reflected by $\beta_i=\beta+\delta\beta_i$ with $\delta\beta_i$ satisfying the normal distribution $\mathcal{N}(0,\sigma^2)$, capturing the intrinsic differences in different cells. The extrinsic noise is represented as $\beta_i=\beta+\xi_i$, wherein $\xi_i$ is the Gaussian white noise. 
	
	Firstly, we find the accurate phase dynamics of the coupled cellular dynamics in Eq~\eqref{eq E.coli1}. By regarding the mean field coupling interaction as the external perturbation, our framework reduces the phase dynamics on the limit cycle of the uncoupled dynamics with $Q=0$. We simulate the temporal trajectories $\{\phi_i(t)\}_{i=1}^N$ of the reduced phase dynamics, and obtain the reconstructed trajectories in the original space as $\{\boldsymbol{\chi}_{\vtheta_{\boldsymbol{\chi}}}(\phi_i(t))\}_{i=1}^N$.We collect the original trajectories by simulating the dynamics in Eq.~\eqref{eq E.coli1}.  We then assess the effectiveness of our framework by calculating the $R^2$ score between the reconstructed trajectories and the original trajectories. However, we note that the performance of this approach is not satisfactory especially when $Q$ increases to $1$. The reason is that, when the coupled oscillators synchronise at a high cell density $Q$, the coupling interaction $Q\bar{S}$ becomes stronger. {\color{black}This forces the oscillations into a new limit cycle, whose dynamical properties are not prescribed in advance.}
	
	To address this issue, we re-organize the coupled dynamics such that the coupling term maintains the original limit cycle under synchronisation,
	\begin{equation}\label{eq E.coli2}
		\dot{S}_i=\underbrace{-k_{s0}S_i+k_{s1}A_i-\eta(1-Q)S_i}_{\text{Self-dynamics}}+\underbrace{\eta Q(\bar{S}-S_i)}_{\text{Coupling}}.
	\end{equation}
	In this way, the coupling term vanishes under synchronisation at the expense of introducing $Q$-dependent self-dynamics. We denote by  $\rmF_Q$ the $Q$-dependent self-dynamics and by $\Gamma_Q$ the corresponding phase response function. We sample the regular grids $\{Q_i\}_{i=1}^K$ in $[0,1]$, and assign an auto-encoder  for each $\rmF_Q$ to learn the functions $\{\phi_{Q},\boldsymbol{\chi}_Q,\rmZ_Q\}$ and the corresponding phase dynamics as,
	\begin{widetext}
		\begin{equation}\label{eq coupled E.coli phase}
			\begin{aligned}
				\dot{\phi}_{Qi}&=w_{Qi}+\dfrac{\eta Q}{N}\sum_{j=1}^N\Gamma_{Q}\left(\phi_{Qj}-\phi_{Qi}\right)\\
				w_{Qi}&=w_Q,\\
				\Gamma_{Q}\left(\phi_{Qj}-\phi_{Qi}\right)&=\dfrac{1}{2\pi}\int_0^{2\pi}\rmZ_Q(\theta+\phi_i)\rmH(\boldsymbol{\chi}_Q(\theta+\phi_i),\boldsymbol{\chi}_Q(\theta+\phi_j))\mathrm{d}\theta,\\
				\rmH(\rvx_i,\rvx_j)&=(0,0,0,0,0,0,S_j-S_i)^\top.
			\end{aligned}
		\end{equation}
	\end{widetext}
	Secondly, we assess the accuracy of the reduced phase dynamics under noisy perturbation. For the intrinsic {\color{black}heterogeneity}, the perturbation leads to non-identical oscillators. Therefore, we employ the derivation in Eq.~\eqref{eq nonidentical} to deduce the reduced phase dynamics under each cell density $Q$ as,
	\begin{equation}\label{eq coupled E.coli phase intrinsic}
		\begin{aligned}
			\dot{\phi}_{Qi}&=w_{Qi}+\dfrac{\eta Q}{N}\Gamma_{Q}(\phi_{Qj}-\phi_{Qi}),\\
			w_{Qi}&=w_{Q}+\delta w_{Qi},\\
			\delta\rmF_i(\rvx_i)&=(0,0,0,\delta\beta_i(a_i-A_i),\delta\beta_i(b_i-B_i),\delta\beta_i(c_i-C_i),0)^\top,\\
			\delta w_{Qi}&=\dfrac{1}{2\pi}\int_0^{2\pi}\rmZ_Q(\theta)\cdot\delta\rmF_i(\boldsymbol{\chi}_Q(\theta))\mathrm{d}\theta.
		\end{aligned}
	\end{equation}
	
	We simulate the phase dynamics in Eq.~\eqref{eq coupled E.coli phase} and obtain the reconstructed trajectories $\{\boldsymbol{\chi}_Q(\phi_i(t))\}$. For measuring the statistics of the trajectories, we employ the synchronisation index, introduced in~\cite{garcia2004modeling}, as
	\begin{equation}\label{eq ind1}
		\begin{aligned}
			R&=\dfrac{\langle M^2\rangle-\langle M\rangle^2}{\frac{1}{N}\sum_{i=1}^N\left(\langle b_i^2\rangle-\langle b_i\rangle^2\right)},\\
			M&=\dfrac{1}{N}\sum_{i=1}^Nb_{i}(t),
		\end{aligned}
	\end{equation}
	where $\langle\cdot\rangle$ denotes the time average.
	We denote the index calculated from the reconstructed  (resp. original) trajectories as $R_{\text{pred}}$ (resp. $R_{\text{orig}}$). Then we test the accuracy of the phase dynamics by calculating the similarity between the  reconstructed and the original trajectories as 
	\begin{equation}
		\text{Sim}(R_{\text{rec}},R_{\text{orig}})=1-\dfrac{|R_{\text{rec}}-R_{\text{orig}}|}{R_{\text{rec}}+R_{\text{orig}}},
	\end{equation} 
	wherein $\text{Sim}(R_{\text{rec}},R_{\text{orig}})=1$ marks the accurate prediction of the original dynamics with the phase dynamics. Since the index $R$ depends  on both the cell density $Q$ and the cell diversity $\sigma$, we test the similarity under different combinations of $Q$ and $\sigma$. We further measure the synchronisation degree of the coupled phase dynamics in Eq.~\eqref{eq coupled E.coli phase} with the order parameter as
	\begin{equation}\label{eq ind2}
		R_{\text{ord}}=\left\langle\bigg|\dfrac{1}{N}\sum_{i=1}^Ne^{\mathrm{i}\phi_{i}(t)}\bigg|\right\rangle.
	\end{equation}
	And we also calculate the similarity between $R_{\text{ord}}$ and $R_{\text{orig}}$.
	
	To proceed, we consider the extrinsic noise driven by the stochastic process $\beta_i = \beta+\xi_{\mu i}$, $\mu=a,b,c$,
	where $\xi_{\mu i}$ are the Ornstein–Uhlenbeck process with zero mean and intensity $\sigma$, i.e., $\xi_{\mu i}(t)\xi_{\nu j}(s)=\delta_{\mu\nu}\delta_{ij}\left(\frac{\sigma}{2\kappa}\right)\exp(-\kappa|t-s|)$. The corresponding formulation of the It$\hat{\text{o}}$'s type SDE is described as 
	\begin{equation}\label{eq E.coli OU1}
		\begin{aligned}
			\mathrm{d}A_i&=\beta(a_i-A_i)\mathrm{d}t+(a_i-A_i)\mathrm{d}\xi_{ai}(t),\\
			\mathrm{d}B_i&=\beta(b_i-B_i)\mathrm{d}t+(b_i-B_i)\mathrm{d}\xi_{bi}(t),\\
			\mathrm{d}C_i&=\beta(c_i-C_i)\mathrm{d}t+(c_i-C_i)\mathrm{d}\xi_{ci}(t),\\
			\mathrm{d}\xi_{\mu i}&=-\kappa\xi_{\mu i}\mathrm{d}t+\sqrt{\sigma}\mathrm{d}W_{\mu i}(t),~\mu=a,b,c,\\
		\end{aligned}
	\end{equation}
	with independent standard Wiener process $W_{ \mu i}$, $\mu=a,b,c$, $1\le i \le N$.
	The corresponding stochastic  phase dynamics is,
	\begin{equation}\label{eq E.coli OU2}
		\begin{aligned}
			\mathrm{d}\phi_{Qi}&=\left[w_{Qi}+\dfrac{\eta Q}{N}\sum_{i=1}^N\Gamma_Q(\phi_j-\phi_i)\right]\mathrm{d}t\\
			&\quad\quad+
			\rmZ_{Q}(\phi_i)\cdot\left[(\mathrm{d}\boldsymbol{\xi})\rmG(\boldsymbol{\chi}(\phi_i))\right],\\
			w_{Qi}&=w_{Q}+\delta w_{Qi},\\
			\mathrm{d}\xi_{\mu i}&=-\kappa\xi_{\mu i}\mathrm{d}t+\sqrt{\sigma}\mathrm{d}W_{\mu i}(t),~\mu=a,b,c,\\
		\end{aligned}
	\end{equation}
	where $\delta w_{Qi} $ in obtained similarly to Eq.~\eqref{eq SDE phase} with $\rmG(\rvx_i)=(0,0,0,\sqrt{\sigma}(a_i-A_i),\sqrt{\sigma}(b_i-B_i),\sqrt{\sigma}(c_i-C_i),0)^\top$, $\mathrm{d}\boldsymbol{\xi}=\mathrm{diag}(0,0,0,\mathrm{d}\xi_{a i},\mathrm{d}\xi_{b i},\mathrm{d}\xi_{c i},0)$. We verify the effectiveness of our stochastic phase reduction by calculating the similarity between $R_{\text{rec}}$ and $R_{\text{orig}}$ under different $Q$ and $\sigma$.
	
	Finally, we address the open question regarding the scaling form of the synchronisation of E. coli cell population under quorum sensing coupling. By directly illustrating the phase response function $\Gamma_Q$ for different $Q$ values and comparing them with the sinusoidal function, we find significant discrepancy between the E. coli dynamics and the classic Kuramoto model~\cite{acebron2005kuramoto}. Since the synchronisation behaviour is naturally reflected in the phase difference of the oscillations, we further investigate the dynamical difference of phase in Eq.~\eqref{eq coupled E.coli phase}. For simplicity, we consider $N=2$, and the dynamics of phase difference $ \delta\phi=\phi_{2}-\phi_{1}$ is
	\begin{equation}
		\begin{aligned}
			\dot{\delta\phi}&=\dfrac{\eta Q}{2}\left[\Gamma_Q(-\delta\phi)-\Gamma_Q(\delta\phi)\right]\triangleq\dfrac{\eta Q}{2}\delta\Gamma_Q(\delta\phi)
		\end{aligned}
	\end{equation}
	
	We focus on the scaling form of the response function $\delta\Gamma_Q$ along $Q$. Specifically, we define the `natural strength' as
	\begin{equation}
		\text{NS}(Q)=\dfrac{Q}{2}\left[\max(\delta\Gamma_Q)-\min(\delta\Gamma_Q)\right].
	\end{equation}
	We note that for the Kuramoto model,
	\begin{equation}
		\dot{\phi}_i=w+\dfrac{Q}{N}\sum_{i=1}^N\sin(\phi_j-\phi_i),
	\end{equation}
	the natural strength $\text{NS}(Q)=Q$ linearly scales on coupling strength.

	To the best of our knowledge, the above results represent the first successful analytical and numerical investigation into the mechanism of quorum sensing underlying the onset of synchronisation in multicellular clocks.  They explicitly elaborate how the cell density controls the collective behaviour of a diverse cell population from the perspective of phase dynamics. We anticipate that  the proposed AI framework, along with its application to E. coli dynamics, becomes a standard paradigm for uncovering the principles of biological rhythms whose {\color{black}underlying} mechanisms remain unresolved.
	
	\section{Details of modulating the genetic circuit}\label{appendix gene}
	
	In the synthetic circuit, the dynamical process of the RNA concentration and the feedback coupling effect among the three components are modelled by the dynamic resistor-capacitor (RC) circuits and the repressor circuits respectively. The dynamic RC circuit consists of a resistor and a capacitor which model the degradation process and the accumulation process of the components respectively \cite{teo2020merging}. The voltage across the capacitor is the modelled concentration of one component. The capacitance ($C_i$) modulates the speed of accumulation, which may vary among the three components. The resistance ($R_c$, same for all three components) modulates the speed of degradation. The outflux of electrical current from the dynamic RC circuits cause the voltages across the capacitors to drop, while the influx of electrical current from repressor circuits lead to the rise of the voltage.
	
	The repressor circuits model the Hill function in the accumulation process of components. It consists of three levels of amplification circuits~\cite{hellen2011electronic,hellen2017electronic}. The first level is an inverting amplifier, where the gain is (ideally): 
	$$ G_1 = -\frac{R_2}{R_1}$$
	The output voltage is:
	$$V_I = G_1 (V_{i-1}-V_{th}) = G_1\Delta V$$
	
	The second level is also an inverting amplifier, where the gain is different between positive and negative input voltages. The gain is (ideally, when not saturated):
	$$ G_2=\left\{ 
	\begin{array}{rcl} 
		G_{+2} = -\frac{R_5}{R_3}, & & {V_I \geq 0},\\\\
		G_{-2} = -\frac{\frac{R_5R_4}{R_5+R_4}}{R_3}, & & {V_I < 0}. \end{array} \right. $$
	Considering the saturation effect of the operational amplifiers, the resulting voltage $V_{II}$ is:
	
	$$ V_{II}=\left\{ 
	\begin{array}{rcl} 
		V_{+sat}, & & {G_1G_{+2}\Delta V \geq V_{+sat}},\\ 
		G_1G_{+2}\Delta V, & & {0 \leq G_1G_{+2}\Delta V < V_{+sat}},\\ 
		G_1G_{-2}\Delta V, & & {V_{-sat} \leq G_1G_{-2}\Delta V < 0},\\
		V_{-sat}, & & {G_1G_{-2}\Delta V < V_{-sat}}. \end{array} \right. $$
	
	Here both $U_1$ and $U_2$ are dual op-amp LF412 powered with $\pm 5 V$. For the first level and second level of inverting amplifier we build here, we have:
	$$V_{+sat} = +4.5 V, V_{-sat} = -3.5 V$$
	
	The third level acts as a current source to charge the following capacitor $C_i$. The pair of $R_{b1}$ and $R_{b2}$ combined is a voltage divider. Therefore, the current through $R_E$ is:
	$$ I_t = \frac{s(5-V_{II})-V_{eb}}{R_E} $$
	Here, $s = \frac{R_{b1}}{R_{b1}+R_{b2}}$, $V_{eb}$ is the emitter-base voltage. The transistor we use here is 2N3906, whose $V_{eb}$ varies from 0.5V to the maximum of approximately $V_{ebmx}=0.7V$. $I_t$ reaches maximum at:
	$$I_{max} = \frac{s(5-V_{-sat})-V_{ebmx}}{R_E}$$
	For the chosen components here, $I_{max} = 2.95 mA$ ($V_{-sat} = -3.5 V, V_{ebmx}=0.7V, s=0.42/2.62=0.16, R_E = 222 \Omega$)
	
	The current $I_t$ charges the capacitor of the down-stream oscillator:
	$$R_cC_0\frac{dx_i}{dt} = \frac{C_0}{C_i}(-x_i+\frac{I_tR_c}{V_{th}})$$
	Here $x_i = \frac{V_i}{V_{th}}$ is the normalized voltage (representing the dimensionless RNA concentration), $R_cC_0$ is the time-scale to normalize the time variable. By defining the dimensionless time and parameters:
	$$ t \gets \frac{t}{R_cC_0}, \beta_i = \frac{C_0}{C_i}, \alpha = \frac{I_{max}R_c}{V_{th}}$$
	We consider to model the dynamics of RNA with Hill function inhibition:
	$$ \frac{dx_i}{dt} = \beta_i(-x_i+\frac{\alpha}{1+x_{i-1}^n})$$
	
	Which depends on:
	$$\frac{1}{1+x_{i-1}^n} \gets \frac{I_t}{I_{max}}$$
	
	We consider that the two derivatives of the two functions to be equal at $x=1$, which gets:
	$$n\alpha = \frac{4sR_cG_1G_{-2}}{R_E}$$
	Therefore, we can model the synthetic oscillator with the electrical circuit here, with parameters in the synthetic oscillator to be determined by components in the circuit. 
	
	In the following, we attribute oscillator 1 to the gene lacI, 2 to the tetR, 3 to the cI.  To modulate the oscillations of the three components, we consider two types of control strategies. Noticing that the voltages ($V_1, V_2, V_3$) represent the corresponding component’s state variables, and that the current through the capacitors ($I_1, I_2, I_3$) represent the velocity, we consider to extract the two as time series by measuring them with voltmeters and ammeters respectively. 
	
	In the first type of static control, we consider to modulate static electrical circuit devices to tune the oscillation frequencies. Here we modulate the value of capacitance of oscillator 1 ($C_1$) from 476.19 nF to 526.12 nF, corresponding to the dimensionless transcription rate from 2.1 to 1.9. The benchmark is set as capacitance $C_1$ = 500 nF, being consistent with the transcription rate 2.0. We calculate the resulted frequency shift from the recorded time series data. By feeding the benchmark data into our framework, we can obtain the PRC $\rmZ(\phi)$ of the genetic circuit, and the phase modulation function as:
	\begin{equation*}
		\begin{aligned}
			\delta\rmF(\boldsymbol{\chi}(\phi))&=\frac{1000\text{nF}}{C1}(-a+\frac{\alpha}{c^n}),\\
			\alpha&=\frac{2950}{60},~n=\frac{2.88\times5.2\times22}{\alpha},\\
			(a,b,c)& = \boldsymbol{\chi}(\phi).
		\end{aligned}
	\end{equation*}
	
	According to our Methods, the predicted frequency shift can be analytically determined as
	\begin{equation*}
		\delta w=\frac{1}{2\pi}\int_0^{2\pi}\rmZ(\phi)\delta\rmF(\boldsymbol{\chi}(\phi))\mathrm{d}\phi.
	\end{equation*}

	In the second type of dynamic control, we introduce an external dynamic current source connecting to the RC circuit of oscillator 1, injecting periodic sinusoidal current to the capacitor. We first calculate the frequency $w$ of the aforementioned benchmark data, then we 
	choose the sinusoidal current stimulation with amplitude 6$\mu$V and essential frequency $\Omega$ ranging from $w~ -0.07$ to $w~+0.07$. The conversion relationship between the dimensionless frequency 
	$\Omega$ and the real frequency $\Omega_{\text{TINA}}$ is
	\begin{equation*}
		\Omega_{\text{TINA}} = 1000  \frac{\Omega}{2\pi}.
	\end{equation*}
	Then we analyse the phase-lock and phase-shift phenomena of the genetic circuit under different external frequencies $\Omega$ using the theoretical framework in Section \textbf{Analytical phase reduction under perturbations}.

	All the electrical circuit simulations are completed with the Tina-TI software of Texas Instrument. The circuit documents for reproducing the results will be available once the paper is published.
	
	{\color{black}
	\section{Non-zero optimal equant non-uniformity in FitzHugh-Nagumo dynamics}\label{appendix proof}
	In this Section, we theoretically analyse the existence of zero solution of $\hat{K}$ for the FitzHugh-Nagumo dynamics.
	We begin by reviewing the previous definition of $\hat{K}$,
	
	\begin{equation*}
		\begin{aligned}
		&\hat{K}(\rvx)=\frac{K(\rvx)}{\bar{A}^2(\rvx)},\\
		&\bar{A}(\rvx)
		=
		\frac{1}{T}\int_0^{T}
		\frac{1}{2}
		\left\|
		\big(\rvx-\rvx(t)\big)\wedge \rmF(\rvx(t))
		\right\|
		\mathrm{d}t,\\
		&	K(\rvx)
		=
		\frac{1}{T}\int_0^{T}
		\left|
		\frac{1}{2}
		\left\|
		\big(\rvx-\rvx(t)\big)\wedge \rmF(\rvx(t))
		\right\|
		-
		\bar{A}(\rvx)
		\right|^2
		\mathrm{d}t.
	\end{aligned}
	\end{equation*}
	
	To prove that there is no equant satisfying exact areal uniformity as $\hat{K}=0$, we argue by contradiction. 
	Suppose that there exists a point $\rvx^\ast=(u,v)$ such that $K(\rvx^\ast)=0$. 
	Then the swept areal velocity$
		\frac{1}{2}
		\left\|
		\big(\rvx^\ast-\rvx(t)\big)\wedge \rmF(\rvx(t))
		\right\|$
	is constant along the whole limit cycle.
	Since the limit cycle is a regular closed orbit, $\mathcal{A}_{\rvx^\ast}(t)$ is continuous in $t$. 
	Moreover, if $\mathcal{A}_{\rvx^\ast}(t)$ is a positive constant, the signed quantity $
		q(t)
		=
		\big(\rvx^\ast-\rvx(t)\big)\wedge \rmF(\rvx(t))$
	cannot change sign along the connected orbit. Therefore, $q(t)$ itself must be constant up to a fixed sign. 
	Taking the derivative with respect to $t$, we obtain
	\begin{equation}
		0
		=
		\frac{\mathrm{d}}{\mathrm{d}t}
		\left[
		\big(\rvx^\ast-\rvx(t)\big)\wedge \rmF(\rvx(t))
		\right].
	\end{equation}
	Using $\dot{\rvx}(t)=\rmF(\rvx(t))$, this gives
	\begin{equation}
		0
		=
		\big(\rvx^\ast-\rvx(t)\big)\wedge \ddot{\rvx}(t),
		\label{eq:necessary_acceleration_condition}
	\end{equation}
	because $\dot{\rvx}(t)\wedge \dot{\rvx}(t)=0$. 
	Equivalently, for every point on the limit cycle, the vector from $\rvx(t)$ to the candidate equant $\rvx^\ast$ must be parallel to the acceleration vector $\ddot{\rvx}(t)$. 
	This is a strong necessary condition for the existence of an equant satisfying exact areal uniformity.
	
	For the FitzHugh-Nagumo system, we have
	\begin{equation*}
		\begin{aligned}
			&f(x,y)=x-\frac{x^3}{3}-y+I,
			~
			g(x,y)=\varepsilon(x+a-by),\\
			&I=0.0,a=0.7,b=0.2,\varepsilon=0.05
		\end{aligned}
	\end{equation*}
	Then
	\begin{equation*}
		\rmF(x,y)
		=
		\begin{pmatrix}
			f(x,y)\\
			g(x,y)
		\end{pmatrix},
	\end{equation*}
	and
	\begin{equation*}
		\ddot{\rvx}
		=
		D\rmF(\rvx)\rmF(\rvx)
		=
		\begin{pmatrix}
			(1-x^2)f-g\\
			\varepsilon f-\varepsilon b g
		\end{pmatrix}.
	\end{equation*}
	
	Now consider the points on the limit cycle at which $\dot{x}=f(x,y)=0$. 
	These points correspond to the extrema of the $x$-coordinate along the periodic orbit. 
	For the relaxation oscillation considered here, the limit cycle intersects the cubic nullcline $
		y=x-\frac{x^3}{3}
$
	at two distinct points, denoted by
	\begin{equation*}
		\rvx_1=(x_1,y_1),
		~
		\rvx_2=(x_2,y_2),
		~
		x_1\neq x_2 .
	\end{equation*}
	At these points, $f(x_i,y_i)=0$. 
	Since they are not equilibrium points, we have $g(x_i,y_i)\neq 0$. 
	Therefore,
	\begin{equation*}
		\ddot{\rvx}_i
		=
		\begin{pmatrix}
			-g(x_i,y_i)\\
			-\varepsilon b g(x_i,y_i)
		\end{pmatrix}
		=
		-g(x_i,y_i)
		\begin{pmatrix}
			1\\
			\varepsilon b
		\end{pmatrix},
		~ i=1,2.
	\end{equation*}
	Hence, the acceleration vectors at both points are parallel to the same direction $(1,\varepsilon b)^\top$.
	
	By the necessary condition in Eq.~\eqref{eq:necessary_acceleration_condition}, the candidate equant $\rvx^\ast=(u,v)$ must lie on the straight line through each $\rvx_i$ with direction $(1,\varepsilon b)^\top$. 
	Thus,
	\begin{equation*}
		v-y_i
		=
		\varepsilon b (u-x_i),
		~ i=1,2.
	\end{equation*}
	Subtracting the two equations gives the necessary condition
	\begin{equation}
		y_2-y_1
		=
		\varepsilon b (x_2-x_1).
		\label{eq:parallel_line_condition}
	\end{equation}
	However, since both points lie on the cubic nullcline $y=x-x^3/3$, we have
	\begin{equation*}
		y_2-y_1
		=
		\left(x_2-\frac{x_2^3}{3}\right)
		-
		\left(x_1-\frac{x_1^3}{3}\right).
	\end{equation*}
	Therefore,
	\begin{equation}
		\frac{y_2-y_1}{x_2-x_1}
		=
		1-\frac{x_1^2+x_1x_2+x_2^2}{3}.
		\label{eq:cubic_secant_slope}
	\end{equation}
	For the parameter setting considered in this work,
	\begin{equation*}
		\varepsilon b = 0.05\times 0.2 = 0.01.
	\end{equation*}
	In contrast, the two intersections of the limit cycle with the cubic nullcline lie on the two slow branches of the relaxation oscillation. 
	For the trajectory used in our numerical experiments, these two points are approximately
	\begin{equation*}
		(x_1,y_1)\approx(-2.07,0.88),
		~
		(x_2,y_2)\approx(1.91,-0.41),
	\end{equation*}
	which gives
	\begin{equation*}
		\frac{y_2-y_1}{x_2-x_1}
		\approx -0.32
		\neq
		0.01
		=
		\varepsilon b.
	\end{equation*}
	Thus, Eq.~\eqref{eq:parallel_line_condition} cannot hold. 
	Consequently, there is no point $\rvx^\ast$ that satisfies the necessary acceleration condition for constant swept areal velocity along the whole limit cycle.
	
	This contradiction shows that no exact  equant exists for the FitzHugh-Nagumo dynamics under the parameters considered here. 
	Therefore,
	\begin{equation}
		K(\rvx)>0,
		~
		\text{for all admissible } \rvx,
	\end{equation}
	and hence
	\begin{equation}
		\hat{K}^\ast>0 .
	\end{equation}

}
%
	\section{Validating the framework with the FitzHugh-Nagumo system}\label{appendix FHN}
	
	In this Section, we systematically validate the effectiveness of the proposed phase reduction framework in the FitzHugh-Nagumo system, especially with respect to phase synchronisation. Phase synchronisation phenomenon of rhythmic oscillations have been extensively studied in the context of collective dynamics from all scales and domains \cite{rosenblum1996phase,rosenblum2001phase,strogatz2018nonlinear,gray1994synchronous,varela2001brainweb,sauseng2008does,palva2005phase}. Finding the phase dynamics of the original complex system is the most natural way to analyse the phase synchronisation of the oscillators. However, among the plentiful natural oscillations, only a few have  their analytical expressions of limit
	cycles and phase dynamics, such as the Stuart-Landau oscillator \cite{zhang2022neural}. 
	
	Here we consider the FitzHugh-Nagumo system as follows:
	\begin{equation}
		\begin{aligned}
			\dot{x}&=x-x^3/3-y+I,\\
			\dot{y}&=\varepsilon(x+a-by),
		\end{aligned}
	\end{equation}
	which describes the physiological phenomena associated with the variation in time of the electrical potential of the neural membrane~\cite{fitzhugh1969mathematical}. Here, $u$,~$v$ respectively represent the electrical potential of the cell membrane and the auxiliary variable depending on the refractory period, $I$ stands for the injected current, $\varepsilon\ll1$ is the time constant leading to the slow-fast dynamics, and the parameters $a$ and $b$ relate to the number of channels of the cell membrane that are opened to the $Na^{+}$ and $K^{+}$ ions. We set $I=0.0$, $a=0.7$, $b=0.2$, $\varepsilon=0.05$, and denote the dynamics by $\dot{\rvx}=\rmF(\rvx)$ with $\rvx=(x,y)$. Then we consider the phase synchronisation against the perturbations using the same notations as those in Section~\ref{Methods phase reduction}.
	\begin{itemize}
		\item \textbf{Network interaction.} The network interaction is described as $\sigma\sum_{j=1}^NA_{ij}\rmH(\rvx_i,\rvx_j)$.
		For the network structure $\rmA$, we select five realistic networks including food-web networks (Rhode $N=20$, gramwet $N=69$, baydry $N=128$, and SmallW $N=396$) and citation network (SmaGri $N=1059$). The data are collected from {\texttt{Pajek Dataset}}. For the coupling function, we consider three practical interactions:
		the electrical coupling {\color{black} $\rmH=(x_j-x_i,0)^\top$}~\cite{li2007effects}, the memristor-based coupling $\rmH=(\left[k_1+k_2(y_2-y_1)^2\right](x_2-x_1),0)^\top$~\cite{volos2015memristor}, and the chemical synaptic coupling $\rmH=\left(\frac{1}{1+\exp(k(\alpha-\theta_i))+\exp(k(\theta_i-\beta))},0\right)^\top$, $\theta_i=\arg\tan\frac{y_i}{x_i}$~\cite{korotkov2019dynamics}. Figure~5A shows that the phase response functions of all the three coupling functions are significantly different from the classic Kuramoto model. The accuracy of the reduced phase dynamics is verified and shown in Fig.~5C(top), wherein the original and the reconstructed trajectory {\color{black}$x_i(t)$} of one randomly picked oscillator $i$ perfectly coincide. To assess the effect of different couplings to synchronisation, we use time series {\color{black}$\{y_i(t)\}_{i=1}^N$} in the original (resp. reconstructed) trajectories to define the indicator $R_{\text{orig}}$ (resp. $R_{\text{rec}}$) in the same way as in Eq.~\eqref{eq ind1}. We vary the coupling strength $\sigma$ in $[0,1]$ and calculate the synchronisation indicators under different networks. Figure~5B shows the reconstructed indicator $R_{\text{rec}}$ under electrical and memristor coupling, showing that they coincide well with the original indicator $R_{\text{orig}}$ among all network structures, while  $R_{\text{rec}}$ under chemical coupling deviates  $R_{\text{orig}}$ as $\sigma$ grows. The reason lies in the synchronisation non-vanishing property of the chemical coupling, i.e., the coupling term is equal to zero when oscillators are synchronised as $\rvx_i=\rvx_1$ for $i=1,...,N$. In this case, the non-zero coupling term changes the limit cycle of the oscillation.
		This property is also reflected in the phase response functions, where $\Gamma(0)=0$ for electrical and memristor coupling and $\Gamma(0)\neq0$ for chemical coupling. Therefore, for synchronisation vanishing coupling, our framework works well regardless of the coupling strength. For synchronisation non-vanishing coupling function, we require a weak coupling strength to restrict the perturbed oscillation in the vicinity of the original limit cycle. In addition, we consider the order parameter $R_{\text{ord}}$ defined in Eq.~\eqref{eq ind2} for the phase dynamics. As shown in Fig.~5C(bottom), the order parameter and the synchronisation indicator of the original trajectories present the same tipping point versus coupling strength, confirming the feasibility of our framework.
		
		\item \textbf{Stochastic force.} For a cluster of oscillators driven by both the common noise and the independent noise described in Eq.~\eqref{eq common+independ}, we respectively consider the influence of the additive noise $\rmG=\rmH=(\sigma,\sigma)^\top$ and the multiplicative noise $\rmG=\rmH=(\sigma x,\sigma y)^\top$. We numerically simulate the trajectories of 5000 oscillators over the time interval $[0, 20]$ using the Euler-Maruyama method~\cite{kloeden1992stochastic}. Since the oscillators are decoupled, we calculate the phase difference $\delta\phi_i(t) = \phi_i(t) - \phi_1(t),~i=1,...,5000,$ for each oscillator with respect to the first oscillator. The theoretical stationary distribution of the phase difference $\delta\phi_i$ is given in Eq.~\eqref{eq FK common+independ}.
		To approximate the stationary distribution of $\delta\phi_i$, we simulate the system 100 times, using random seeds from 0 to 99. We then calculate the stationary state by taking the time average over the interval $[19.6, 20]$, i.e., $\delta\phi_i^\ast = \langle \delta\phi_i(t)\rangle_{[19.6,20]}$.
		Next, we validate the reduced stochastic phase dynamics by comparing the theoretical predictions with the numerical approximations. As shown in Fig.~\ref{figFHN}D, our theoretical predictions match the numerical results very well under both additive and multiplicative noise.

		\item \textbf{Periodic force.} We add the periodic force $f(t)$ to the first variable $x$, and we consider sine wave $f=\sin(\Omega t)$ and square wave {\color{black} $f=\text{Heaviside}(M(\Omega t))$, $M(\Omega t)=\Omega t ~\text{mod}~2\pi\in[-\pi,\pi)$}. We set $\sigma=0.05$ to keep the oscillation near the original limit cycle. According to Eq.~\eqref{eq phase period}, we obtain the dynamics of the phase difference $\delta\phi=\phi-\Omega t$ as
		$\dot{\delta\phi}=\Delta+\Gamma(\delta \phi)$. The asymptotic behaviour is determined by the relation between $\Delta$ and $\Gamma$, as shown in Fig.~\ref{figFHN}E. We find a saddle-node bifurcation for dynamics under both periodic force. As the $\Delta$ grows, the dynamics first appears as a saddle-node, then this saddle-node evolves to a stable equilibrium (solid dot) and an unstable equilibrium (hollow dot), and finally they merge back to a saddle-node before the equilibrium disappears. 
		The phase dynamics has phase-lock phenomenon when there exist a stable equilibrium (blue region), i.e., the frequency of oscillation converges to external frequency $\lim_{t\to\infty}\dot{\phi}(t)=\Omega$.  Otherwise, the phase frequency diverges and is characterised as a phase-slip (grey region). To validate the partition regions of phase-lock and phase-slip, we uniformly sample various $\Delta$ in each region and simulate the corresponding original trajectories $\{\rvx(t):0\le t\le1500\}$. We obtain the phase trajectories using the learned auto-encoder as $\phi(t)=\rvh(\text{Enc}(\rvx(t)))$. Then the phase difference is obtained as $\delta\phi(t)=\phi(t)-\Omega t~(\text{mod} 2\pi)$. We calculate the temporal variance of $\delta\phi(t)$ in the last $15$ time unit to see whether the phase locks or slips. The numerical results in Fig.~\ref{figFHN}E fit with the theoretical prediction well.
	\end{itemize}

	\clearpage
\onecolumngrid


\section*{Supplemental Information}

\setcounter{section}{0}
\setcounter{subsection}{0}
\setcounter{equation}{0}
\setcounter{figure}{0}
\setcounter{table}{0}

\renewcommand{\thesection}{S\arabic{section}}
\renewcommand{\thesubsection}
{S\arabic{section}.\arabic{subsection}}
\renewcommand{\theequation}
{S\arabic{section}.\arabic{equation}}
\renewcommand{\thefigure}{S\arabic{figure}}
\renewcommand{\thetable}{S\arabic{table}}
\section*{Contents}
\startcontents[supp]
\printcontents[supp]{supp}{1}[2]{
	\contentsmargin{2.5em}
}



\newcommand{\revtex}{REV\TeX\ }
\newcommand{\classoption}[1]{\texttt{#1}}
\setlength{\textheight}{9.5in}
\setcounter{section}{0}
\renewcommand{\thesection}{S{\arabic{section}}}
\renewcommand{\thesubsection}{\arabic{section}.\arabic{subsection}}
\renewcommand{\thesubfigure}{(\alph{subfigure})}
\renewcommand{\theequation}{S\arabic{section}.\arabic{equation}}
\renewcommand\thefigure{S\arabic{figure}}
\renewcommand\thetable{S\arabic{table}}
\renewcommand{\thefootnote}{\fnsymbol{footnote}}
\renewcommand{\baselinestretch}{1.5}
\renewcommand{\bm}[1]{\mbox{\boldmath{$#1$}}}
\renewcommand{\thesubfigure}{(\alph{subfigure})}
\linespread{1.05}

\newcommand{\rvu}{\mathbf{u}}
\newcommand{\rva}{\mathbf{a}}

\newpage

\section{A brief overview of Methods}

\begin{figure}[htp]
		\centering
		\includegraphics[width=16.5cm]{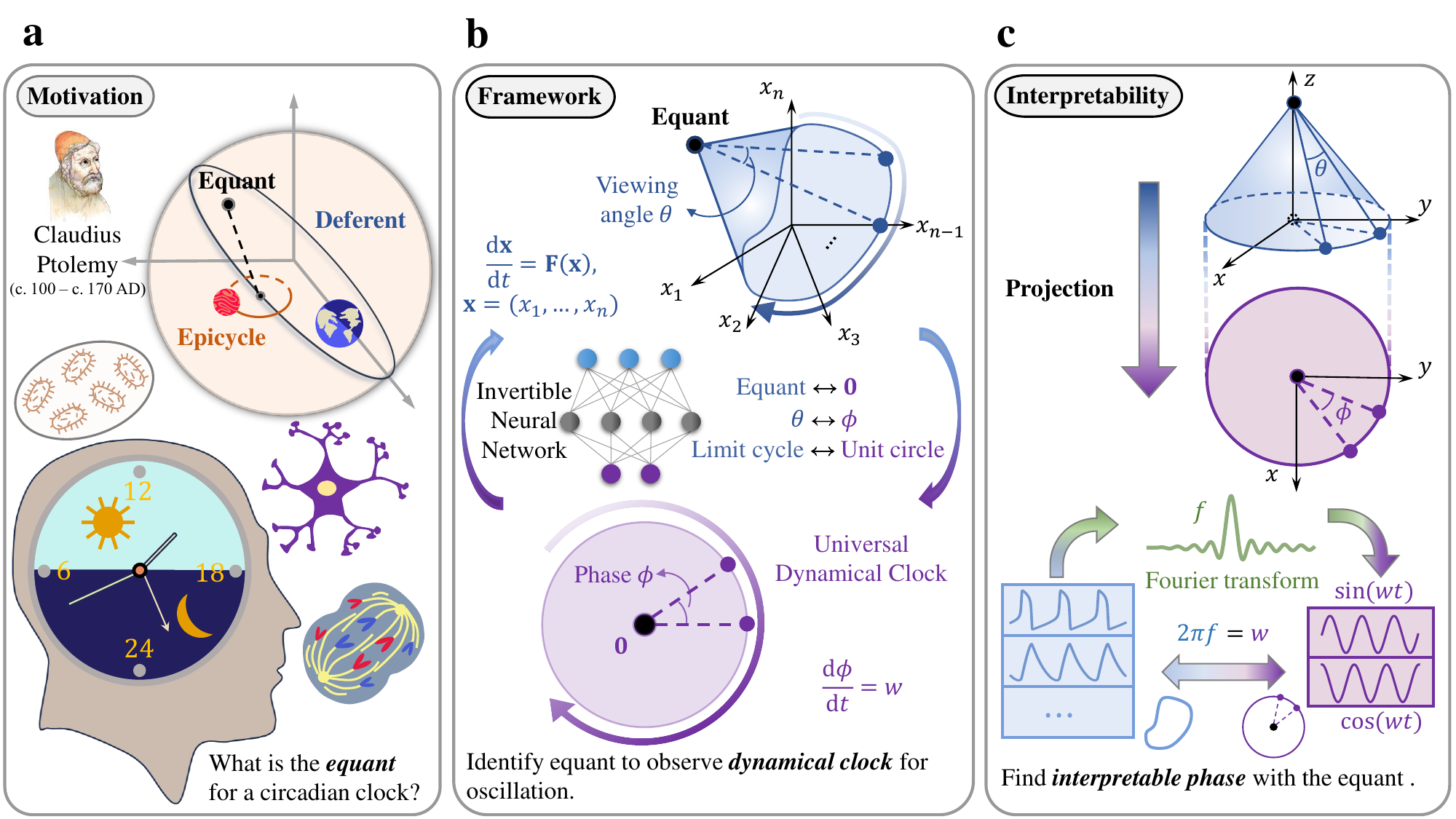}
		\caption{\textbf{Main ideas of the machine learning framework}. \textbf{a}, Ptolemy developed the concept of equant to explain that the path of heavenly body
 is uniform about the equant and circular around the Earth. Does there exist the equant for circadian clock,
 e.g., genetic clock, neuronal rhythm and mitosis?\textbf{b}, We specify the universal dynamical clock over the unit circle and equate it with a general oscillatory system through machine learning. \textbf{c}, (Top) Projecting the viewing angle $\theta$ observed from equant on z-axis to the phase $\phi$. (Bottom) The dominant Fourier frequency
 of oscillator’stemporal trajectory closely relates to the natural frequency of the universal dynamical clock.}
		\label{figS1}
	\end{figure}

\textbf{Identifying the Equant.} An observer situated at the equant is expected to perceive a uniform motion of Mars as it orbits the Earth, according to Ptolemy's theory. Based on the same idea, we aim at identifying an equivalent equant from which an imaginary observer would obtain the most uniform clock representation of a general oscillator, such as the circadian clock for biorhythms illustrated in Fig.~\ref{figS1}A. In this way, the phase of the oscillator gains physical significance from the perspective of the equant, defined as the interpretable phase.
	
	Specifically, we consider an oscillatory dynamics $\mathrm{d}{\rvx}/\mathrm{d}t=\rmF(\rvx)$ with $\rvx\in\mathbb{R}^n$ and 
	a stable periodic solution as a limit cycle. Generally, the limit cycle is an irregular closed orbit in the high-dimensional space without an analytical expression, and the phase $\phi$ is a scalar variable whose value represents a unique point on the limit cycle. As the oscillator moves along the limit cycle for one complete period, the associated phase dynamically increases by $2\pi$. Our goal is to identify the equant and establish the equivalent relationship between the phase $\phi$ and the viewing angle observed from the equant to the limit cycle. We require that the phase should vary uniformly with time, i.e., $\mathrm{d}{\phi}/\mathrm{d}t=w$
	where $w$ denotes the natural frequency. 
	
	Before solving this mathematical problem for general oscillatory dynamics, we consider a universal dynamical clock in which the oscillator, i.e., the pointer, rotates uniformly along the unit circle (see Fig.~\ref{figS1}B bottom). For the universal dynamical clock, the centre of the circle is the equant and the arc angle naturally becomes the phase. Beginning from the universal dynamical clock, we identify the equant for general oscillation through the following two-step approach.
	
	Firstly, we find a function that uniquely associates each point on the limit cycle to a point on the unit circle. Analytically constructing such a function is intractable due to the irregularity and asymmetry of the limit cycle. To circumvent this difficulty, we devise a parametrised neural network on circular manifold to learn the function. This step aims to equate the limit cycle of the oscillatory system with the unit circle of the universal dynamical clock.
	
	Secondly, we use the equant to define a surface having the limit cycle as its boundary, such that any straight line connecting the equant to the limit cycle lies within this surface. For the universal dynamical clock, such a surface is naturally the disk with the unit circle as its boundary. In the first step we find a function that equates the boundaries of the surface defined by limit cycle and the disk. In the second step, we extend the function from the boundary to the surface such that the extended function uniquely maps points on this surface to points on the disk. Among all candidate extensions, we require the optimal extension to map line segments emanating from the unknown equant to line segments emanating from the centre of the circle. Therefore, the unknown equant is mapped by the optimal extended function to the centre of the unit circle, yielding a connecting relationship for identifying the equant. To find the optimal extension, we employ the invertible neural network~\cite{dinh2016density}, a machine learning method that naturally guarantees the invertible structure of the extended function.  A high-level summary of our two-step approach is provided in Fig.~\ref{figS1}B.\\
	
	\noindent\textbf{Interpreting phase from equant}.
	Merely identifying the equant is insufficient to complete the depiction of an interpretable phase, since the viewing angle $\theta$ directly defined by the intersection angle between line segments from the equant to the limit cycle potentially possesses a period different from $2\pi$ and does not represent the phase $\phi$. 
	
	To provide a meaningful definition of phase from the viewpoint of the equant, we first consider the following dynamics in a three-dimensional space $(x,y,z)$, in which the oscillator uniformly rotates about a circle on the $(x,y)$-plane (see Eq.~\eqref{eq circle}). All the points on the $z$-axis satisfy the definition of the equant, but only the viewing angle observed from the origin possesses period $2\pi$ and matches the phase. Nonetheless, we can simply transform the viewing angle defined by any equant into the phase by multiplying it by a constant. This is possible because there exists a linear projection from the viewing angle of any equant on the $z$-axis to the phase, as shown in Fig.~\ref{figS1}C (top). 
	
	For general oscillatory dynamics with irregular and asymmetrical limit cycle, the map from the viewing angle to the phase is essentially non-linear and intractable analytically. To address such a challenge, we recall that the extended function obtained from our aforementioned two-step strategy maps line segments emanating from the equant to
	line segments emanating from the centre of the disk. Therefore, the extended function maps the viewing angle from the equant to the phase and endows the phase with physical meaning by equating it with the viewing angle observed from the equant.

	We then shift our focus to the dynamic property of the interpretable phase: the frequency. In time series analysis, the dominant Fourier frequency $f$ of a periodic temporal trajectory is the reciprocal of the period~\cite{hamilton2020time}. Since we equate the oscillation on the limit cycle with the universal dynamical clock, the universal dynamical clock must have the same period as the oscillatory motion. We note that the phase frequency $w$ is the angular speed of the universal dynamical clock about the unit circle. Thus, the period of the universal dynamical clock is ${2\pi}/{|w|}$. Therefore, we conclude that the frequency derived from the interpretable phase is precisely a rescaling of the dominant Fourier frequency
	as $|w|=2\pi f$. Figure~\ref{figS1}C (bottom) illustrates the relationship between the natural frequency and the Fourier frequency.\\

\begin{figure}[htp]
		\centering
		\includegraphics[width=16.5cm]{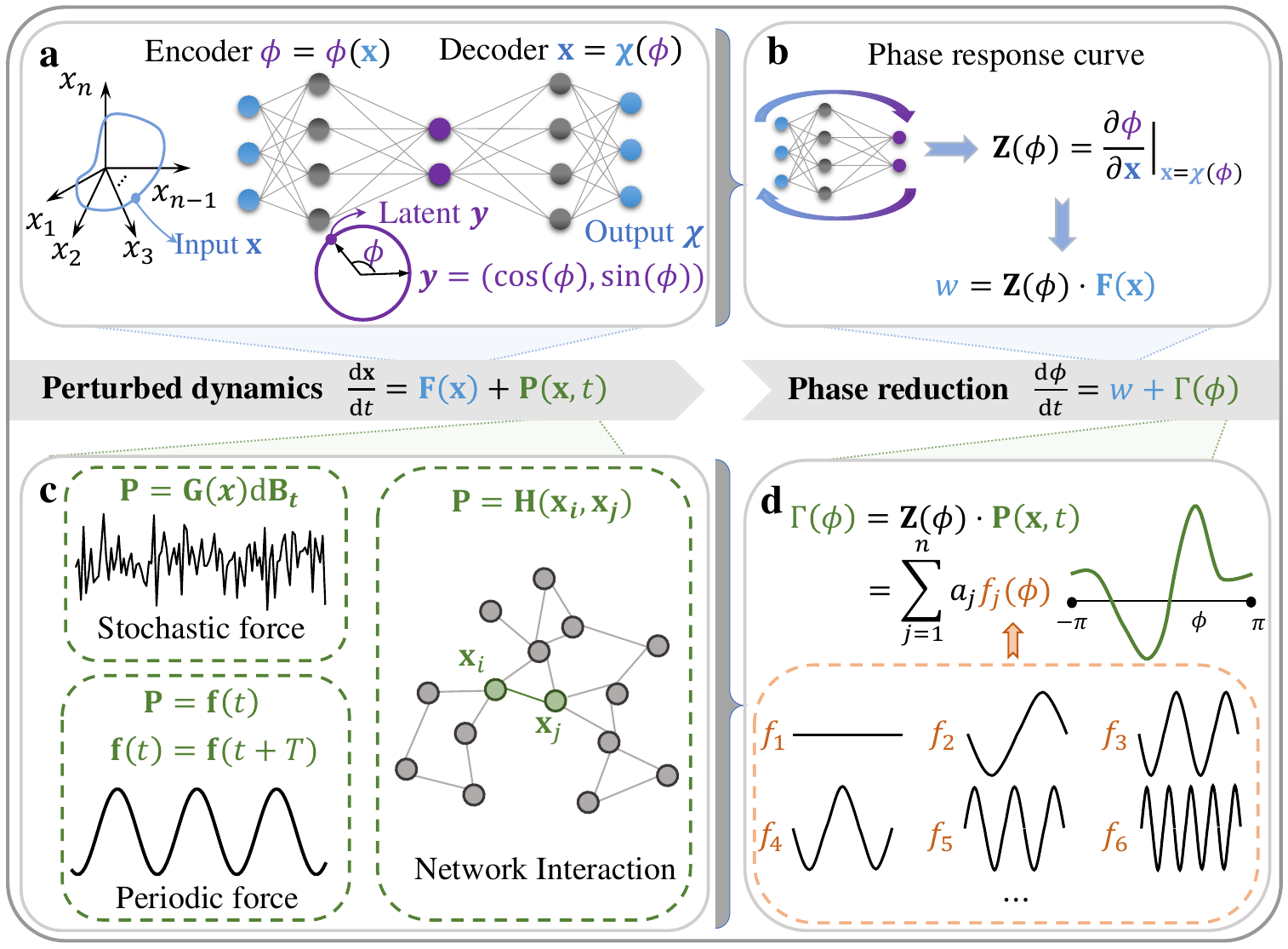}
		\caption{\textbf{Illustration of analytical phase reduction.} \textbf{a}, The neural network equates the limit cycle and the unit circle.
 \textbf{b},The phase response curve and the natural frequency is deduced from the neural network. \textbf{c}, The oscillations
 are perturbed, respectively, by the noise, periodic force, or network interaction. \textbf{d}, The phase response function to the perturbation is expanded with analytical function basis.}
		\label{figS2}
	\end{figure}

    
	\noindent\textbf{Analytical phase reduction}.
	Understanding the precise effects of external forces on oscillations is critically important in  life sciences, such as medical treatment design~\cite{altinok2007cell,levi2008implications} and circadian biology~\cite{goldbeter2022multi,leloup2003toward}.
	We delve into identifying the analytical formulation of phase dynamics for complex oscillatory systems perturbed by external forces, including the interplay among oscillators connected in a network~\cite{pecora1998master}, the stochastic force driven by intrinsic or extrinsic noise~\cite{gonze2002robustness}, and the periodic force from existing rhythms in nature~\cite{yan2019robust}, as shown in Fig.~\ref{figS2}C. 
	
	We aim at reducing perturbed oscillatory dynamics ${\mathrm{d}\rvx}/{\mathrm{d}t}=\rmF(\rvx)+\rmP(\rvx,t)$ to phase dynamics ${\mathrm{d}\phi}/{\mathrm{d}t}=w+\Gamma(\phi,t)$, where $\rmP(\rvx,t)$ represents the perturbation and $\Gamma(\phi)$ is the phase response function to the perturbation. According to the chain derivative rule, the phase response to the perturbation is the inner product of the gradient of phase function and external force, i.e., $\Gamma(\phi,t)=\frac{\partial \phi}{\partial\rvx}\cdot\rmP(\rvx,t)|_{\rvx=\boldsymbol{\chi}(\phi)}$. Here, $\frac{\partial \phi}{\partial \rvx}|_{\rvx=\boldsymbol{\chi}(\phi)}$  is the phase response curve~\cite{monga2019phase}, and ${\boldsymbol{\chi}(\phi)}$ maps the phase defined on the universal dynamical clock to points on the limit cycle .

	Recall that in our two-step approach, we identify the function that relates the unit circle to the limit cycle. We then compose this function with the Euler map, which connects the phase to points on the unit circle, to obtain the phase function. This establishes an equivalent relationship between the phase and the limit cycle. As a result, we derive the phase function $\phi(\rvx)$ and its inverse $\boldsymbol{\chi}(\phi)$. Next, we calculate the phase response curve by directly computing the gradient of the phase function, as shown in Figs.~\ref{figS2}A,B. Finally, by calculating the inner product of the phase response curve with the external force, we successfully reduce the original perturbed system to the analytical phase dynamics.
	To improve the interpretability of our phase reduction framework, we find a simple mathematical expression for the phase response function using a predefined 
	dictionary of basis functions~(see~Fig.~\ref{figS2}D).

    \subsection{Problem setup}
	We consider a globally (or locally) attractive limit cycle oscillator described by a general non-linear dynamical system 
	\begin{equation}\label{eq1}
		\dfrac{\mathrm{d}\rvx}{\mathrm{d}t}=\rmF(\rvx),~\rvx\in\mathbb{R}^n,
	\end{equation} 
	wherein we denote by $\rmC$ the stable limit cycle.
	We assume there exists a periodic scalar function $\phi(\rvx)$ on $\rmC$ with period $2\pi$ such that its time derivative is a constant $w$, i.e.,
	\begin{equation}
		\dfrac{\mathrm{d}\phi(\rvx)}{\mathrm{d}t}=w,~\phi\in[0,2\pi),
	\end{equation}
	where $\phi$ and $w$ are known, respectively, as the phase function and the natural frequency~\cite{monga2019phase}, respectively. Such a phase function establishes a homeomorphism between the limit cycle and the unit circle, and
	we denote by $\boldsymbol{\chi}(\phi)\in\rmC$ the inverse map of phase. According to chain derivation, we have,
	\begin{equation}
		\frac{\mathrm{d}\phi(\rvx)}{\mathrm{d}t}=\frac{\partial \phi}{\partial\rvx}\cdot\rmF(\rvx)=w,
	\end{equation}
	which is used to construct loss function in our machine learning framework. In general systems the limit cycle $\rmC$ is a complex $1$-dimensional manifold without analytical formulation, which causes challenge for identifying $\phi$ and $\boldsymbol{\chi}$. For brevity, we denote $\frac{\mathrm{d}\rvx}{\mathrm{d}t}$ by $\dot{\rvx}$  below.
	
	{\color{black}
	To physically interpret the phase, otherwise a mathematical coordinate, we introduce the equant as an observational point $\rvx^\ast$ from which the observer at the equant finds the oscillator uniformly rotating along $\rmC$. Motivated by Ptolemy’s equant and to avoid the non-uniqueness of purely topological viewing-angle conditions, we select this point through an areal-uniformity criterion inspired by Kepler’s second law: the swept area along the limit cycle should be as uniform as possible.

	\noindent\textbf{Definition 1.1}\textit{
	We define the equant for the system in Eq.~\eqref{eq1} as a point $\rvx^\ast$ such that $\mathrm{(i)}$ for the viewing line $l_{\rvx^\ast\rvx(t)}$ connecting $\rvx^\ast$ and $\rvx(t)\in\rmC$, there
	exists a homeomorphism between the phase $\phi$ and the viewing angle $\theta$ induced by the viewing lines,
	$
	\theta(t, t')
	:=
	\angle\left(l_{\rvx^\ast\rvx(t)},l_{\rvx^\ast\rvx(t')}\right)
	$; and  $\mathrm{(ii)}$ it minimises the normalised areal variance, i.e.,
	\begin{equation}
		\rvx^\ast \in \arg\min_{\rvx} \hat{K}(\rvx),
		~
		\hat{K}(\rvx)=\frac{K(\rvx)}{\bar{A}^2(\rvx)} .
	\end{equation}
	Here, $\hat{K}(\rvx)$ is a dimensionless measure of the non-uniformity of the swept areas observed from $\rvx$. Specifically,
	\begin{equation*}
		\bar{A}(\rvx)
		=
		\frac{1}{T}\int_0^{T}
		\frac{1}{2}
		\left\|
		\big(\rvx-\rvx(t)\big)\wedge \rmF(\rvx(t))
		\right\|
		\mathrm{d}t
	\end{equation*}
	denotes the mean areal velocity over one period, and
	\begin{equation*}
		K(\rvx)
		=
		\frac{1}{T}\int_0^{T}
		\left|
		\frac{1}{2}
		\left\|
		\big(\rvx-\rvx(t)\big)\wedge \rmF(\rvx(t))
		\right\|
		-
		\bar{A}(\rvx)
		\right|^2
		\mathrm{d}t
	\end{equation*}
	denotes the variance of the swept areal velocity over one period.   $\rvx(t)\subset \rmC$ is the periodic solution with period $T$.
	We define the optimal equant non-uniformity as
	$
	\hat{K}^\ast=\min_{\rvx}\hat{K}(\rvx)$.
}

According to the Definition, $\hat{K}$ is non-negative and continuous. Hence, provided that the admissible search domain of $\rvx$ is bounded, there exists an equant that minimises $\hat{K}$. Moreover, if the minimiser set $\arg\min_{\rvx}\hat{K}(\rvx)$ is a singleton, the equant is unique.	

}
    
%

	In order to take perturbations ensuing from the real-world scenarios into consideration, we extend the domain of $\phi$ to the vicinity of $\rmC$. Given any phase value $\psi\in[0,2\pi)$, we define the isochron $\rmI_{\psi}$ as the manifold that intersects $\rmC$ with point $\boldsymbol{\chi}(\psi)$ such that any trajectory initiated from $\rmI_{\psi}$ converges to the same limit, i.e.,
	\begin{equation}
	\lim_{t\to\infty}\Vert\rvx(t;0,\rvx_1)-\rvx(t;0,\rvx_2)\Vert=0,~\forall \rvx_1,\rvx_2\in \rmI_\psi.
	\end{equation}
	The isochron can be regarded as the orthogonal space of the limit cycle, and thus it is a $(n-1)$-dimensional manifold.
	We can naturally extend $\phi$ to the isochrons $\{\rmI_\psi:\psi\in[0,2\pi)\}$ as $\phi(\rvx)=\psi$ for any $\rvx\in\rmI_\psi$~\cite{winfree1980geometry,winfree1974patterns,guckenheimer1975isochrons}. Then we consider the perturbed limit cycle oscillator described as,
	\begin{equation}	\dot{\rvx}=\rmF(\rvx)+\rmP(\rvx,t),~\rvx\in\mathbb{R}^n,t\ge0,
	\end{equation}
	where $\rmP(\rvx,t)$ is the external force representing the perturbation. The corresponding perturbed phase dynamics becomes
	\begin{equation}
		\dot{\phi}=\dfrac{\partial \phi}{\partial\rvx}\cdot\left(\rmF(\rvx)+\rmP(\rvx,t)\right)=w+\dfrac{\partial \phi}{\partial\rvx}\cdot\rmP(\rvx,t).
	\end{equation}
	By denoting the phase response curve as
	\begin{equation}\label{eq PRC}
		\rmZ(\phi)=\frac{\partial \phi}{\partial \rvx}\bigg|_{\rvx=\boldsymbol{\chi}(\phi)},
	\end{equation}
	we obtain the closed-form of phase dynamics on $\rmC$ as
	\begin{equation}\label{eq phase}
		\dot{\phi}=w+\Gamma(\phi,t),~\Gamma(\phi,t)=\rmZ(\phi)\cdot\rmP(\boldsymbol{\chi}(\phi),t).
	\end{equation}
	We require that the strength of perturbation is restricted to a certain range such that the perturbed dynamics fluctuate in the vicinity of $\rmC$, which is a fundamental assumption in phase reduction theory~\cite{brown2004phase,goldobin2010dynamics,NOVICENKO20121090}.
	We consider three major forms of perturbation:  $\text{(i)}$ the stochastic force driven by noise; $\text{(ii)}$ the periodic external force; and $\text{(iii)}$ and the interplay among oscillators $\{\rvx_i\}_{i=1}^N$ connected in a network. The detailed formulations of these perturbations are provided in Section~\ref{Methods phase reduction}.
	
	Given the general oscillatory dynamics in Eq.~\eqref{eq1}, the goal of our framework is threefold,
	\begin{itemize}
		\item [(i)] finding analytical expression of the phase $\phi(\rvx)$ and its inverse $\boldsymbol{\chi}(\phi)$ to obtain the homeomorphism between the general oscillation and the universal dynamical clock,
		\item [(ii)] discovering the equant $\rvx^\ast$ to interpret the physical meaning of the phase, and
		\item [(iii)] investigating the influence of perturbations to oscillation by obtaining the analytical formulation of perturbed phase dynamics in Eq.~\eqref{eq phase}.
	\end{itemize}

	\subsection{First step: Equating the limit cycle with the unit circle}
	In order to obtain the phase function, we note that the limit cycle $\rmC$ is a $1$-dimensional manifold, and the phase $\phi(\rvx)\in[0,2\pi)$ belongs to the quotient space $\mathbb{R}/_\sim$ with an equivalent relation as 
	\begin{equation}
		x\sim y\Leftrightarrow x-y=2k\pi,~k\in\mathbb{Z}.
	\end{equation}
	We seek the phase function $\phi:\rmC\to \mathbb{R}/_\sim$ and its inverse $\boldsymbol{\chi}:\mathbb{R}/_\sim\to\rmC$ via using an auto-encoder neural network.  To integrate the quotient manifold structure into the machine learning framework, we note that $\mathbb{R}/_\sim$ is equivalently related to  the unit circle $\rmS^1=\{(x,y):x^2+y^2=1\}$ via the exponent map $e^{\text{i}\phi}$.  Therefore, we seek the phase function $\phi$ and its inverse $\boldsymbol{\chi}$ via an auto-encoder neural network with the latent space $\rmS^1$ as: 
	\begin{equation}
		\begin{aligned}
			\text{(Encoder)}~&\rvx\in\boldsymbol{C}\to \rvy\in \rmS^1: \text{Enc}(\rvx)=\dfrac{\text{NN}_{\boldsymbol{\theta}_\phi}(\rvx)}{\|\text{NN}_{\boldsymbol{\theta}_\phi}(\rvx)\|},\\
			\text{(Decoder)}~&\rvy\in \rmS^1\to \rvx\in\rmC: \text{Dec}(\rvy)=\text{NN}_{\boldsymbol{\theta}_{\boldsymbol{\chi}}}(\rvy),\\
			&\phi_{\boldsymbol{\theta}_\phi}(\rvx)= \rvh \left(\text{Enc}(\rvx)\right),\\
			&\boldsymbol{\chi}_{\boldsymbol{\theta}_{\boldsymbol{\chi}}}(\phi)= \text{Dec}\left(\rvh^{-1}(\phi)\right),
		\end{aligned}
	\end{equation}
	where Enc: $\rmC\to\rmS^1$ is the neural network (NN) parametrised by $\vtheta_\phi$, Dec: $\rmC\to\rmS^1$ is the neural network parametrised by $\vtheta_{\boldsymbol\chi}$, and $\rvh$ is the homeomorphism between $\rmS^1$ to $\mathbb{R}/_\sim$ as
	\begin{equation}
		\begin{aligned}
			\rvh(x,y)&=\left\{\begin{aligned}
				&\arcsin y,~x\ge0,~y\ge0\\
				&\arccos x,~x<0,~y\ge0,\\
				&\pi-\arcsin y,~x\le0,~y<0\\
				&2\pi-\arccos x,~x>0,y<0.\\
			\end{aligned}\right.   \\
			\rvh^{-1}(\phi)&=(\cos(\phi),\sin(\phi)).
		\end{aligned}
	\end{equation}
	We also set the natural frequency as a learnable quantity $w$.  Then, we train the parameters $\boldsymbol{\theta}_{\phi}$, $\boldsymbol{\theta}_{\boldsymbol{\chi}}$ and $w$ using the loss function $\mathcal{L}_{\text{phase}}$ as follows, 
	\begin{equation}
		\mathcal{L}_{\text{phase}}(\boldsymbol{\theta}_{\phi},\boldsymbol{\theta}_{\boldsymbol{\chi}},w)=\lambda_1\mathcal{L}_1(\boldsymbol{\theta}_{\phi},w)+\lambda_2\mathcal{L}_2(\boldsymbol{\theta}_{\boldsymbol{\chi}}),
	\end{equation}
	with
	\begin{equation}\label{eq loss1}
		\begin{aligned}
			\mathcal{L}_1(\boldsymbol{\theta}_{\phi},w)&=\mathbb{E}_{\rvx\sim\mathcal{U}(\rmC)}\bigg\Vert\dfrac{\partial\phi_{\boldsymbol{\theta}_{\phi}}}{\partial\rvx}\cdot \rmF(\rvx)-w\bigg\Vert_2^2,\\
			\mathcal{L}_2(\boldsymbol{\theta}_{\boldsymbol{\chi}})&=\mathbb{E}_{\rvx\sim\mathcal{U}(\rmC)}\big\|\rvx-\boldsymbol{\chi}_{\boldsymbol{\theta}_{\boldsymbol{\chi}}}\left(\phi_{\boldsymbol{\theta}_{\phi}}(\rvx)\right)\big\|_2^2.\\
		\end{aligned}
	\end{equation}
	Here $\mathcal{U}(\rmC)$ represents the uniform distribution over the limit cycle $\rmC$ and $\lambda_{1,2}$ are the predefined regularisation parameters. In practice, we collect the trajectory data $\{\rvx_i\}_{i=1}^K$ over the limit cycle $\rmC$ as training data and use the following empirical approximation of loss function in Eq.~\eqref{eq loss1} to train the parameters. 
	\begin{equation}\label{eq mc loss1}
		\begin{aligned}
			\mathcal{L}_1(\boldsymbol{\theta}_{\phi},w)&=\dfrac{1}{K}\sum_{i=1}^K\bigg\Vert\dfrac{\partial\phi_{\boldsymbol{\theta}_{\phi}}}{\partial\rvx}\cdot \rmF(\rvx_i)-w\bigg\Vert_2^2,\\
			\mathcal{L}_2(\boldsymbol{\theta}_{\boldsymbol{\chi}})&=\dfrac{1}{K}\sum_{i=1}^K\|\rvx_i-\boldsymbol{\chi}_{\boldsymbol{\theta}_{\boldsymbol{\chi}}}\left(\phi_{\boldsymbol{\theta}_{\phi}}(\rvx_i)\right)\|_2^2,\\
			\mathcal{L}_{\text{phase}}(\boldsymbol{\theta}_{\phi},\boldsymbol{\theta}_{\boldsymbol{\chi}},w)&=\lambda_1\mathcal{L}_1(\boldsymbol{\theta}_{\phi},w)+\lambda_2\mathcal{L}_2(\boldsymbol{\theta}_{\boldsymbol{\chi}}),
		\end{aligned}
	\end{equation}
	where the first term  $\mathcal{L}_1$ in the loss function measures whether the auto-encoder maps the original dynamics near the limit cycle to the universal dynamical clock, and the second term establishes the analytical homeomorphism between the limit cycle and $\rmS^1$.

	\subsection{Second step: Discovering the equant with the invertible neural network}
	Given the equant $\rvx^\ast$ and the limit cycle $\rmC$, the line segments $l_{\rvx^\ast\rvx(t)}$ between the equant and the oscillator state on the limit cycle form a manifold  $\mathbb{M}=\{\tilde{\rvx}:\tilde{\rvx}\in l_{\rvx^\ast\rvx},\rvx\in\rmC\}$. The viewing angle $\theta$ observed by the equant is not exactly the phase $\phi$ because the period of $\theta$ may not be $2\pi$. For universal dynamical clock on the unit circle, the centre of the unit disk $\mathbb{D}=\{(x,y):x^2+y^2\le1\}$ is the equant and naturally equates the phase as its viewing angle.  To capture the relation between the phase and the equant for general oscillation, we extend the homeomorphism between the limit cycle and $\rmS^1$ established by an auto-encoder to the invertible map between $\mathbb{M}$ and $\mathbb{D}$. Such an extension $\rvg$ is defined as 
 {\color{black}
	\begin{equation*}
		\begin{aligned}
			&\rvg\big|_{\rmC=\partial\mathbb{M}} = \text{Enc},~\rvg^{-1}\big|_{\rmS^1=\partial\mathbb{D}} = \text{Dec},\\
			&\rvg(\rvx^\ast)=\boldsymbol{0},~\rvx^\ast\in\arg\min K(\rvx),\\
			&\rvg\big|_{\mathbb{M}}: r\boldsymbol{\chi}(\phi)+(1-r)\rvx^\ast\mapsto r(\cos\phi,\sin\phi), \\
		\end{aligned}
	\end{equation*}
}
	
	We employ  a special class of generative model, the invertible NN~\cite{dinh2016density}, to parametrise the extension as $\rvg_{\vtheta_{\rvg}}$. The structure of our NN is comprised of the following three modules.
	\begin{itemize}
		\item [(i)]
		\textbf{Invertible mapping.}~We respectively split the input $\rvx=(x_1,...,x_n)^\top$ and the output $\rvy=(y_1,...,y_n)^\top=\rvf_I(\rvx)$ of the invertible non-linear mapping into two halves, respectively, as $\rvx=(\rvx_{1:d}^\top,\rvx_{d+1:n}^\top)^\top$ and $\rvy=(\rvy_{1:d}^\top,\rmY_{d+1:n}^\top)^\top$ , in which $d<n$ is a hyperparameter. The split variables obey the following equations,
		\begin{equation}
			\begin{aligned}
				\rvy_{1:d} &= \rvx_{1:d},\\
				y_{j} &= x_j\exp(s_j(\rvx_{1:d}))+t_j(\rvx_{1:d}),~d+1\le j\le n,
			\end{aligned}
		\end{equation}
		where $\rvs=(s_{d+1},...,s_n)$ and $\rvt=(t_{d+1},...,t_n)$ are  parametrised NN functions. The corresponding formulation of inverse mapping is trivially calculated as follows,
		\begin{equation*}
			\begin{aligned}
				\rvx_{1:d} &= \rvy_{1:d},\\
				x_{j} &= \exp\left(-s_j(\rvy_{1:d})\right)\left[y_j-t_j(\rvy_{1:d})\right],~d+1\le j\le n.
			\end{aligned}
		\end{equation*}
		\item [(ii)]
		\textbf{Permutation mapping.}~Since the invertible mapping keeps half of the state variables unchanged, we introduce the permutation mapping to transfer the unchanged variables to the next layer with a non-linear activation function. The permutation mapping and its inverse mapping are given by 
		\begin{equation}\label{eq permutation}
			\begin{aligned}
				&\rvf_P((x_1,\cdots,x_n)^\top)=(x_{\alpha(1)},\cdots,x_{\alpha(n)})^\top,~\alpha\in \rmS_n\\
				&\rvf_P^{-1}((x_1,\cdots,x_n)^\top)=(x_{\beta(1)},\cdots,x_{\beta(n)})^\top,~\beta=\alpha^{-1},
			\end{aligned}
		\end{equation}
		where $\rmS_n$ is the permutation group on $\{1,\cdots,n\}$. In implementation, we randomly sample the permutation element and obtain its inverse element via the efficient \texttt{randperm} and \texttt{argsort} methods in Pytorch package~\cite{paszke2019pytorch} as $	\alpha=\texttt{randperm}(n),~\beta=\texttt{argsort}(\alpha)$. 
		
		By alternatively combining the invertible mapping and the permutation mapping, we obtain the parametrised invertible function, 
		\begin{equation}\label{eq INN}
			\rvg_{\boldsymbol{\theta}_{\rvg}}=\rvf_I\circ\rvf_P\cdots\rvf_P\circ\rvf_I.
		\end{equation}
		\item [(iii)]
		\textbf{Dimension alignment.} Although both $\mathbb{M}\subset\mathbb{R}^n$ and $\mathbb{D}\subset\mathbb{R}^2$ are $2$-dimensional manifolds, the above parametrized function cannot be directly employed to find the target bijection since the input and the output of $\rvg_{\boldsymbol{\theta}_{\rvg}}$ belong to $\mathbb{R}^n$. To align the dimension between $\rvg_{\boldsymbol{\theta}_{\rvg}}$ and $\rvg$, we augment the unit disk $\mathbb{D}$ to an $n$-dimensional space as $\mathbb{D}_n=\text{Aug}(\mathbb{D},n-2)$, where $\text{Aug}(\cdot,\cdot)$ is an augmenting operation defined as
		\begin{equation}	
			\text{Aug}(\mathbb{D},k)=\mathbb{D}\times\{0\}^{k}.	
		\end{equation}  
		In this way, the parametrised invertible function relates to the target bijection as
		\begin{equation}
			\rvg_{\boldsymbol{\theta}_{\rvg}}(\rvx)=\text{Aug}(\rmG(\rvx),n-2)\in\mathbb{D}_n,~\text{for all}~\rvx\in\mathbb{M}.
		\end{equation} 
		
	\end{itemize}
	Finally, we train the NNs according to the following loss function $\mathcal{L}_{\text{equant}}$ as
	\begin{equation}
		\mathcal{L}_{\text{equant}}(\boldsymbol{\theta}_{\rmG})=\lambda_3\mathcal{L}_3(\vtheta_{\rvg})+\lambda_4\mathcal{L}_4(\vtheta_{\rvg})+\lambda_5\mathcal{L}_5(\vtheta_{\rvg}),
	\end{equation}
	
	with
	\begin{equation}\label{eq mc loss2}
		\begin{aligned}
			\mathcal{L}_{3}(\vtheta_{\rvg})&=\dfrac{1}{K}\sum_{i=1}^K\|\rvg_{\boldsymbol{\theta}_{\rvg}}(\rvx_i)-\text{Enc}(\rvx_i)\|_2^2,\\
		\mathcal{L}_4(\vtheta_{\rvg})&=	{\color{black}\frac{\sum_{i=1}^{K}\Delta_i\left(A_i-\sum_{k=1}^{K}\Delta_kA_k\right)^2}{\left(\sum_{i=1}^K\Delta_iA_i\right)^2},~A_i=\frac{1}{2}\|\left(\rvg_{\boldsymbol{\theta}_{\rvg}}^{-1}(\boldsymbol{0})-\rvx_i\right)\wedge\rmF(\rvx_i)\|,~\Delta_i=\frac{(t_{i+1}-t_{i})\omega}{2\pi}},\\
			\mathcal{L}_5(\vtheta_{\rvg})&=\dfrac{1}{K}\sum_{i=1}^K\sum_{j=1}^M\bigg\|\rvg_{\boldsymbol{\theta}_{\rvg}}\left(\dfrac{j}{M}\rvg_{\boldsymbol{\theta}_{\rvg}}^{-1}(\boldsymbol{0})+\dfrac{M-j}{M}\rvx_i\right)-\dfrac{M-j}{M}\rvg_{\vtheta_{\rvg}}(\rvx_i)\bigg\|_2^2.
		\end{aligned}
	\end{equation}
	{\color{black}Here we require the dataset $\{\rvx_i=\rvx(t_i)\}_{i=1}^K\subset\rmC$ contain temporally ordered limit-cycle samples with known sampling intervals, and $\lambda_{3,4,5}$ to be predefined parameters.} The loss term $\mathcal{L}_3$ guarantees that $\rvg_{\vtheta_{\rvg}}$ is an extension of the homeomorphism between the limit cycle and $\rmS^1$, {\color{black} $\mathcal{L}_4$ is used to find the  equant that has minimal $\hat{K}$ value}, and $\mathcal{L}_5$ measures whether the line of sight $l_{\rvg_{\boldsymbol{\theta}_{\rvg}}^{-1}(\boldsymbol{0})\rvx_i}$ in the original space corresponds to the line of sight $l_{\boldsymbol{0}\rvg_{\vtheta_{\rvg}}(\rvx_i)}$ of universal dynamical clock, i.e.,
	\begin{equation}
		\rvg_{\boldsymbol{\theta}_{\rvg}}\left(l_{\rvg_{\boldsymbol{\theta}_{\rmG}}^{-1}(\boldsymbol{0})\rvx_i}\right)=l_{\boldsymbol{0}\rvg_{\vtheta_{\rvg}}(\rvx_i)},~\rvg_{\boldsymbol{\theta}_{\rvg}}\left(l_{\rvg_{\boldsymbol{\theta}_{\rmG}}^{-1}(\boldsymbol{0})\rvx_i}\right)=\left\{\rvg_{\boldsymbol{\theta}_{\rvg}}(\rvx):\rvx\in l_{\rvg_{\boldsymbol{\theta}_{\rvg}}^{-1}(\boldsymbol{0})\rvx_i}\right\}.
	\end{equation}
    
\noindent\textbf{Definition of the interpretable phase}
	
	To comprehend the phase function from a physical viewpoint with the equant, we consider the following system for example,
	\begin{equation}\label{eq circle}
		\begin{aligned}
			&\dot{x}=y,~\dot{y}=-x,~\dot{z}=-z,\\
			&x(0)=\frac{\sqrt{2}}{2},~y(0)=\frac{\sqrt{2}}{2}~z(0)=z_0.\\
		\end{aligned}
	\end{equation}
	The limit cycle of this system is the unit circle $\rmS^1$ in the $(x,y)$ plane, and the oscillator rotates uniformly about $\rmS^1$. Therefore, any point on $z$-axis is an equant denoted by $\rvx^\ast(z)$, and the viewing angle observed from the origin  $\rvx^\ast(0)$ is the phase.
	Since there exists a linear projection from line $l_{\rvx^\ast(z)\rvx}$ to line $l_{\rvx^\ast(z)\rvx}$ for any $\rvx\in\rmS^1$, we can rescale the viewing angle of any equant to the phase by multiplying by a constant. However, such linear operation does not exist for a general oscillator with irregular limit cycle. In our framework, we address this issue by identifying an invertible map between the surface $\mathbb{M}$ and the domain $\mathbb{D}$. This map acts similarly to a projection, as demonstrated in the example. Finally we come up with the definition of interpretable phase with the equant in our framework as follows.
	
	\noindent\textbf{Definition 1.2}~\textit{
		For a phase function $\phi$, if there exists a bijection $\rvg$ between $\mathbb{M}$ and $\mathbb{D}$ such that $\mathrm{(i)}$ $\rvg$ maps the lines of sight emanating from the equant $\rvx^\ast$ to the lines of sight emanating from the centre of $\mathbb{D}$, i.e., $	\rvg(l_{\rvx^\ast\rvx})=l_{\boldsymbol{0}\rvg(\rvx)}$ for any $\rvx\in\mathbb{M}$, and $\mathrm{(ii)}$ for the temporal trajectory $\rvx(t)$ over the limit cycle, the intersection angle between the mapped lines $\rvg\left(l_{\rvx^\ast\rvx(t)}\right)$ exactly is equal to the phase, i.e., $\phi(\rvx(t))=\angle(\rvg(l_{\rvx^\ast\rvx(t)}),l_{\vzero\rve_1})$ with $\rve_1=(1,0)^\top$. Then we define $\phi$ as the interpretable phase.}

	\subsection{Analytical phase reduction against perturbations }\label{Methods phase reduction}
	
	After the training stage, we obtain  $\phi_{\vtheta_\phi}(\rvx)$ and $\boldsymbol{\chi}_{\vtheta_{\boldsymbol\chi}}(\phi)$. Hence, we obtain the phase response curve as 
	\begin{equation}
		\rmZ(\phi)=\dfrac{\partial\phi_{\boldsymbol{\theta}_{\phi}}}{\partial\rvx}\bigg|_{\rvx=\boldsymbol{\chi}_{\boldsymbol{\theta}_{\boldsymbol{\chi}}}(\phi)},
	\end{equation} 
	which can be efficiently  computed by \textsc{autograd} method in \textsc{PyTorch} \cite{paszke2017automatic}.
	We utilise these functions to deduce the analytical form of phase response function $\Gamma(\phi,t)$ under three different perturbations. For simplicity, we omit the subscripts $\vtheta_\phi$ and $\vtheta_{\boldsymbol{\chi}}$ in the following derivations.

	\begin{itemize}
		\item [(i)]\textbf{Stochastic force.} Firstly, we focus on the stochastic oscillation governed by a stochastic differential equation (SDE)
		\begin{equation}
			\begin{aligned}
				\mathrm{d}{\rvx}&=\rmF(\rvx)\mathrm{d}t+\rmG(\rvx)\mathrm{d}\rmW(t),\\
			\end{aligned}
		\end{equation}
		where $\rmG(\rvx)\in\mathbb{R}^{n\times m}$ is the diffusion term, and $\rmW(t)\in\mathbb{R}^m$ is the standard Wiener process. By applying I$\hat{\text{t}}$o's formula~\cite{oksendal2013stochastic} and the averaging method, we obtain the phase dynamics on the limit cycle as
		\begin{equation}
			\begin{aligned}
				\mathrm{d}{\phi}(\rvx) &=\rmZ(\phi)\cdot\mathrm{d}\rvx+\frac{1}{2}\mathrm{d}\rvx^\top\underbrace{\nabla^2\phi(\rvx)\big|_{\rvx=\boldsymbol{\chi}(\phi)}}_{\text{denoted by}~\rmY(\phi)}\mathrm{d}\rvx\\
				&=\rmZ(\phi)\cdot\left[\rmF(\boldsymbol{\chi}(\theta))\mathrm{d}t+\rmG(\boldsymbol{\chi}(\theta))\mathrm{d}\rmW(t)\right]+\frac{1}{2}\text{Tr}\left[\rmG(\boldsymbol{\chi}(\phi))^\top \rmY(\phi)\rmG(\boldsymbol{\chi}(\phi))\right]\mathrm{d}t\\
				&=\left\{w+\frac{1}{2}\text{Tr}\left[\rmG(\boldsymbol{\chi}(\phi))^\top \rmY(\phi)\rmG(\boldsymbol{\chi}(\phi))\right]\right\}\mathrm{d}t+\boldsymbol\Gamma(\phi)\cdot\mathrm{d}\rmW(t),\\
				\boldsymbol\Gamma(\phi)&=(\rmZ(\phi)\cdot\rmG_1(\boldsymbol{\chi}(\phi)),...,\rmZ(\phi)\cdot\rmG_m(\boldsymbol{\chi}(\phi))^\top,
			\end{aligned}
		\end{equation}
		here $\text{Tr}[\cdot]$ is the trace operator. 	The corresponding Fokker-Planck equation of the phase distribution $p(\phi,t)$ is,
		\begin{equation}\label{eq PR noise FK}
			\begin{aligned}
				\partial_t p(\phi,t)&=-\partial_\phi\left\{\left(w+\frac{1}{2}\text{Tr}\left[\rmG(\boldsymbol{\chi}(\phi))^\top \rmY(\phi)\rmG(\boldsymbol{\chi}(\phi))\right]p(\phi,t)\right)\right\} \\
				&+\frac{1}{2} \sum_{i=1}^m\partial_{\phi}^2\left[\rmZ(\phi)\cdot\rmG_i(\boldsymbol{\chi}(\phi))p(\phi,t)\right].
			\end{aligned}
		\end{equation}
		The influence of noise to oscillation is reflected in the stationary distribution of the above Fokker-Planck equation. 
		
		To further simplify the phase dynamics, we consider the phase decomposition $\phi=wt+\psi$, with $\psi$ representing the effect of noise to the phase. Then we have,
		\begin{equation}\label{eq PR noise slow}
			\begin{aligned}
				\mathrm{d}{\phi} &=w+\mathrm{d}{\psi},\\
				\mathrm{d}{\psi}&=\frac{1}{2}\text{Tr}\left[\rmG(\boldsymbol{\chi}(wt+\psi))^\top \rmY(wt+\psi)\rmG(\boldsymbol{\chi}(wt+\psi))\right]\mathrm{d}t+\boldsymbol\Gamma(wt+\psi)\cdot\mathrm{d}\rmW(t),\\
			\end{aligned}
		\end{equation}
		When the noise perturbation is weak, $\psi$ is a slow component and approximately stays
		constant over one cycle of oscillation with frequency $w$, the drift term on the right hand side of Eq.~\eqref{eq PR noise slow} changes the $\psi$ over one period $T=\frac{2\pi}{w}$ as
		\begin{equation}\label{eq PR noise delta psi}
			\delta\psi=\int_0^{T}\frac{1}{2}\text{Tr}\left[\rmG(\boldsymbol{\chi}(wt+\psi))^\top \rmY(wt+\psi)\rmG(\boldsymbol{\chi}(wt+\psi))\right]\mathrm{d}s.\\
		\end{equation}
		Then the dynamics of $\delta\psi$ during one period is equivalent to a constant,
		\begin{equation}\label{eq PR noise delta w}
			\begin{aligned}
				\delta w&=\frac{\delta\psi}{T}\\
				&=\frac{1}{T}\int_0^{T}\frac{1}{2}\text{Tr}\left[\rmG(\boldsymbol{\chi}(wt+\psi))^\top \rmY(wt+\psi)\rmG(\boldsymbol{\chi}(wt+\psi))\right]\mathrm{d}s\\
				&=\frac{1}{2\pi}\int_0^{2\pi}\frac{1}{2}\text{Tr}\left[\rmG(\boldsymbol{\chi}(\theta))^\top \rmY(\theta)\rmG(\boldsymbol{\chi}(\theta))\right]\mathrm{d}\theta.
			\end{aligned}
		\end{equation}
		Putting Eqs.~\eqref{eq PR noise delta psi},\eqref{eq PR noise delta w} into Eq.~\eqref{eq PR noise slow} leads to the following dynamics,
		\begin{equation}
			\begin{aligned}
				\mathrm{d}{\psi} &= \delta w\mathrm{d}t+\boldsymbol\Gamma(wt+\psi)\cdot\mathrm{d}\rmW(t).
			\end{aligned}
		\end{equation}
		Similarly, we perform the time-average operation to the diffusion term and obtain the simplified phase dynamics as,
		\begin{equation}
			\begin{aligned}
				\mathrm{d}{\psi} &= \delta w\mathrm{d}t+\sum_{i=1}^m\Gamma_i\mathrm{d}W_i(t),\\
				\Gamma_i&=\frac{1}{2\pi}\int_0^{2\pi}\rmZ(\theta)\cdot\rmG_i(\boldsymbol{\chi}(\theta))\mathrm{d}\theta.
			\end{aligned}
		\end{equation}
		Therefore, we obtain the simplified phase dynamics as,
		\begin{equation}\label{eq SDE phase}
			\begin{aligned}
				\mathrm{d}{\phi} = \left(w+\delta w\right)\mathrm{d}t+\sum_{i=1}^m\Gamma_i\mathrm{d}W_i(t).
			\end{aligned}
		\end{equation}
		%
		Although we use the standard Wiener process in our derivation, the results presented above can be easily extended to the more general case involving correlated Wiener processes.
		
		To proceed, we investigate the synchronisation behaviour of collective dynamics driven by scalar common noise and independent noise as,
		\begin{equation}\label{eq common+independ}
			\mathrm{d}{\rvx}_i=\rmF(\rvx)\mathrm{d}t+\rmG(\rvx)\mathrm{d}\xi(t)+\rmH(\rvx)\mathrm{d}\eta_i(t),~i=1,\cdots,N.
		\end{equation}
		We denote the corresponding phase dynamics as,
		\begin{equation}
			\begin{aligned}
				\mathrm{d}{\phi}_i &= \left(w+\delta w\right)\mathrm{d}t+ \Gamma_1(\phi_i)\mathrm{d}\xi+\Gamma_2(\phi_i)\mathrm{d}\eta_i,\\
				\delta w &= \dfrac{1}{2\pi}\int_0^{2\pi}\frac{1}{2}\text{Tr}\left[(\rmG+\rmH)^\top \rmY(\theta)(\rmG+\rmH)\right]\mathrm{d}\theta.\\
			\end{aligned}
		\end{equation}
		Here $\Gamma_1(\phi)=\rmZ(\phi)\cdot \rmG(\boldsymbol{\chi}(\phi))$, $\Gamma_2(\phi)=\rmZ(\phi)\cdot \rmH(\boldsymbol{\chi}(\phi))$, and the Wiener processes $\xi$ and $\{\eta_i\}_{i=1}^N$ are common noise and independent noise, respectively.
		The synchronisation behaviour of the noise driven dynamics is reflected in the phase difference  $\phi_i-\phi_j,~1\le i\neq j\le N$. We follow the approach in~\cite{nakao2007noise} to deduce the stationary distribution of the phase difference. We first consider the covariance matrix $\rmD=(D_{ij})\in\mathbb{R}^{N\times N}$ of $\{\mathrm{d}\phi_i\}_{i=1}^N$ as,
		\begin{equation}
			\begin{aligned}
				D_{ij}(\boldsymbol{\phi})&=\Gamma_1(\phi_i)\Gamma_1(\phi_j)+\delta_{ij}\Gamma_2(\phi_i)\Gamma_2(\phi_j),\\
				\boldsymbol{\phi}&=(\phi_1,...,\phi_N)^\top.
			\end{aligned}
		\end{equation}
		Next, we consider the Fokker-Planck equation of distribution $q$ of slow phase component $\psi_i=\phi_i-(w+\delta w)t$ as, 
		\begin{equation}
			\partial_{t}q(\boldsymbol{\psi},t)=\frac{1}{2}\sum_{i=1}^N\sum_{j=1}^N\partial_{\psi_i}\partial_{\psi_j}\left[D_{ij}(\boldsymbol{\psi}+(w+\delta w)t)q(\boldsymbol{\psi},t)\right].
		\end{equation}
		By time-averaging the covariance $D_{ij}$ over one period $T=\frac{2\pi}{w+\delta w}$, we obtain the simplified Fokker-Planck equation,
		\begin{equation}\label{eq PR noise diff FK}
			\begin{aligned}
				\partial_{t}q(\boldsymbol{\psi},t)&=\frac{1}{2}\sum_{i=1}^N\sum_{j=1}^N\partial_{\psi_i}\partial_{\psi_j}\left[\left(\gamma_1(\psi_i-\psi_j)+\delta_{ij}\gamma_2(\psi_i-\psi_j) \right)q(\boldsymbol{\psi},t)\right].\\
				\gamma_i(\psi)&=\frac{1}{2\pi}\int_0^{2\pi}\Gamma_i(\theta)\Gamma_i(\psi+\theta)\mathrm{d}\theta,~i=1,2.
			\end{aligned}
		\end{equation}
		
		Without loss of generality, we focus on the phase difference $\delta\phi=\phi_1-\phi_2=\psi_1-\psi_2$. The evolution equation for distribution of $\psi_1$, $\psi_2$ can be derived
		from $q(\boldsymbol{\psi})$ by integrating over all other phase variables as, 
		\begin{equation}
			R(\psi_1,\psi_2,t)=\int q(\boldsymbol{\psi})\mathrm{d}\psi_3\cdots\psi_N.
		\end{equation}
		According to Eq.~\eqref{eq PR noise diff FK}, the evolution equation of $R(\psi_1,\psi_2)$ is, 
		\begin{equation}\label{eq PR noise 1,2}
			\partial_tR(\psi_1,\psi_2,t)=\frac{\gamma_1(0)+\gamma_2(0)}{2}\partial_{\psi_1}^2R+\frac{\gamma_1(0)+\gamma_2(0)}{2}\partial_{\psi_2}^2R+\partial_{\psi_1}\partial_{\psi_2}\left[\gamma_1(\delta\phi)R\right].
		\end{equation}
		By transforming the two variables $\psi_1,~\psi_2$ to the mean $\psi=\frac{\psi_1+\psi_2}{2}$ and the difference $\delta\phi=\psi_1-\psi_2$, we have 
		\begin{equation}
			\begin{aligned}
				R(\psi_1=\psi+\frac{\delta\phi}{2},\psi_2=\psi-\frac{\delta\phi}{2},t)=Q(\psi,t)P(\delta\phi,t).
			\end{aligned}
		\end{equation}
		
		Then we have 
		\begin{equation}
			\begin{aligned}
				\partial_tR&=P\partial_tQ+Q\partial_tP\\
				\partial_{\psi_1}^2R&=\frac{1}{4}P\partial_\psi^2Q+(\partial_\psi Q)(\partial_{\delta\phi} P)+Q\partial_{\delta\phi}^2 P,\\
				\partial_{\psi_2}^2R&=\frac{1}{4}P\partial_\psi^2Q-(\partial_\psi Q)(\partial_{\delta\phi} P)+Q\partial_{\delta\phi}^2 P,\\
				\partial_{\psi_1}\partial_{\psi_2}\left[\gamma_1(\delta\phi)R\right]&=\frac{1}{4}\gamma_1(\delta\phi)P\partial_\psi^2Q-Q\partial_{\delta\phi}^2 \left[\gamma_1(\delta\phi)P\right].
			\end{aligned}
		\end{equation}
		
		Thus, the Eq.~\eqref{eq PR noise 1,2} is decoupled as,
		
		\begin{equation}
			\begin{aligned}
				\partial_t Q&=\frac{1}{4}\left[\gamma_1(0)+\gamma_1(\delta\phi)+\gamma_2(0)\right]\partial_{\psi}^2Q,\\
				\partial_t P&=\partial_{\delta\phi}^2\left\{\left[\gamma_1(0)-\gamma(\delta\phi)+\gamma_2(0)\right]P\right\}.
			\end{aligned}
		\end{equation}
		
		Then the stationary distribution of phase difference $P(\delta \phi)$ is determined by, 
		\begin{equation}\label{eq FK common+independ}
			P^\ast(\delta\phi)=\dfrac{p_0}{\gamma_1(0)-\gamma_1(\delta\phi)+\gamma_2(0)},
		\end{equation}
		where $P_0$ is the normalisation constant. Therefore, the long term collective behaviour of the noise perturbed oscillators can be readily determined by calculating $P^\ast$ using our framework. 
		

		\item [(ii)]\textbf{Periodic force.}~We are concerned whether an oscillator driven by the external periodic force
		$\rvf(t)$ with the frequency $\Omega$ shows phase lock phenomenon, i.e., if the frequency of the driven
		oscillator converges to $\Omega$. The perturbed dynamics is described as
		\begin{equation}
			\dot{\rvx}=\rmF(\rvx)+ \rvf(t),~\rvf(t)=\rvf(t+T).
		\end{equation}
		The corresponding phase dynamics is,
		\begin{equation}
			\begin{aligned}
				\dot{\phi}&=\rmZ(\phi)\cdot\left[\rmF(\rvx)+\rvf(t)\right]\\
				&=w+\rmZ(\phi)\cdot\rvf(t).
			\end{aligned}
		\end{equation}
		To investigate the influence of the periodic force to the phase frequency, we consider the phase difference $\delta \phi=\phi-\Omega t$ whose dynamics is governed by, 
		\begin{equation}\label{eq phase period}
			\begin{aligned}
				\dot{\delta \phi}&=\dot{\phi}-\Omega\\
				&=\Delta+\rmZ(\Omega t+\delta\phi)\cdot\rvf(t),\\
				\Delta&=w-\Omega.~
			\end{aligned}
		\end{equation}
		By time-averaging the weak term over one cycle of oscillation of frequency $\Omega$, we obtain the closed-form dynamics of phase difference as,
		\begin{equation}
			\begin{aligned}
				\dot{\delta \phi}&=\Delta+\Gamma(\delta \phi),\\
				\Gamma(\delta \phi) &= \dfrac{1}{2\pi}\int_0^{2\pi}\rmZ(\delta \phi+\theta)\rvf\left(\frac{\theta}{\Omega}\right)\mathrm{d}\theta.
			\end{aligned}
		\end{equation} 
		The phase lock occurs once the vector field $\Delta+\Gamma(\delta\phi)$ has a stable zero solution. Therefore, by viewing the intersection of constant $\Delta$ and function $\Gamma(\delta \phi)$ in the diagram, our framework efficiently discerns whether the periodic force leads to phase lock or divergence. 
		
		\item [(iii)] \textbf{Network interaction.} 
		We consider the collective dynamics of coupled oscillators, expressed in a general form as
		\begin{equation}
			\begin{aligned}	\dot{\rvx}_i&=\rmF(\rvx_i)+\sum_{j=1}^NA_{ij}\rmH(\rvx_i,\rvx_j),\\
			\end{aligned}
		\end{equation}
		where $\rmA = (A_{ij})\in\mathbb{R}^{N\times N}$ captures the interacting structure among the oscillators, and  $\rmH(\rvx_i,\rvx_j)$ represents the influence of $j$th oscillator to $i$th oscillator. Then we have the phase dynamics in the vicinity of $\rmC$ as
		\begin{equation}\label{eq PR network}
			\begin{aligned}
				\dot{\phi}_i &= \rmZ(\phi_i)\cdot\rmF(\rvx_i)+\sum_{j=1}^NA_{ij}\rmZ(\phi_i)\cdot\rmH(\rvx_i,\rvx_j)\\
				&=w+\sum_{j=1}^NA_{ij}\rmZ(\phi_i)\cdot\rmH(\boldsymbol{\chi}(\phi_i),\boldsymbol{\chi}(\phi_j))\\
			\end{aligned}
		\end{equation}
		from which the phase response function of network interaction is identified by a bivariate function $\tilde{\Gamma}(\phi_i,\phi_j)=\rmZ(\phi_i)\cdot\rmH(\boldsymbol{\chi}(\phi_i),\boldsymbol{\chi}(\phi_j))$. Even though the form $\tilde{\Gamma}(\phi_i,\phi_j)$ is clear and precise, a more predominant way is to simplify the phase response function as a univariate function $\Gamma(\phi_j-\phi_i)$ by averaging $\Gamma$ over one period~\cite{kuramoto1984chemical,nakao2016phase}. In this way, the phase response only  relies on the phase difference and the reduced phase dynamics is described as,
		\begin{equation}\label{eq PR network final}
			\begin{aligned}
				\dot{\phi}_i &=w+
				\sum_{j=1}^{N}A_{ij}\Gamma(\phi_j-\phi_i),\\
				\Gamma(\phi_j-\phi_i)&=\dfrac{1}{2\pi}\int_0^{2\pi}\rmZ(\theta+\phi_i)\rmH\left(\boldsymbol{\chi}(\theta+\phi_i),\boldsymbol{\chi}(\theta+\phi_j)\right)\mathrm{d}\theta.
			\end{aligned}
		\end{equation}
		The proof is mainly based on the method in~\cite{kuramoto2019concept}, we provide the detailed derivation here for completeness. 
		Specifically, for Eq.~\eqref{eq PR network} we split the phase variable as $\phi_i=wt+\psi_i$. Since we require the perturbation, i.e., the network interaction to be weak, the contribution of the perturbation to the phase frequency is weaker than natural frequency $w$. Thus, $wt$ is the fast component of the phase and $\psi_i$ represents the slow component. Then we have,
		\begin{equation}\label{eq PR network slow}
			\begin{aligned}
				\dot{\phi}_i&=w+\dot{\psi}_i,\\
				\dot{\psi}_i&=\sum_{j=1}^NA_{ij}\rmZ(wt+\psi_i)\cdot\rmH(\boldsymbol{\chi}(wt+\psi_i),\boldsymbol{\chi}(wt+\psi_j)).\\
			\end{aligned}
		\end{equation}
		Since the component $\psi_i$ varies slowly during one cycle of oscillation with frequency $w$, we treat the $\psi_i$, $\psi_j$ on the right hand side of Eq.~\eqref{eq PR network slow} as constant over one period of oscillation. Then the vector field of $\psi_i$ on the left hand side of Eq.~\eqref{eq PR network slow} solely varies along $wt$ with period $T=\frac{2\pi}{w}$, the dynamics of which is equal to the time-average over this period under fixed $\psi_i$, $\psi_j$,
		\begin{equation}\label{eq PR network average}
			\begin{aligned}
				\dot{\psi}_i&=\frac{1}{T}\int_0^T\sum_{j=1}^NA_{ij}\rmZ(ws+\psi_i)\cdot\rmH(\boldsymbol{\chi}(ws+\psi_i),\boldsymbol{\chi}(ws+\psi_j))\mathrm{d}s\\
				&=\frac{1}{2\pi}\int_0^{2\pi}\sum_{j=1}^NA_{ij}\rmZ(\theta+\psi_i)\cdot\rmH(\boldsymbol{\chi}(\theta+\psi_i),\boldsymbol{\chi}(\theta+\psi_j))\mathrm{d}\theta\\
				&=\sum_{j=1}^NA_{ij}\left[\frac{1}{2\pi}\int_0^{2\pi}\rmZ(\theta+\psi_i)\cdot\rmH(\boldsymbol{\chi}(\theta+\psi_i),\boldsymbol{\chi}(\theta+\psi_j))\mathrm{d}\theta\right]\\
				&\triangleq \sum_{j=1}^NA_{ij}\Gamma(\psi_j-\psi_i).
			\end{aligned}
		\end{equation}
		We note that $\phi_j-\phi_i=\psi_j-\psi_i$, and finally the Eq.~\eqref{eq PR network final} is obtained by putting Eq.~\eqref{eq PR network average} into Eq.~\eqref{eq PR network slow}.

		We consider the homogeneous oscillators in the above analysis, next we extend our framework to heterogeneous oscillators as follows,
		\begin{equation}\label{eq nonidentical}
			\begin{aligned}				 
				\dot{\rvx}_i&=\rmF_i(\rvx_i)+\sum_{j=1}^NA_{ij}\rmH(\rvx_i,\rvx_j)\\
				\rmF_i(\rvx_i)&=\rmF(\rvx_i)+\delta\rmF_i(\rvx_i),~\dfrac{\|\delta\rmF_i(\rvx_i)\|}{\|\rmF_i(\rvx_i)\|}\ll1.
			\end{aligned}
		\end{equation}
		The corresponding phase dynamics takes the form as,
		\begin{equation}
			\dot{\phi}_i=w+\rmZ(\phi_i)\cdot\delta \rmF_i(\boldsymbol{\chi}(\phi_i))+\sum_{j=1}^{N}A_{ij}\Gamma(\phi_j-\phi_i).
		\end{equation}
		Similarly we consider the dynamics of slow component $\psi_i=\phi_i-wt$ as,
		
		\begin{equation}\label{eq PR network heter phase}
			\dot{\psi}_i=\rmZ(wt+\psi_i)\cdot\delta \rmF_i(\boldsymbol{\chi}(wt+\psi_i))+\sum_{j=1}^{N}A_{ij}\Gamma(\psi_j-\psi_i).
		\end{equation}
		The first term on the right-hand side of Eq.~\eqref{eq PR network heter phase} is replaced by the time-average over one period of oscillation,
		\begin{equation}\label{eq PR network heter average}
			\begin{aligned}
				\dot{\psi}_i&=\frac{1}{T}\int_0^T\rmZ(ws+\psi_i)\cdot\delta \rmF_i(\boldsymbol{\chi}(ws+\psi_i))\mathrm{d}s+\sum_{j=1}^{N}A_{ij}\Gamma(\psi_j-\psi_i)\\
				&=\frac{1}{2\pi}\int_0^{2\pi}\rmZ(\theta+\psi_i)\cdot\delta \rmF_i(\boldsymbol{\chi}(\theta+\psi_i))\mathrm{d}s+\sum_{j=1}^{N}A_{ij}\Gamma(\psi_j-\psi_i)\\
				&=\underbrace{\frac{1}{2\pi}\int_0^{2\pi}\rmZ(\theta)\cdot\delta \rmF_i(\boldsymbol{\chi}(\theta))\mathrm{d}s}_{(\rmZ,\boldsymbol{\chi}~\text{is periodic by $2\pi$})}+\sum_{j=1}^{N}A_{ij}\Gamma(\psi_j-\psi_i).
			\end{aligned}
		\end{equation}
		Then we obtain the reduced phase dynamics of Eq.~\eqref{eq nonidentical} as follows,
		\begin{equation}
			\begin{aligned}
				\dot{\phi}_i&=(w+\delta w_i)+\sum_{i=1}^{N}A_{ij}\Gamma(\phi_j-\phi_i),\\
				\delta w_i&=\rmZ(\phi_i)\cdot\delta\rmF(\rvx_i)=\dfrac{1}{2\pi}\int_0^{2\pi}\rmZ(\theta)\cdot\delta\rmF_i(\boldsymbol{\chi}(\theta))\mathrm{d}\theta.
			\end{aligned}
		\end{equation}
		
	\end{itemize}

	Although we derive the phase reduction for each of the three perturbations separately, the phase dynamics resulting from the combination of these perturbations can  easily be obtained by summing the corresponding phase response functions.
	
	For the analytical expression of the phase dynamics, we note that the NN functions are analytical when the parameters are fixed after training. To  enhance the interpretability of the established phase response functions, we approximate the NN-based function by a linear combination of known basis functions as
	\begin{equation}
		\Gamma(\phi)=\sum_{i=1}^{M}a_if_i(\phi),
	\end{equation}
	where $\{f_i:\mathbb{R}\to\mathbb{R}\}_{i=1}^{M}$ is the predefined dictionary of basis functions. Since the phase is a periodic function, the basis function can be further required as $f_i:\mathbb{R}/_\sim\to\mathbb{R}$, e.g., the Trigonometric basis.\\

\section{Experimental configuration of the machine learning framework}
	
	In this Section, we provide a detailed description for the experimental configurations of the numerical results.  The computing device that we use for calculating our examples includes a single i7-10870 CPU with 16GB memory, and we train all the parameters with Adam optimiser.   
    Our code will be available once the paper is published. 
    We summarise the experimental configuration in our training process as follows.
	
	\begin{itemize}
		\item \textbf{Auto-encoder.} For constructing the auto-encoder, we employ the forward NN as:
		\begin{equation*}
			\begin{aligned}
				\rvz_1 &= \mathcal{F}(\rmW_0\rvx+\rmB_0),\\
				\rvz_{i+1} &= \mathcal{F}(\rmW_i\rvz_{i}+\rmB_i),\ i=1,\cdots,k-1,\\
				\rvz_{k+1} &= \rmW_k\rvz_{k},\\
			\end{aligned}
		\end{equation*}
		where $\mathcal{F}(\cdot)$ is the activation function, and $\rmW_i\in\mathbb{R}^{h_{i+1}\times h_i}$ and $\rmB_i\in\mathbb{R}^{h_{i+1}}$ represent the weight and bias of each layer, respectively. The default activation function is $\mathcal{F}(\cdot)=\tanh$ unless otherwise specified. We denote the auto-encoder parametrized by this function as Enc$(h_0,h_1,\cdots,h_{k+1})$ and Dec$(h_0,h_1,\cdots,h_{k+1})$.

		\item \textbf{Invertible Neural Network.} We employ the standard invertible NN `RealNVP' in \textsc{Python} package \texttt{INNLab} (https://github.com/ELIFE-ASU/INNLab) to parametrise $\rvf_{I}$ in Eq.~\eqref{eq INN}, and we parametrize $\rvf_{P}$ as we described in Eq.~\eqref{eq permutation}. We denote our invertible NN-based function as Surface$(a,b)$, where $a$ refers to the number of the invertible functions $\rvf_{I}$ in Eq.~\eqref{eq INN}, and $b$ refers to the dimension of the invertible function. We use ``Surface'' to refer to that the function maps the surface $\mathbb{M}$ to $\mathbb{D}$ in Definition~1.2.
		
		\item \textbf{Dataset.} To collect the training data, we first simulate the original dynamics Eq.~\eqref{eq1} for sufficient long time $[0,T_{\text{end}}]$ such that the trajectory converges to the attractive limit cycle and rotates along the limit cycle over multiple times. Then we intercept the trajectory as $[T_{\text{intercept}},T_{\text{end}}]$ and uniformly sample the data on it. We set the time interval of the trajectory as the constant $\delta t$ to collect regular time series data $\{\rvx_i=\rvx(T_{\text{intercept}}+i\delta t)\}_{i=1}^K$. We choose small $\delta t$ for {\color{black} guaranteeing the dense sampling over the limit cycle. }
		
		\item \textbf{Initialization.}  We use the Normal initialisation for our NNs. We initialise the learnable phase frequency $w$ as $2\pi/\tilde{T}$, where $\tilde{T}$ is the rough approximation of the period $T$ calculated from the time series data.
		
		\item \textbf{Preprocessing.} We aim at finding the homeomorphism between the irregular limit cycle and the unit circle using the auto-encoder. To improve the training efficiency, we select the affine transformation $\rvy=\rmA\rvx+\rvb$ to adjust the position and the shape of the limit cycle, such that its projection to each plane is a closed curve having similarities with the unit circle, having low eccentricity and covering the origin. The corresponding dynamics of $\rvy$ is $\dot{\rvy}=\rmA\rmF(\rmA^{-1}(\rvy-\rvb))$. And after training, we use the inverse affine transformation $\rvx=\rmA^{-1}(\rvy-\rvb)$ to pull the learned equant and the functions back to the original space. 
		
		\item \textbf{Regularisation.} To improve the performance of the trained models, we introduce some auxiliary regularisations into the loss function. 
		\begin{itemize}
			\item To ensure the auto-encoder based phase function to satisfy the condition that the phase varies $2\pi$ when the training data go all the way over the limit cycle, we introduce the regularisation as $\mathcal{L}_{\text{ptp}}=|2\pi-\max_{\rvx_i}\phi_{\vtheta_\phi}(\rvx_i)-\min_{\rvx_i}\phi_{\vtheta_\phi}(\rvx_i)|^2$. 
			\item One trivial minimiser to the loss function $\mathcal{L}_{4}(\vtheta_\rvg)$ is the equant at infinity. To avoid such a case, we introduce the $L_2$ regularisation of the equant as $\mathcal{L}_{\text{norm}}=\|\rvg_{\vtheta_{\rvg}}^{-1}(\vzero)\|^2$. 
			
		\end{itemize}
		The above regularisation terms are not necessary for all the systems, we find that they improve the training performance for certain systems, such as the Mitotic dynamics. 
		\item \textbf{Training hyperparameters.} We train the auto-encoder for $10000$ iterations with learning rate $0.005$, and train the invertible neural network for $10000$ iterations with learning rate $0.001$. The learning rate is determined by grid search method.
	\end{itemize}

    \section{List of limit cycle oscillators}\label{sec details}
	Most of the oscillatory dynamics investigated in our numerical studies have been used previously to describe the rhythmic phenomena in nature. In this Section, we list the dynamical equations, the preprocessing matrices $(\rmA,\rvb)$, and the experimental configurations. {\color{black}We stress that the preprocessing matrices $(\rmA,\rvb)$ are used only for numerical conditioning and do not impose an equant-like geometry. Specifically, the affine transformation is chosen so that the origin lies inside the limit cycle in the relevant two-dimensional projections. The role of $\rmA$ is to rescale the trajectory, especially for systems with small-amplitude limit cycles, so that the machine-learning model can better resolve the oscillatory structure rather than treating it as a small perturbation around a point.
		Importantly, this preprocessing does not bias the intrinsic geometry of the dynamics, because we transform not only the limit-cycle trajectory but also the vector field consistently under the same affine coordinate change. That is, for $\rvy=\rmA\rvx+\rvb$, the transformed dynamics are given by
		\[
		\dot{\rvy}=\rmA\rmF(\rmA^{-1}(\rvy-\rvb)).
		\]
		Thus, the preprocessing is simply a linear change of coordinates rather than an artificial deformation of the dynamics. It preserves the underlying dynamical structure while improving numerical conditioning. This preprocessing is not essential for the framework, but empirically facilitates convergence.
	
 For the criterion of using the additional regularisation terms $L_{\mathrm{ptp}}$ and $L_{\mathrm{norm}}$. We clarify that these terms are not part of the default training objective and are only introduced when the basic loss converges to a trivial or degenerate solution, such as $w=0$ or a collapsed latent representation. In such cases, the optimisation landscape makes the basic framework difficult to train reliably. The role of $L_{\mathrm{ptp}}$ and $L_{\mathrm{norm}}$ is therefore to prevent these degenerate solutions and improve convergence. When the basic loss converges successfully without these terms, we do not use them.}
	
	\noindent \textbf{Morris-lecar model}~\cite{morris}:
	\begin{equation}\label{eq morris}
		\begin{aligned}
			C_M\dot{v}&=I_b-g_L(v-E_L)-g_Kn(v-E_K)-g_{Ca}m_{\infty}(v)(v-E_{Ca}),\\
			\dot{n}&=\phi\frac{n_{\infty}(v)-n}{\tau_n(v)},\\
			m_\infty(v)&=0.5\left(1+\tanh\left(\frac{v-v_1}{v_2}\right)\right),\\
			\tau_n(v)&=\frac{1}{\cosh(\frac{v-v_3}{2v_4})},\\
			n_\infty(v)&=0.5\left(1+\tanh\left(\frac{v-v_3}{v_4}\right)\right),\\
			\phi&=0.067,~g_{Ca}=4,~g_K=8,~g_L=2,~E_{Ca}=120,~E_K=-84,\\
			E_L&=-60,~v_1=-1.2,~v_2=18,~v_3=12,~v_4=17.4,~C_M=20,\\
			\rmA&=\mathrm{diag}\left(\frac{1}{10},10\right),~\rvb=(0,-2)^\top.
		\end{aligned}
	\end{equation}
	We set the network structure as Enc$(2,32,32,2)$, Dec$(2,32,32,2)$, and Surface$(7,3)$. For the dataset, we set $T_{\text{end}}=400$, $\delta t=0.04$, and $T_{\text{intercept}}=379.04$. We employ the auxiliary loss $\mathcal{L}_{\text{ptp}}$ and $\mathcal{L}_{\text{norm}}\times1$e-$4$ in the training process.\\
	
	\noindent \textbf{FitzHugh-Nagumo model}~\cite{fitzhugh1969mathematical}:
	\begin{equation}\label{eq FithHugh-Nagumo}
		\begin{aligned}
			\dot{x}&=x-x^3/3-y+I,\\
			\dot{y}&=\varepsilon(x+a-by),\\
			I&=0.0,~a=0.7,~b=0.2,~\varepsilon=0.05,\\
			\rmA&=\rmI,~\rvb=\vzero.
		\end{aligned}
	\end{equation}
	We set the network structure as Enc$(2,40,40,2)$, Dec$(2,40,40,2)$, and Surface$(7,3)$. For the dataset, we set $T_{\text{end}}=20$, $\delta t=0.002$ and $T_{\text{intercept}}=16.94$. We employ the auxiliary loss $\mathcal{L}_{\text{norm}}\times1$e-$4$ in the training process. For training we employ the learning rate $0.01$.\\
	
	\noindent \textbf{Wilson Cowan model}~\cite{wilson1972excitatory}:
	\begin{equation}\label{eq Wilson Cowan}
		\begin{aligned}
			\tau\dot{E}& =-E + (S_{E_{\max}} - E)S_e(c_1  E - c_2  I + P),\\
			\tau\dot{I}& = -I + (S_{I_{\max}} - I)S_i(c_3 E - c_4  I),\\
			S_e(x)&= \frac{1}{ (1 + \exp(-a_e (x - \theta_e)))} - \frac{1}{ (1 + \exp(a_e  \theta_e))},\\
			S_i(x)&= \frac{1}{ (1 + \exp(-a_i(x - \theta_i)))} - \frac{1}{ (1 + \exp(a_i  \theta_i))},\\
			\tau& = 4,~t_d = 0,~c_1 = 16,~c_2 = 12,~c_3 = 15,\\
			c_4 &= 3,~a_e = 1,~a_i = 2,~\theta_e = 4,~\theta_i = 3,~P = 1,\\
			\rmA&=\mathrm{diag}(10,10),~\rvb=(-2,-1)^\top.
		\end{aligned}
	\end{equation}
	We set the network structure as Enc$(2,32,32,2)$, Dec$(2,32,32,2)$, and Surface$(7,2)$. For the dataset, we set $T_{\text{end}}=200$, $\delta t=0.02$ and $T_{\text{intercept}}=186.72$. We employ the auxiliary loss $\mathcal{L}_{\text{ptp}}$ in the training process.\\
	
	\noindent  \textbf{Selkov glycolysis}~\cite{sel1968self}:
	\begin{equation}\label{eq Selkov}
		\begin{aligned}
			\dot{x}& = 1 - xy^\gamma,\\
			\dot{y}& =\alpha y(xy^{\gamma-1}-1),\\
			\alpha&=1.1,~\gamma=2,\\
			\rmA&=\mathrm{diag}(2,2),~\rvb=(-2.5,-2.5)^\top.
		\end{aligned}
	\end{equation}
	We set the network structure as Enc$(2,32,32,2)$, Dec$(2,32,32,2)$, and Surface$(7,2)$. For the dataset, we set $T_{\text{end}}=100$, $\delta t=0.02$ and $T_{\text{intercept}}=93.36$. We employ the auxiliary loss $\mathcal{L}_{\text{ptp}}$ in the training process.\\

    \noindent  \textbf{Lotka-Volterra}~\cite{bomze1983lotka}:
	\begin{equation}\label{eq Lotka-Volterra}
		\begin{aligned}
			\dot{N}& = N(1-\frac{N}{\kappa}-\frac{\rho P}{N+k_N}),\\
			\dot{P}& =\varepsilon P(\frac{\rho N}{N+k_N}-d_P),\\
			       \kappa &= 1.0,~\varepsilon = 0.5,~k_N = 0.5,~d_P = 0.1,~\rho = 0.5 \\
			\rmA&=\mathrm{diag}(10,10),~\rvb=(-2.5,-10)^\top.
		\end{aligned}
	\end{equation}
	We set the network structure as Enc$(2,32,32,2)$, Dec$(2,32,32,2)$, and Surface$(7,2)$. For the dataset, we set $T_{\text{end}}=400$, $\delta t=0.08$ and $T_{\text{intercept}}=358.32$. \\

    \noindent  \textbf{Laser}~\cite{dolcemascolo2020effective}:
	\begin{equation}\label{eq laser}
		\begin{aligned}
			\dot{x}& = x(y-1),\\
			\dot{y}& = \gamma(\delta-y+k(w+A\log(1+\alpha x))-xy),\\
            \dot{w}& =-\varepsilon(w+A\log(1+\alpha x)),\\
			    \gamma &= 0.04,~
        \varepsilon = 0.1,~
        \delta = 1.2,~
        \alpha = 2.0,~
        k = 0.7,
        A = 1/0.7, \\
			\rmA&=\rmI,~\rvb=(-0.6,-0.95,0.8)^\top.
		\end{aligned}
	\end{equation}
	We set the network structure as Enc$(3,32,32,2)$, Dec$(2,32,32,3)$, and Surface$(7,3)$. For the dataset, we set $T_{\text{end}}=1000$, $\delta t=0.2$ and $T_{\text{intercept}}=881.6$.  We employ the auxiliary loss $\mathcal{L}_{\text{ptp}}$ in the training process.\\
    
	\noindent \textbf{Mitotic oscillator}~\cite{mitotic}:
	\begin{equation}\label{eq Mitotic}
		\begin{aligned}
			\dot{C}&=v_i-v_dX\frac{C}{K_d+C}-k_dC,\\
			\dot{M}&=V_1\frac{1-M}{K_1+1-M}-V_2\frac{M}{K_2+M},\\
			\dot{X}&=V_3\frac{1-X}{K_3+1-X}-V_4\frac{X}{K_4+X},\\
			V_1&=\frac{V_{M1}C}{K_c+C},\\
			V_3&=MV_{M3},\\
			v_i &= 0.023,~v_d = 0.1,~K_d = 0.02,~k_d = 0.0033,~V_{M1} = 0.5,~V_2 = 0.167,\\
			V_{M3} &= 0.2,~V_4 = 0.1,~K_C = 0.3,~K_1 = 0.1,~K_2 = 0.1,~K_3 = 0.1,~K_4 = 0.1,\\
			\rmA&=\mathrm{diag}(5,5,5),~\rvb=(-0.675,-2.0,-1.5)^\top.
		\end{aligned}
	\end{equation}
	We set the network structure as Enc$(3,32,32,2)$, Dec$(2,32,32,3)$, and Surface$(7,3)$. For the dataset, we set $T_{\text{end}}=300$, $\delta t=0.03$ and $T_{\text{intercept}}=264.3$. We employ the auxiliary loss $\mathcal{L}_{\text{ptp}}$ and $\mathcal{L}_{\text{norm}}\times1$e-$5$ in the training process.\\

	\noindent \textbf{Thalamic neuron model}~\cite{rubin2004high}
	\begin{equation}\label{eq Thalamic}
		\begin{aligned}
			\dot{v}&=\frac{-I_L-I_{Na}-I_K-I_T+I_b}{C_m},\\
			\dot{h}&=\frac{h_\infty-h}{\tau_h},\\
			\dot{r}&=\frac{r_\infty-r}{\tau_r},\\
			h_\infty&=\frac{1}{1+\exp(\frac{v+41}{4})},\\
			r_\infty&=\frac{1}{1+\exp(\frac{v+84}{4})},\\
			\alpha_h&=0.128\exp\left(-\frac{v+46}{18}\right),\\
			\beta_h&=\frac{4}{1+\exp(-\frac{v+23}{5})},\\
			\tau_h&=\frac{1}{\alpha_h+\beta_h},\\
			\tau_r&=28+\exp\left(-\frac{v+25}{10.5}\right),\\
			m_\infty&=\frac{1}{1+\exp(-\frac{v+37}{7})},\\
			p_\infty&=\frac{1}{1+\exp(-\frac{v+60}{6.2})},\\
			I_L&=g_L(v-e_L),\\
			I_{Na}&=g_{Na}m_\infty^3h(v-e{Na}),\\
			I_K&=g_K\left(0.75(1-h)\right)^4(v-e_K),\\
			I_T&=g_Tp_\infty^2r(v-e_T),\\
			C_m&=1,~g_L=0.05,~e_L=-70,~g_{Na}=3,~e_{Na}=50,~\\
			g_K&=5,~e_K=-90,~g_T=5,~e_T=0,~I_b=5,\\
			\rmA&=\mathrm{diag}\left(\frac{1}{10},10,1000\right),~\rvb=(3,-4,-5)^\top.
		\end{aligned}
	\end{equation}
	We set the network structure as Enc$(3,32,32,2)$, Dec$(2,32,32,3)$, and Surface$(7,3)$. For the dataset, we set $T_{\text{end}}=50$, $\delta t=0.02$ and $T_{\text{intercept}}=41.24$. We employ the auxiliary loss $\mathcal{L}_{\text{ptp}}$ and $\mathcal{L}_{\text{norm}}\times1$e-$4$ in the training process.\\

	\noindent \textbf{Escherichia coli cell}~\cite{garcia2004modeling}:
	\begin{equation}\label{eq Escherichia coli}
		\begin{aligned}
			\dot{a}&=-a+\dfrac{\alpha}{1+C^n},\\
			\dot{b}&=-b+\dfrac{\alpha}{1+A^n},\\
			\dot{c}&=-c+\dfrac{\alpha}{1+B^n}+\dfrac{\kappa S}{1+S},\\
			\dot{A}&=\beta(a-A),\\
			\dot{B}&=\beta(b-B),\\
			\dot{C}&=\beta(c-C),\\
			\dot{S}&=-k_{s0}S+k_{s1}A-\eta S,\\
			\alpha&=216, n=2, \kappa=20, \beta=1, k_{s0}=1, k_{s1}=0.01, \eta=2,\\
			\rmA&=\rmI,~\rvb=\vzero.
		\end{aligned}
	\end{equation}
	We set the network structure as Enc$(7,140,140,2)$, Dec$(2,140,140,7)$, and Surface$(7,7)$. For the dataset, we set $T_{\text{end}}=200$, $\delta t=0.02$ and $T_{\text{intercept}}=185.6$. We employ the auxiliary loss $\mathcal{L}_{\text{ptp}}$ and $\mathcal{L}_{\text{norm}}\times1$e-$4$ in the training process.\\
	
	\noindent  \textbf{Skeleton model for the cyclin-dependent
		kinases network}~\cite{gerard2011skeleton}:
	\begin{equation}\label{eq Skeleton}
		\begin{aligned}
			\dot{E2F}&=V_{1e2f}\frac{E2F_{tot}-E2F}{K_{1e2f}+E2F_{tot}-E2F}(Md+Me)-V_{2e2f}\frac{E2F}{K_{2e2f+E2F}}Ma,\\
			\dot{Me}&=v_{se}E2F-V_{de}Ma\frac{Me}{K_{de}+Me},\\
			\dot{Ma}&=v_{sa}E2F-V_{da}Cdc20\frac{Ma}{K_{da}+Ma},\\
			\dot{Mb}&=v_{sb}Ma-V_{db}Cdc20\frac{Mb}{K_{db}+Mb},\\
			\dot{Cdc20}&=V_{1cdc20}Mb\frac{Cdc20_{tot}-Cdc20}{K_{1cdc20}+Cdc20_{tot}-Cdc20}-V_{2cdc20}\frac{Cdc20}{K_{2cdc20}+Cdc20},\\
			Md&=\frac{K_{dd}v_{sd}\frac{GF}{K_{gf}+GF}}{V_{dd}-\frac{v_{sd}GF}{K_{gf}+GF}},\\
			Cdc20_{tot}&=5,~E2F_{tot}=3,~GF=1,~K_{da}=0.1,~K_{db}=0.005,~K_{dd}=0.1~K_{de}=0.1,\\
			K_{gf}&=0.1,~K_{1cdc20}=1,~K_{2cdc20}=1,~K_{1e2f}=0.01,~K_{2e2f}=0.01,~V_{da}=0.245,\\
			V_{db}&=0.28,~V_{dd}=0.245,~V_{de}=0.35,~v_{sa}=0.175,~v_{sb}=0.21,~v_{sd}=0.175,\\
			v_{se}&=0.21,~V_{1cdc20}=0.21,~V_{2cdc20}=0.35,~V_{1e2f}=0.805,~V_{2e2f}=0.7,\\
			\rmA&=\rmI,~\rvb=\vzero.
		\end{aligned}
	\end{equation}
	
	We set the network structure as Enc$(5,50,50,2)$, Dec$(2,50,50,5)$, and Surface$(7,5)$. For the dataset, we set $T_{\text{end}}=200$, $\delta t=0.04$ and $T_{\text{intercept}}=179.44$. We employ the auxiliary loss $\mathcal{L}_{\text{ptp}}$ and $\mathcal{L}_{\text{norm}}\times1$e-$5$ in the training process.\\
	
	\noindent \textbf{Mammalian Circadian Clock}~\cite{leloup2003toward}: 
	
	The governing equation of the mammalian circadian clock has $19$ variables and $53$ meaningful parameters, we omit the equation and parameters here for the sake of simplicity. The code for simulating this 16-dimensional ODE will be available at our code repository \textsc{torchequant} once the paper is published.
    To preprocess the data, we set $\rvy=\rmA\rvx+\rvb$ with $\rmA=\rmI$, 
	\begin{equation}
		\begin{aligned}
			\rvb&=(-2,-1.5,-8.5,0.6,4,0.09,0.6,2.5,2,0,0,0,0,0,0,0)^\top.
		\end{aligned}
	\end{equation}
	We set the network structure as Enc$(16,32,32,2)$, Dec$(2,32,32,16)$, and Surface$(7,16)$. For the dataset, we set $T_{\text{end}}=300$, $\delta t=0.03$ and $T_{\text{intercept}}=276.12$. We employ the auxiliary loss $\mathcal{L}_{\text{ptp}}$ and $\mathcal{L}_{\text{norm}}\times1$e-$5$ in the training process.\\
	
	\noindent \textbf{Cyclin-dependent
		kinases (Cdks) network for the cell cycle}~\cite{gerard2009temporal}: 
	
	Similarly, we omit this 44-dimensional ODE and provide the code for reproducing the dynamical trajectory at \textsc{torchequant}, which will be available after the paper is accepted.
    To preprocess the data, we set $\rvy=\rmA\rvx+\rvb$ with $\rmA=\rmI$, 
	\begin{equation}
		\begin{aligned}
			\rvb&=(-5.5,-2,-1.2,-14,-4,-0.012,-6,-0.3,-0.08,-0.015,-0.6,-0.6,\\
			&-0.01,-0.02,-6,-0.06,-0.12,-1.5,-0.1,-0.2,-0.2,-0.03,-4,-0.1,\\
			&-0.4,-0.5,-0.4,-0.7,-2,-0.2,-0.2,-0.2,-0.5,-1,-0.3,-0.6,-0.2,\\
			&-0.2,-0.2,-0.1,-0.05,-0.004,-0.04)^\top.
		\end{aligned}
	\end{equation}
	We set the network structure as Enc$(44,88,88,2)$, Dec$(2,88,88,44)$, and Surface$(7,44)$. For the dataset, we set $T_{\text{end}}=200$, $\delta t=0.038$ and $T_{\text{intercept}}=171.0936$. We employ the auxiliary loss $\mathcal{L}_{\text{ptp}}$ and $\mathcal{L}_{\text{norm}}\times1$e-$5$ in the training process.
	

\section{Extension of the framework to data-driven scenario}\label{sed NODE}
	{\color{black}In the established two-step approach, knowledge of the data $\{\rvx_i,\rmF(\rvx_i)\}$ is required.} However, we emphasize that our primary concern is obtaining an interpretable phase based on the equant. Therefore, prior knowledge of the {\color{black}governing equation} of dynamics is not necessary. In this Section, we aim at enhancing the generalisation capability of our framework by extending it to {\color{black}governing-equation-free} scenarios, where only the data of the limit cycle is available.
	
\begin{figure}[htp]
	\centering
	\includegraphics[width=13.69cm]{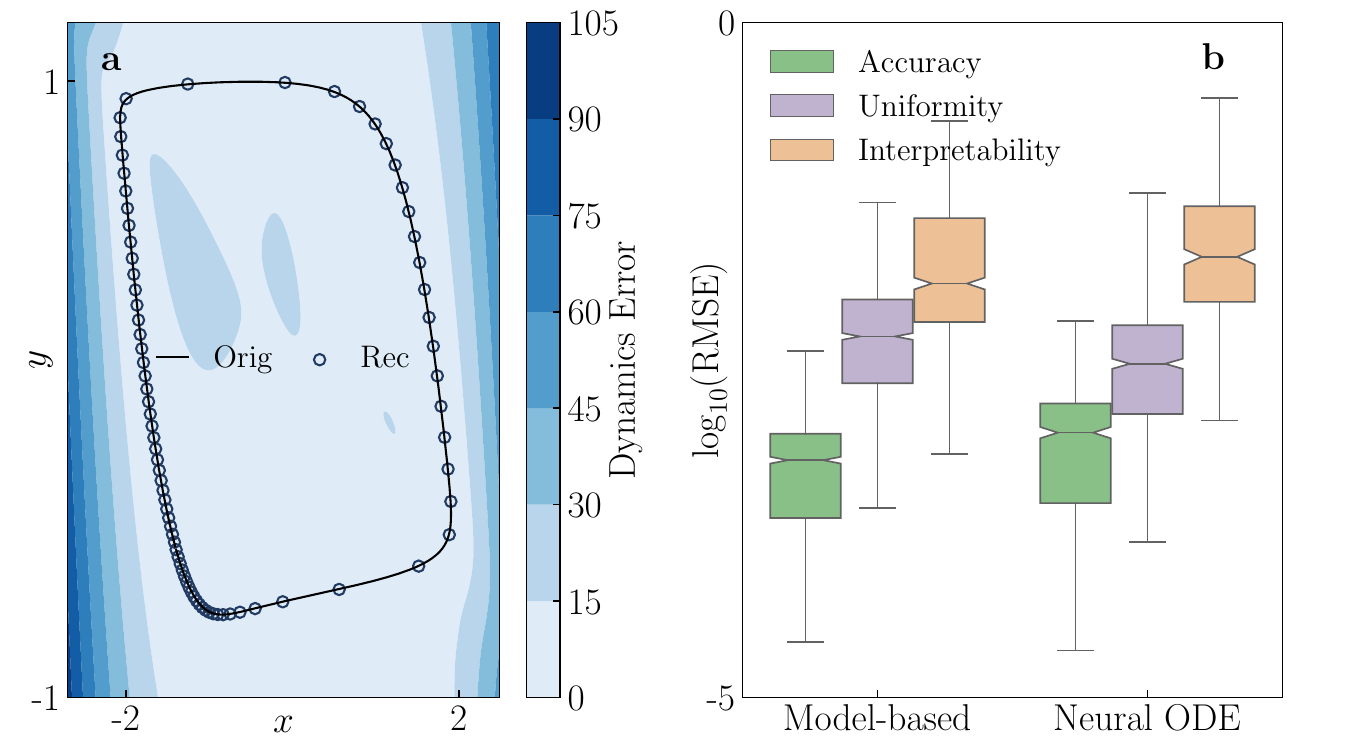}
	\caption{\textbf{Test Neural ODE in the FitzHugh-Nagumo dynamics}. a, The comparison between the real limit cycle (solid line) and the predicted limit cycle (hollow circle). The shadow region represents the error between the real dynamics and the learned dynamics. b, The comparison of the performance results of the real dynamics and the learned dynamics in our framework.}
	\label{figS4}
\end{figure}

	In the model reconstruction community, several benchmark methods are used to identify the governing equations behind time series data. Those methods include both parametrised approaches, such as SINDy~\cite{brunton2016discovering}, and non-parametrised methods like Neural ODE~\cite{chen2018neural} and reservoir computing~\cite{pathak2018model}.
	As an example, we apply our framework to a data-driven scenario using Neural ODE. We begin by considering the time series data $\{\rvx(t_i),~t_i\}_{i=0}^K$ collected from the underlying dynamics, $\dot{\rvx}=\rmF(\rvx)$ (the time interval between adjacent points is not necessarily equidistant). The Neural ODE parametrises the potential vector field as $\rmF_{\boldsymbol{\theta}}(\rvx)$, where $\boldsymbol{\theta}$ represents the learnable parameters of the neural network. By setting each data $\rvx(t_i)$ in the dataset as the initial position, the candidate dynamics predict the data at the next time step $\tilde{\rvx}(t_{i+1})$ as,
	\begin{equation}
		\tilde{\rvx}(t_{i+1})=\rvx(t_i)+\int_{t_i}^{t_{i+1}}\rmF_{\boldsymbol{\theta}}(\rvx(s))\mathrm{d}s.
	\end{equation}
	Then the parameters are updated according to the following loss function, such that the predicted data is consistent with the real data,
	\begin{equation}
		\mathcal{L}_{\text{Neural ODE}} = \frac{1}{K}\sum_{i=1}^K\|\rvx(t_i)-\tilde{\rvx}(t_{i+1})\|^2.
	\end{equation}
{\color{black} We notice that here we do not require the dataset is collected equidistantly over time, but can be irregular time series data with known time intervals $\{t_{i+1}-t_i\}$.}
	Based on the uniqueness of the solution, the candidate dynamics $\rmF_{\boldsymbol{\theta}}$ is expected to converge to the underlying dynamics $\rmF$ once the training loss converges to zero. We then apply our framework to the learned dynamics to obtain the equant and the interpretable phase. To demonstrate the effectiveness and performance of the extended framework, we apply it to the FitzHugh-Nagumo model, using only the data from the limit cycle to learn the dynamics.
	Once the dynamics are learned, we input them into our framework to derive the universal dynamical clock and the equant. As shown in Fig.~\ref{figS4}, the Neural ODE accurately learns the dynamics near the limit cycle, and the performance results from our framework align closely with those obtained from the actual dynamics. Therefore, we conclude that our framework naturally integrates with data-driven scenarios by combining with existing dynamics reconstruction methods.
	
\begin{figure}[htp]
	\centering
	\includegraphics[width=13.69cm]{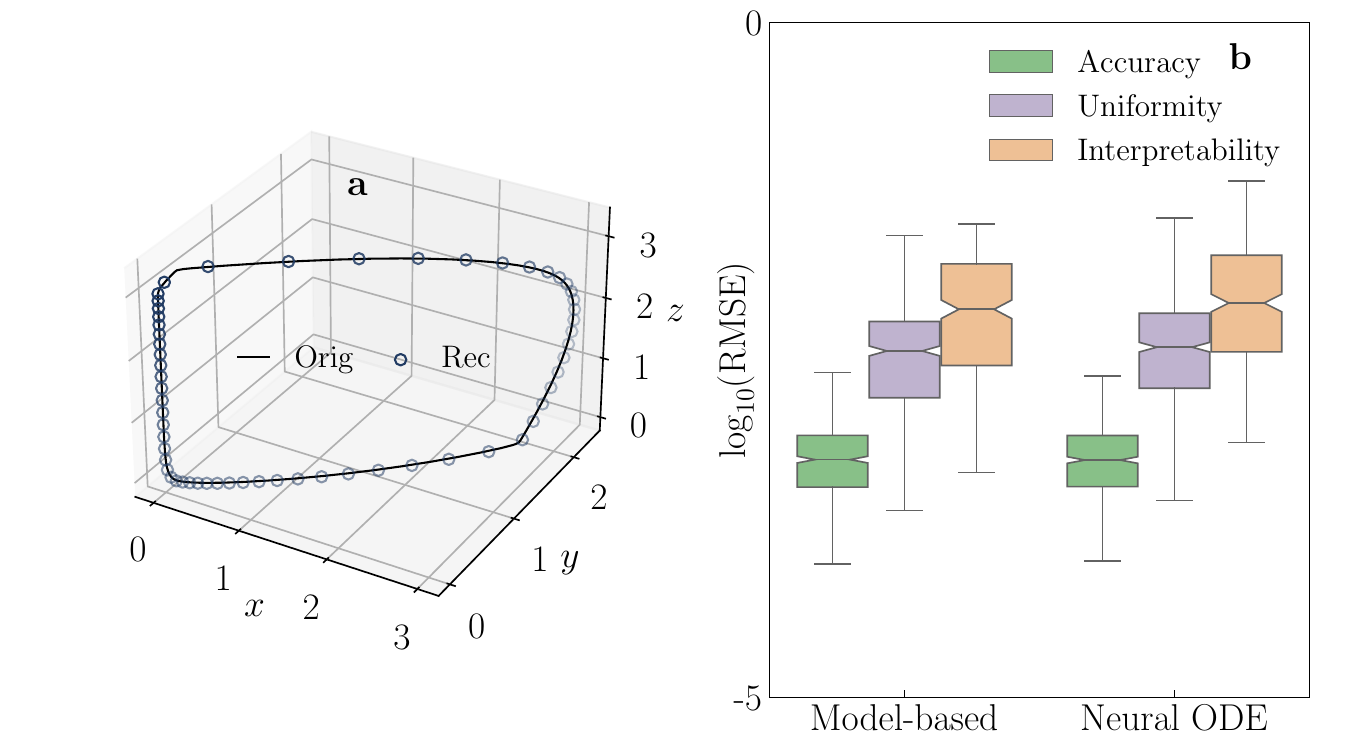}
	\caption{\textbf{Test Neural ODE in the Skeleton dynamics}. a, The comparison between the real limit cycle (solid line) and the predicted limit cycle (hollow circle). b, The comparison of the performance results of the real dynamics and the learned dynamics in our framework.}
	
	\label{fig_skeleton}
\end{figure}

\begin{figure}[htp]
	\centering
	\includegraphics[width=13.69cm]{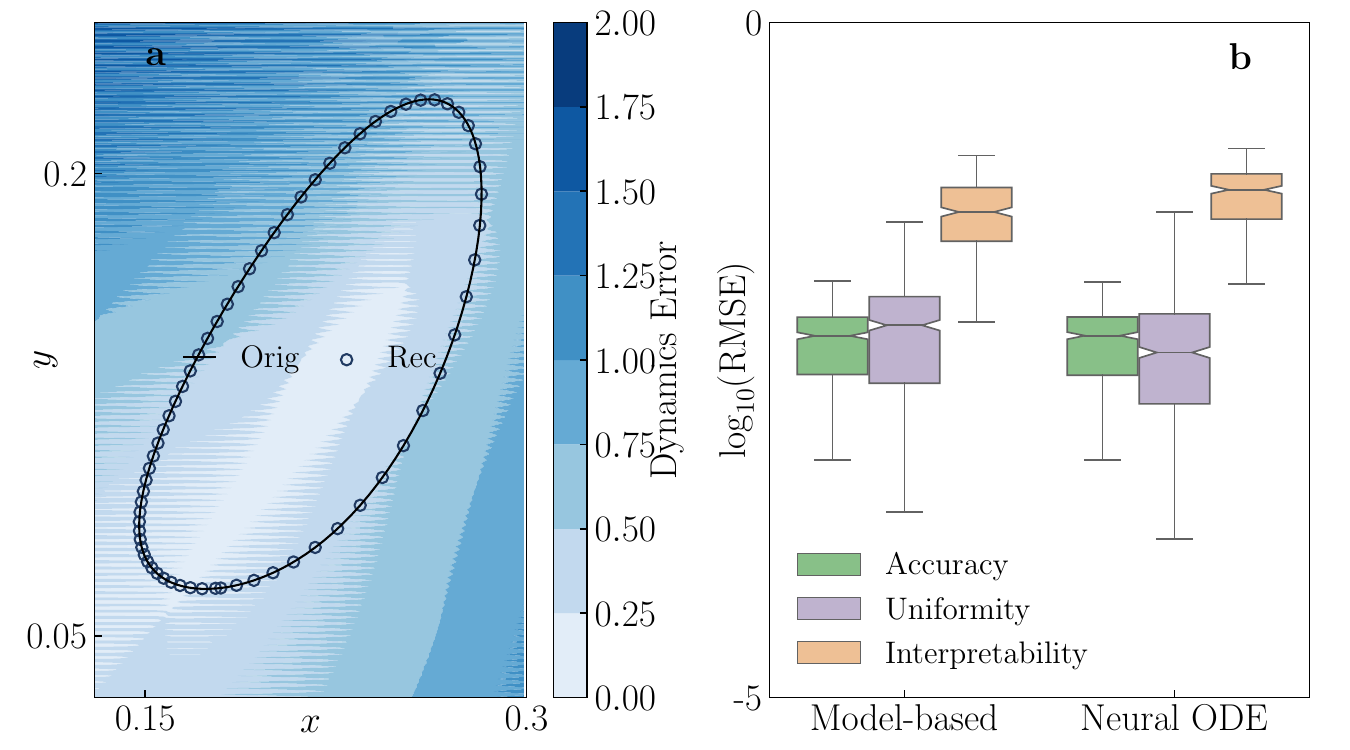}
	\caption{\textbf{Test Neural ODE in the Wilson-Cowan dynamics}. a, The comparison between the real limit cycle (solid line) and the predicted limit cycle (hollow circle). The shadow region represents the error between the real dynamics and the learned dynamics. b, The comparison of the performance results of the real dynamics and the learned dynamics in our framework.}
	
	\label{fig_wc}
\end{figure}

	{\color{black}
	Furthermore, to demonstrate the generality of the proposed Neural ODE based extension, we have supplemented the neural ODE extension with two additional numerical examples beyond FitzHugh--Nagumo: a 5-dimensional Skeleton system and a 2-dimensional Wilson--Cowan system. The corresponding results, shown in Figs.~\ref{fig_skeleton} and~\ref{fig_wc}, demonstrate that the neural-ODE-based extension can recover effective velocity fields and construct accurate phase dynamics across systems of different dimensions and dynamical structures. These additions make clear that the data-driven capability of our framework is not restricted to the FitzHugh--Nagumo model, but extends to both experimental circuit data and additional unknown-dynamics benchmarks.
}


{\color{black}
	\section{Effect of dimension augmentation on the optimal equant non-uniformity}
	
	We further discuss the effect of dimension augmentation on the optimal equant non-uniformity $\hat K^\ast$. Given a limit-cycle trajectory $\rvx(t)\in\mathbb R^n$ with vector field $\rmF(\rvx(t))$, we define its augmented representation in $\mathbb R^{n+m}$ as
	\[
	\tilde{\rvx}(t)=\left(\rvx(t),\mathbf 0_m\right),
	~
	\tilde{\rmF}(\tilde{\rvx}(t))
	=
	\left(\rmF(\rvx(t)),\mathbf 0_m\right).
	\]
	This augmentation does not change the underlying limit-cycle dynamics. It only embeds the same trajectory into a higher-dimensional ambient space. Consequently, the admissible set of observer points, and hence of possible equants, becomes larger:
	\[
	\tilde{\rvx}^{\ast}=(\rva,\rvb)\in\mathbb R^{n+m},
	~
	\rva\in\mathbb R^n,~ \rvb\in\mathbb R^m .
	\]
	The swept areal velocity observed from $\tilde{\rvx}^{\ast}$ is
	\[
	\tilde A_{\tilde{\rvx}^{\ast}}(t)
=	
	\frac{1}{2}
	\|
	\left(\tilde{\rvx}^{\ast}-\tilde{\rvx}(t)\right)
	\wedge
	\tilde{\rmF}(\tilde{\rvx}(t))
	\|.
	\]
	Using the identity $|\rvu\wedge\rvv|^2=|\rvu|^2|\rvv|^2-(\rvu\cdot\rvv)^2$, we obtain
	\begin{equation*}
			\begin{aligned}
			4\tilde A_{\tilde{\rvx}^{\ast}}^2(t)
			&=
			\left(|\rva-\rvx(t)|^2+|\rvb|^2\right)
			|\rmF(\rvx(t))|^2
			-	
			\left[(\rva-\rvx(t))\cdot \rmF(\rvx(t))\right]^2 \\
			&=
			4A_{\rva}^2(t)+|\rvb|^2|\rmF(\rvx(t))|^2,
		\end{aligned}
	\end{equation*}
	where
	$
	A_{\rva}(t)
=
	\frac{1}{2}
	\|
	\left(\rva-\rvx(t)\right)\wedge \rmF(\rvx(t))
	\|
	$
	is the swept areal velocity observed from $\rva$ in the original space.
	
	Therefore, zero-padding enlarges the family of admissible equants by allowing observer points with non-zero components $\rvb$ in the added dimensions. Optimising over this larger space may reduce the normalised areal non-uniformity:
	\[
	\hat K_{n+m}^{\ast}
=
	\min_{\tilde{\rvx}^{\ast}\in\mathbb R^{n+m}}
	\hat K_{n+m}(\tilde{\rvx}^{\ast})
	\leq
	\min_{\rvx^{\ast}\in\mathbb R^n}
	\hat K_{n+m}\left((\rvx^\ast,\mathbf 0_m)\right).
	\]
	Thus, the augmented representation can make the choice of equant more flexible and may yield a smaller value of $\hat K^\ast$.
	
	However, this effect should not be interpreted as changing the underlying oscillatory dynamics. The zero-padded trajectory and vector field remain confined to the original $n$-dimensional subspace. Hence, if the original system does not admit an observer point with exactly uniform swept areal velocity, zero-padding alone does not remove this geometric obstruction. Exact areal uniformity in the augmented space would require the existence of $\tilde{\rvx}^{\ast}=(\rva,\rvb)$ such that
	\[
	4A_{\rva}^2(t)+|\rvb|^2|\rmF(\rvx(t))|^2
	\equiv \mathrm{constant},
	~ \text{for all } t\in[0,T],
	\]
	which imposes a special cancellation relation between the original areal velocity $A_{\rva}(t)$ and the speed $|\rmF(\rvx(t))|$. Such a relation is not generically created by adding zero coordinates. Consequently, zero-padding can lower the optimal equant non-uniformity by expanding the admissible observer space, but it does not generically turn a system with $\hat K_n^\ast>0$ into one with $\hat K_{n+m}^\ast=0$.

\section{Comparison with Adjoint-based method}

We compare our original framework and the extension proposed in this section with adjoint-computed PRCs using the FitzHugh--Nagumo model. The adjoint method requires knowledge of both the vector field $\mathbf{F}(\mathbf{x})$ and its Jacobian $\nabla \mathbf{F}(\mathbf{x})$ along the limit cycle. It also relies on regularly sampled time-series data, since the phase is parametrized by time along one period of the orbit. For comparison, we consider three adjoint-based settings:
\begin{itemize}
	\item[(1)] Reference: The dataset consists of regularly sampled time-series data $D_1=\{\mathbf{x}_i,\nabla\mathbf{F}(\mathbf{x}_i)\}$ on the limit cycle, where $\nabla\mathbf{F}(\mathbf{x}_i)$ is computed analytically from the known governing equations.
	\item[(2)] Limit-cycle only: The dataset consists only of regularly sampled limit-cycle states $D_2$. In this setting, $\nabla\mathbf{F}(\mathbf{x}_i)$ is numerically approximated from limit-cycle samples. We consider two sampling densities, denoted by $D_{2,\mathrm{low}}$ and $D_{2,\mathrm{high}}$, with $\#D_{2,\mathrm{low}}=102$ and $\#D_{2,\mathrm{high}}=153$, respectively.
	
	\item[(3)] Off-cycle: The dataset consists of regularly sampled states $D_3$ in a neighbourhood of the limit cycle. The Jacobian $\nabla\mathbf{F}(\mathbf{x}_i)$ is then numerically approximated from these off-cycle samples. We again consider two sampling densities, corresponding to low- and high-density sampling around the limit cycle.
\end{itemize}
Similarly, for our ML-based method with a neural ODE model in Section~\ref{sed NODE}, we test the performance using $D_{2,\mathrm{low}}$ and $D_{2,\mathrm{high}}$ separately.  As shown in Fig.~\ref{fig_adjoint}, our ML-based method achieves performance comparable to the off-cycle adjoint method while only using the  limit-cycle-only data. In contrast, the limit-cycle-only adjoint method exhibits substantially larger errors. These results demonstrate that the ML-based method is more reliable than the data-driven adjoint method when only limit-cycle data are available.

The reason is that estimating $\nabla \mathbf{F}(\mathbf{x})$ numerically from limit-cycle samples is ill-conditioned, because such data only cover a one-dimensional manifold and provide insufficient information about transverse variations of the vector field. The resulting Jacobian errors are further amplified when solving the boundary-value problem in the adjoint method. By contrast, our ML framework fits the trajectory data through the solution flow of a neural ODE, thereby avoiding direct numerical estimation of $\nabla \mathbf{F}$ from sparse limit-cycle samples.

\begin{figure}[htp]
	\centering
	\includegraphics[width=\textwidth]{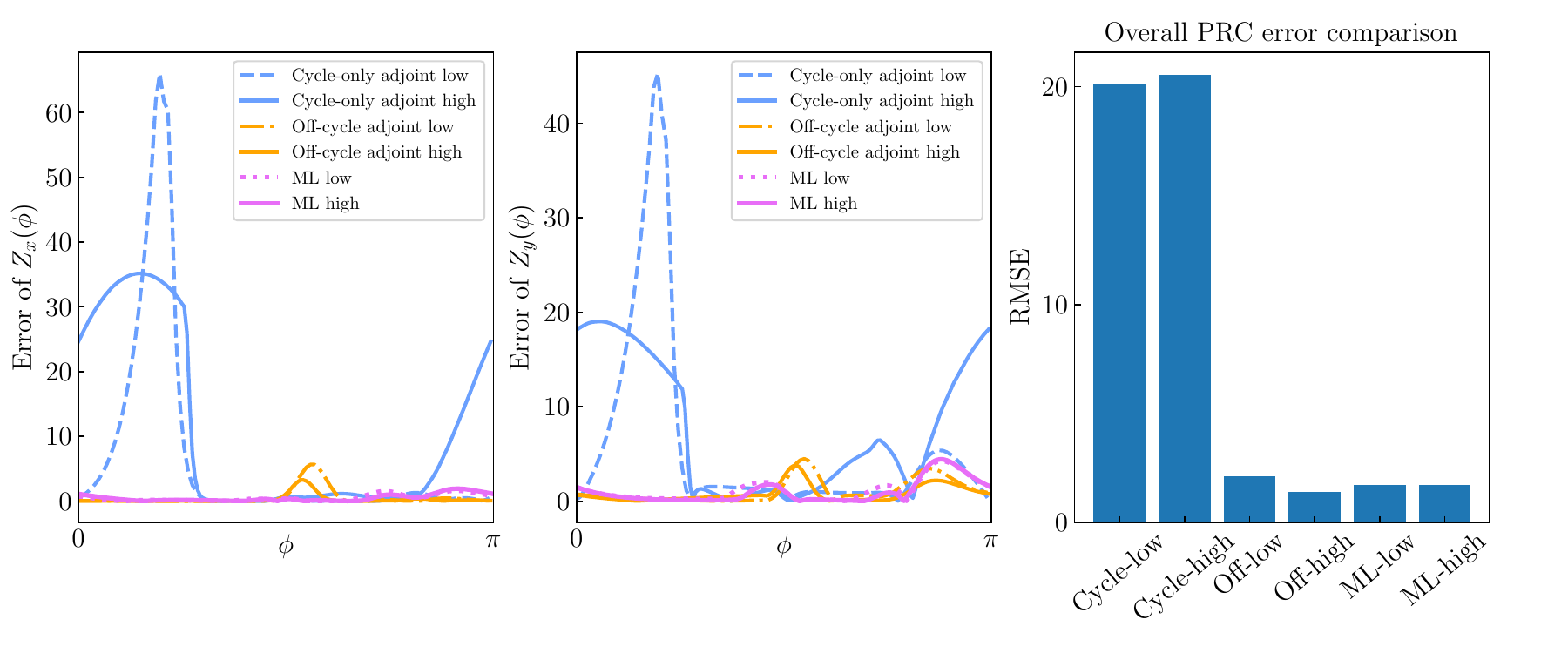}
	\caption{Comparison between ML approach and adjoint methods on FitzHugh-Nagumo model.}
	
	\label{fig_adjoint}
\end{figure}

\section{Extension of the framework with adjoint regularisation}

In the main framework, the phase map $\phi(\rvx)$, its inverse map $\chi(\phi)$, and the phase response curve $
	\rmZ(\theta)
	=
	\frac{\partial \phi}{\partial \rvx}
	\bigg|_{\rvx=\chi(\theta)}$	
are obtained from the auto-encoder trained on the limit-cycle data. This construction is sufficient for deriving the reduced phase dynamics considered in this work. Nevertheless, for completeness, we further introduce an adjoint-regularised extension of our framework, in which the transversal gradient information encoded by the Jacobian of the vector field is explicitly incorporated into the training process.

For a smooth limit-cycle oscillator
\begin{equation}
	\dot{\rvx} = \rmF(\rvx),
\end{equation}
the classical adjoint equation for the PRC along the limit cycle can be written as~\cite{monga2019phase}
\begin{equation}
	\frac{\mathrm{d} \rmZ(\theta)}{\mathrm{d} \theta}
	=
	-\nabla \rmF(\chi(\theta)) \rmZ(\theta),
	\label{eq:adjoint_equation}
\end{equation}
together with the normalisation condition
\begin{equation}
	\rmZ(\theta) \cdot \rmF(\chi(\theta)) = w.
	\label{eq:adjoint_normalisation}
\end{equation}
Here $\nabla \rmF(\chi(\theta))$ denotes the Jacobian of the vector field evaluated along the limit cycle. In our framework, $\rmZ(\theta)$ is directly computed from the learned phase map $\phi_{\vtheta_\phi}$ by automatic differentiation, and thus Eq.~\eqref{eq:adjoint_equation} can be imposed as an additional regularisation term.

Specifically, given phase samples $\{\theta_i\}_{i=1}^{K}$ on the unit circle and the corresponding reconstructed limit-cycle states $\chi_{\theta_\chi}(\theta_i)$, we define the adjoint residual as
\begin{equation}
	\mathcal{R}_{\rm adj}(\theta_i)
	=
	\frac{\partial \rmZ(\theta)}{\partial \theta}
	\bigg|_{\theta=\theta_i}
	+
	\nabla \rmF(\chi_{\theta_\chi}(\theta_i))
	\rmZ(\theta_i).
	\label{eq:adjoint_residual}
\end{equation}

The corresponding adjoint loss is defined as the empirical $L^2$ norm of the residual:
\begin{equation}
	\mathcal{L}_{\rm adj}
	=
	\frac{1}{K}
	\sum_{i=1}^{K}
	\left\|
	\mathcal{R}_{\rm adj}(\theta_i)
	\right\|_2^2.
	\label{eq:adjoint_loss}
\end{equation}
Then, the auto-encoder training objective in the first step can be extended as
\begin{equation}
	\mathcal{L}_{\text{phase}}(\boldsymbol{\theta}_{\phi},\boldsymbol{\theta}_{\boldsymbol{\chi}},w)=\lambda_1\mathcal{L}_1(\boldsymbol{\theta}_{\phi},w)+\lambda_2\mathcal{L}_2(\boldsymbol{\theta}_{\boldsymbol{\chi}})+\lambda_{\text{adj}}\mathcal{L}_{\text{adj}},
	\label{eq:adjoint_regularised_loss}
\end{equation}
where $\lambda_{\rm adj}$ is a predefined regularisation parameter. The first two terms preserve the original role of our framework: $\mathcal{L}_1$ learns a phase map with constant velocity along the limit cycle, and $\mathcal{L}_2$ establishes the homeomorphism between the limit cycle and the unit circle. The additional term $\mathcal{L}_{\rm adj}$ further constrains the learned PRC to satisfy the adjoint equation along the limit cycle. The loss function can be employed in different scenarios as follows.

\textbf{Model-driven scenario.}
When the governing equation $\dot{\rvx}=\rmF(\rvx)$ is known, the Jacobian $\nabla \rmF(\rvx)$ in Eq.~\eqref{eq:adjoint_residual} can be computed either analytically or through automatic differentiation. This setting corresponds to the classical adjoint formulation, but the PRC is not solved from a boundary-value problem. Instead, it is induced by the learned phase map and regularised to satisfy the adjoint equation.

\textbf{Data-driven scenario.}
In the governing-equation-free case, the analytical form of $\rmF$ is unavailable and only the trajectory data $\{\rvx(t_i),t_i\}_{i=0}^{K}$ are given. Following the Neural ODE extension introduced above, we first learn a surrogate vector field $\rmF_{\vtheta}(\rvx)$ from data. After the Neural ODE has been trained, we replace $\rmF$ in the phase-learning loss with $\rmF_{\vtheta}$. The Jacobian required by the adjoint regularisation is then computed as
\begin{equation*}
	\nabla \rmF_{\vtheta}(\rvx)
	=
	\frac{\partial \rmF_{\vtheta}(\rvx)}{\partial \rvx},
\end{equation*}
which can be efficiently obtained through automatic differentiation. Therefore, even when the governing equation is unknown, the transversal gradient information can still be incorporated into our framework through the Jacobian of the Neural ODE-learned vector field. This provides a data-driven counterpart of the classical adjoint method while preserving the original structure of our approach, in which the PRC is induced by the learned phase map rather than solved as a separate adjoint boundary-value problem.

We further numerically investigate whether the proposed extension improves the prediction of the PRC $Z(\phi)$ compared with the original method, under the same experimental setting described in Section~\ref{sec details}. For the extended framework, we additionally provide the Jacobian data $\nabla \rmF(\rvx_i)$ along the limit cycle. We evaluate the adjoint residual error of the learned $Z(\phi)$ in Eq.~\eqref{eq:adjoint_loss} on the test phase grid $\{\phi_i=2\pi i/1000\}_{i=1}^{1000}$. The results show that the proposed extension reduces the normalised adjoint residual error to $0.12$ of that obtained by the original framework.

\section{Ablation study for training robustness}

In this section, we provide additional ablation studies to examine the robustness of the proposed machine-learning framework with respect to random seeds, network architectures, learning rates, training iterations, and preprocessing choices. 

\subsection{Ablation study for Step 1: learning the dynamical clock}

We first investigate the robustness of Step 1, where the auto-encoder learns the phase map $\phi$, the inverse map $\chi$, and the natural frequency $w$ by minimising the phase-learning loss. We use the FitzHugh--Nagumo system as the benchmark model and vary the main training configurations, including random seed from $\{1,2,3,4,5\}$, network width from $\{40,60,80\}$, learning rate (lr) from $\{0.01,0.001,0.0001\}$, and number of training iterations until $10000$. We report the reconstruction loss $\mathcal{L}_1$ and phase-dynamics loss $\mathcal{L}_2$ during training.

As shown in Fig.~\ref{fig_ablation_step1}, the training curves converge across different random seeds, network widths, learning rates, and training iterations. In all tested settings where the learning rate is in a suitable range, both the reconstruction loss and the phase-dynamics loss decrease consistently and converge to small values. Among them we also found the optimal learning rate can be set as $\mathrm{lr}=0.001$. These results indicate that the phase map and inverse map obtained in Step 1 are not sensitive to a specific random initialization or a particular network architecture.

\begin{figure}[htp]
	\centering
	\includegraphics[width=0.8\textwidth]{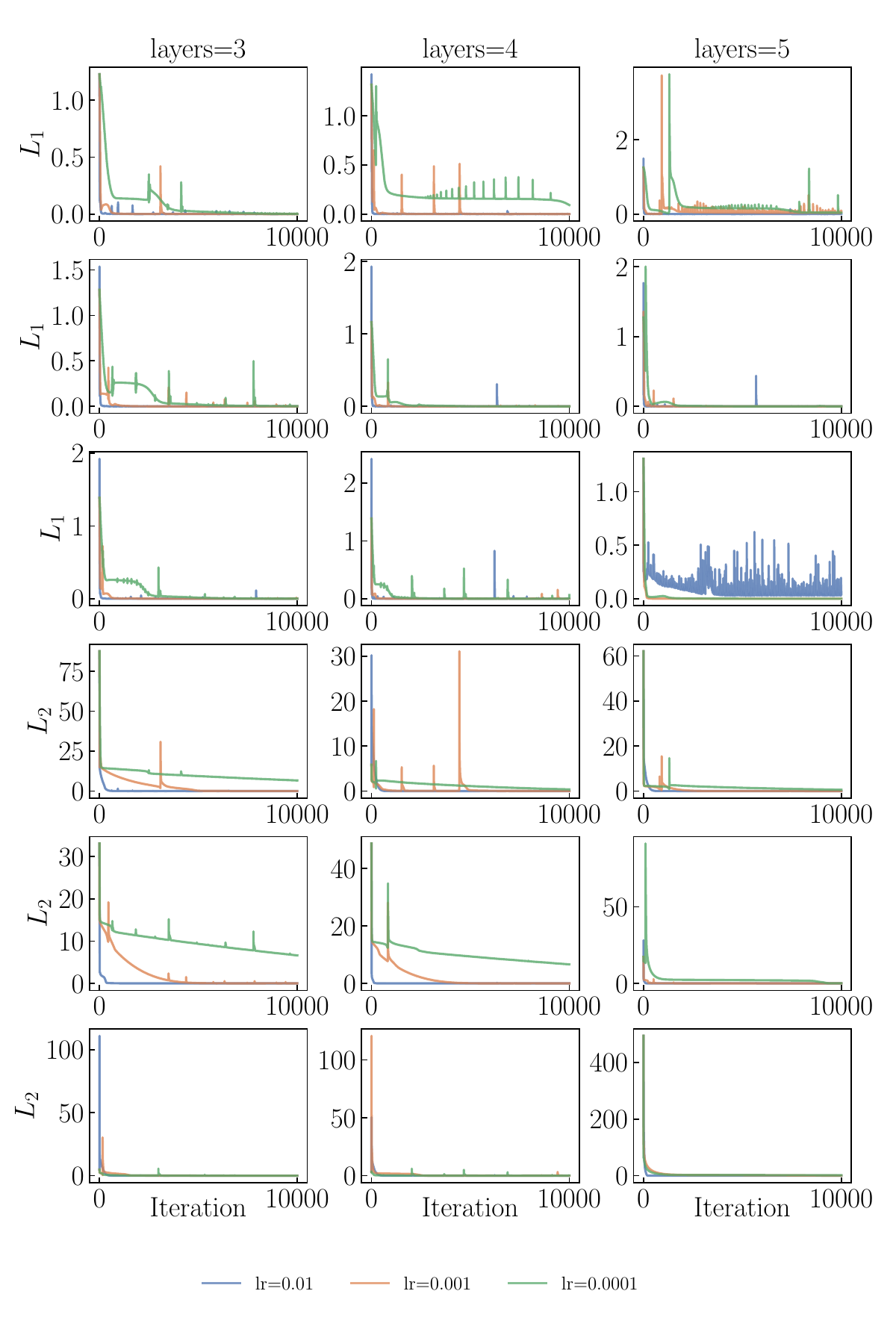}
	\caption{{\color{black}Ablation study of training FitzHugh-Nagumo model in step 1. The results are reported under different layers, widths of neural networks and learning rates (lr) of training loss. The first and fourth rows correspond to width of hidden layer $H=40$, the second and fifth rows correspond to $H=60$, the third and sixth rows correspond to $H=80$.  In each figure, the results are reported based on the average value under $5$ random seeds.} }
	
	\label{fig_ablation_step1}
\end{figure} 

\subsection{Ablation study for Step 2: identifying the equant}

We next examine the robustness of Step 2, where the invertible neural network is used to identify the equant. As discussed in the main text,  the equant is selected by minimising the normalised areal non-uniformity $\hat{K}$, which provides an additional geometric and dynamical criterion.

To test whether the equant is robustly identified, we train the invertible neural network under different random seeds in $\{1,2,3,4,5\}$. We repeat this experiment for the FitzHugh--Nagumo system with different values of the time-scale parameter $\varepsilon$. For each $\varepsilon$, we report the optimal  equant and the corresponding optimal value $\hat{K}^{\ast}$.

The results are summarised in Table~\ref{tab_ablation_step2}. Across different training runs, the learned equants converge to nearly identical locations for each parameter value. Moreover, as $\varepsilon$ approaches the Hopf transition, the optimal equant moves continuously towards the equilibrium, while $\hat{K}^{\ast}$ decreases. This behaviour is consistent with the results reported in the main text and supports the robustness of the equant as a geometric diagnostic.

\begin{table}[htp]
	\centering
	\caption{{\color{black}Ablation results for Step 2 on the FitzHugh--Nagumo system. For each value of the time-scale parameter $\varepsilon$, the table reports the learned  equant $\rvx^\ast$ and the corresponding optimal equant non-uniformity $\hat{K}^{\ast}$. Different random seeds yield nearly identical results.}}
	\label{tab_ablation_step2}
	\begin{tabular}{c c c}
		\hline
		$\varepsilon$ & equant $\rvx^\ast$ & $\hat{K}^{\ast}$ \\
		\hline
		1.41 & $(-0.3450,\; 0.5554)$ & $5.2453\times 10^{-4}$ \\
		1.42 & $(-0.3373,\; 0.4405)$ & $4.7338\times 10^{-4}$ \\
		1.43 & $(-0.3380,\; 0.3262)$ & $4.2839\times 10^{-4}$ \\
		1.44 & $(-0.3371,\; 0.2180)$ & $3.9994\times 10^{-4}$ \\
		1.45 & $(-0.3347,\; 0.1078)$ & $3.6858\times 10^{-4}$ \\
		1.46 & $(-0.3328,\; 0.0018)$ & $3.3322\times 10^{-4}$ \\
		1.47 & $(-0.3283,\; -0.0994)$ & $2.7469\times 10^{-4}$ \\
		1.48 & $(-0.3174,\; -0.1983)$ & $2.5048\times 10^{-4}$ \\
		1.49 & $(-0.3134,\; -0.2964)$ & $2.4673\times 10^{-4}$ \\
		1.50 & $(-0.3073,\; -0.3922)$ & $2.4532\times 10^{-4}$ \\
		1.51 & $(-0.3015,\; -0.4857)$ & $2.2077\times 10^{-4}$ \\
		1.52 & $(-0.2968,\; -0.5741)$ & $1.7187\times 10^{-4}$ \\
		1.53 & $(-0.2952,\; -0.6644)$ & $1.1765\times 10^{-4}$ \\
		1.54 & $(-0.2979,\; -0.7503)$ & $1.2706\times 10^{-4}$ \\
		\hline
	\end{tabular}
\end{table}

\subsection{Role of preprocessing and auxiliary regularisation}

We also clarify the role of the affine preprocessing $(\rmA,\rvb)$ and the auxiliary regularisation terms. The affine preprocessing is used only for numerical conditioning. In particular, it is chosen so that the origin lies inside the limit cycle in the relevant two-dimensional projections, and $\rmA$ rescales small-amplitude oscillations so that the learning algorithm can better resolve the oscillatory structure. Importantly, this preprocessing is applied consistently to both the limit-cycle trajectory and the vector field. For the affine coordinate transformation
\[
\rvy=\rmA\rvx+\rvb,
\]
the transformed dynamics are
\[
\dot{\rvy}=\rmA\rmF(\rmA^{-1}(\rvy-\rvb)).
\]
Thus, the preprocessing is simply an invertible affine change of coordinates, rather than an artificial deformation of the dynamics. After training, the learned equant and the associated functions are mapped back to the original coordinate system using the inverse transformation. The ablation results show that different admissible choices of $(\rmA,\rvb)$ lead to consistent learned phase maps and the equants, provided that the transformed trajectory is numerically well conditioned.

Finally, we specify the criterion for using the auxiliary regularisation terms $L_{\mathrm{ptp}}$ and $L_{\mathrm{norm}}$. These terms are not part of the default training objective. They are introduced only when the basic loss converges to a trivial or degenerate solution, such as $w=0$ or a collapsed latent representation. In such cases, $L_{\mathrm{ptp}}$ encourages the learned phase to cover a full $2\pi$ period along the limit cycle, while $L_{\mathrm{norm}}$ prevents the equant from escaping to an unbounded region. When the basic loss converges successfully without these auxiliary terms, they are not used. Therefore, these regularisation terms serve only to improve numerical convergence and do not determine the final equant.

Overall, these ablation studies show that the learned dynamical clock in Step 1 and the  equant in Step 2 are robust to changes in random seed, architecture, optimisation hyperparameters, and preprocessing choices. This supports the interpretation that the  equant is a stable geometric feature selected by the areal-uniformity criterion, rather than a training artifact.

\section{Look-up table for common systems}	
	In this Section, we systematically summarise the phase response curves (PRCs) defined in Eq.~\eqref{eq PRC} of all the studied systems based on their type and shape. The PRC of the system offers profound insights into how external stimuli modulate the phase of oscillation, either advancing or delaying it, thereby providing a key mechanism for understanding dynamic responses to perturbations.
	
	For vector-valued dynamics, the PRC is also a vector-valued function. In this study, we focus on a potentially controllable variable for each system, which lays the groundwork for future regulation of this variable using our results. We emphasise that the PRC for any other variable in the dynamics can also be determined using our framework.
	
	For quick reference, we provide all the results in Tables~\ref{table1},\ref{table2}, allowing for easy querying.

	\clearpage



\clearpage
	\begin{table*}[t]
		\centering
		\caption{Results of the phase response curve (PRC) with respect to specific variable.} \label{table1}
		\resizebox{\linewidth}{!}{
			\begin{tabular}{cccc}
				\toprule
				\toprule
				\multicolumn{1}{c}{Systems}   & \multicolumn{1}{c}{Variable  } & \multicolumn{1}{c}{Meaning} & \multicolumn{1}{c}{Plot}  \\
				\hline
				\hline
				{Morris-Lecar~\cite{morris}}   &  \makecell[c]{$n$ in Eq.~\eqref{eq morris}}
				
				&  \makecell[c]{Fraction of open\\ K$^{+}$ channels} & \begin{minipage}[H]{0.22\columnwidth}
					\centering
					\raisebox{-.3\height}{\includegraphics[width=0.55\linewidth]{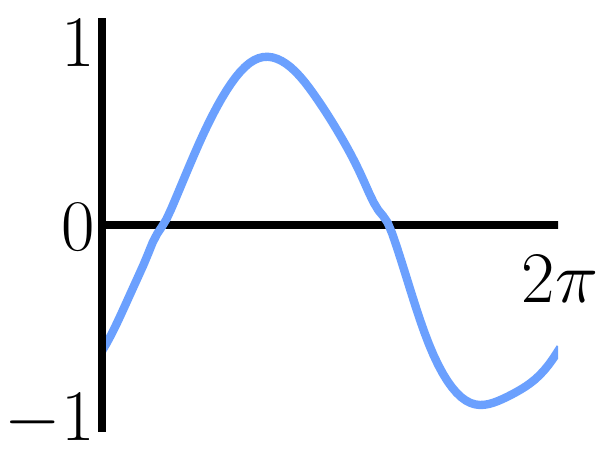}}
				\end{minipage} \\
				\hline 
				{FitzHugh-Nagumo~\cite{fitzhugh1969mathematical}}   & \makecell[c]{ $x$ in Eq.~\eqref{eq FithHugh-Nagumo}}
				
				& \makecell[c]{Membrane potential}& \begin{minipage}[H]{0.22\columnwidth}
					\centering
					\raisebox{-.3\height}{\includegraphics[width=0.55\linewidth]{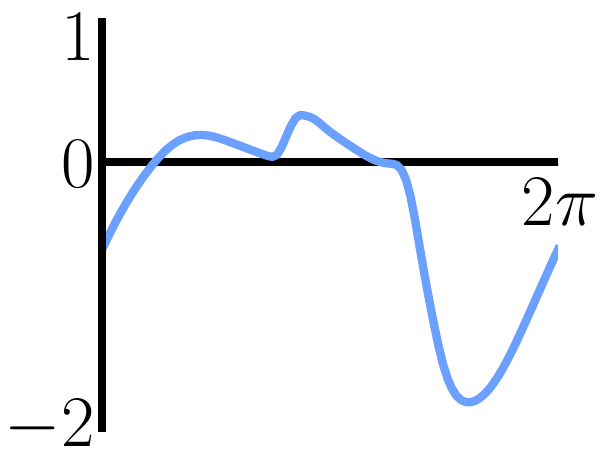}}
				\end{minipage} \\
				\hline 
				{Wilson Cowan~\cite{wilson1972excitatory}}   & \makecell[c]{$E$ in Eq.~\eqref{eq Wilson Cowan}}
				
				& \makecell[c]{Proportion of \\excitatory cells firing \\per unit time} & \begin{minipage}[H]{0.22\columnwidth}
					\centering
					\raisebox{-.3\height}{\includegraphics[width=0.55\linewidth]{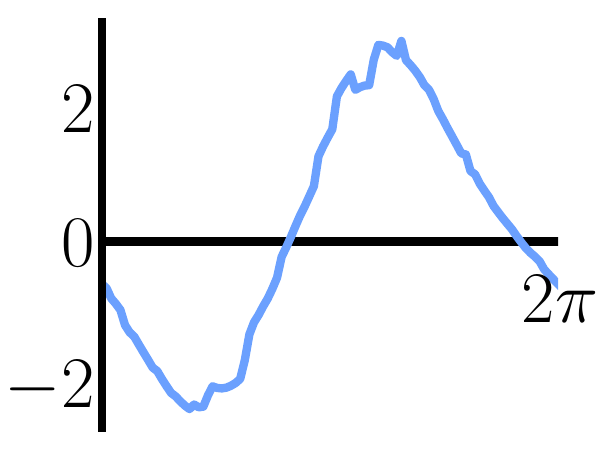}}
				\end{minipage} \\
				\hline 
				{Selkov~\cite{sel1968self}}   &  \makecell[c]{$x$ in Eq.~\eqref{eq Selkov}}
				
				& \makecell[c]{
					Relative concentration\\ of the ATP} & \begin{minipage}[H]{0.22\columnwidth}
					\centering
					\raisebox{-.3\height}{\includegraphics[width=0.55\linewidth]{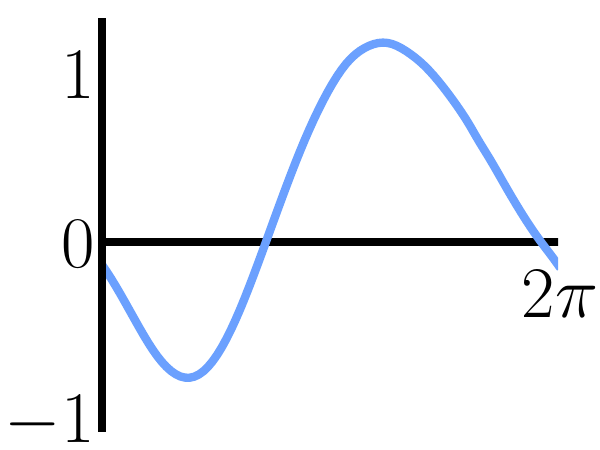}}
				\end{minipage} \\
				\hline 

                 {Lotka–Volterra~\cite{bomze1983lotka}}   &  \makecell[c]{$N$ in Eq.~\eqref{eq Lotka-Volterra} }
				
				& \makecell[c]{
					One of the \\competing species\\} & \begin{minipage}[H]{0.22\columnwidth}
					\centering
					\raisebox{-.3\height}{\includegraphics[width=0.55\linewidth]{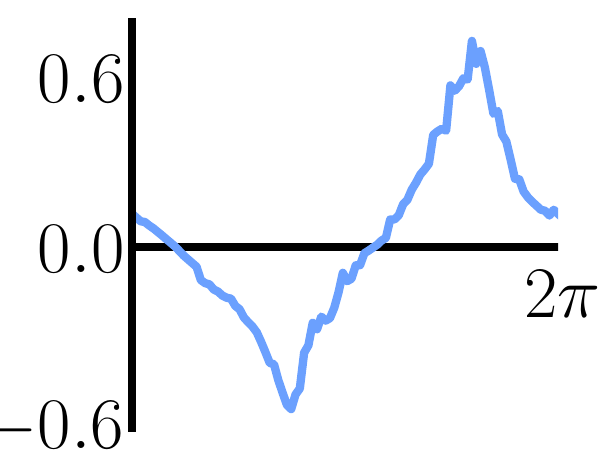}}
				\end{minipage} \\
				\hline 
                {Laser~\cite{dolcemascolo2020effective}}   & \makecell[c]{$w$ in Eq.~\eqref{eq laser}}
				
				& \makecell[c]{High-pass filtered \\feedback current  } & \begin{minipage}[H]{0.22\columnwidth}
					\centering
					\raisebox{-.3\height}{\includegraphics[width=0.55\linewidth]{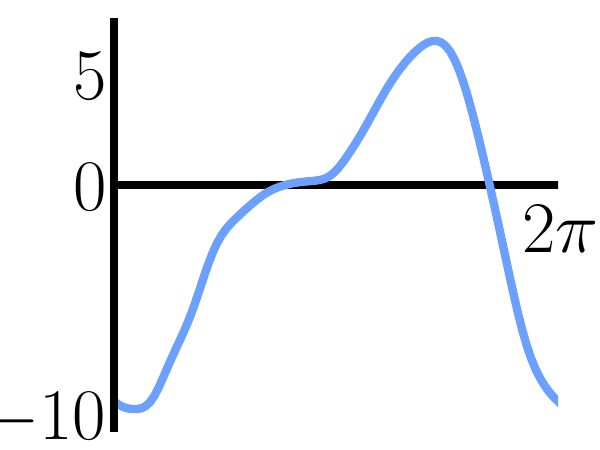}}
				\end{minipage} \\
                \hline
                {Mitotic~\cite{mitotic}}   &  \makecell[c]{$C$ in Eq.~\eqref{eq Mitotic}}
				
				&\makecell[c]{Cyclin concentration} & \begin{minipage}[H]{0.22\columnwidth}
					\centering
					\raisebox{-.3\height}{\includegraphics[width=0.55\linewidth]{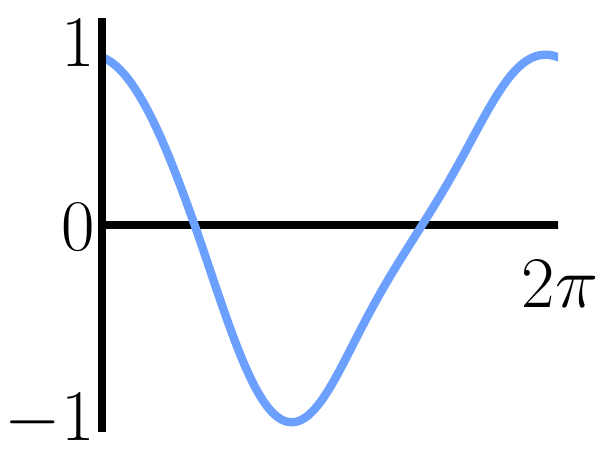}}
				\end{minipage} \\
	
				\hline  
				\bottomrule
			\end{tabular}
		}
	\end{table*}

     \clearpage
  \begin{table*}[t]
		\centering
		\caption{Results of the phase response curve (PRC) with respect to specific variable.} \label{table2}
		\resizebox{\linewidth}{!}{
			\begin{tabular}{cccc}
				\toprule
				\toprule
				\multicolumn{1}{c}{Systems}   & \multicolumn{1}{c}{Variable  } & \multicolumn{1}{c}{Meaning} & \multicolumn{1}{c}{Plot}  \\
				\hline
				\hline 
				{Thalamic Neuron~\cite{rubin2004high}}   &  
				\makecell[c]{$r$ in Eq.~\eqref{eq Thalamic}}
				
				& \makecell[c]{Gating variable\\ of the neuron} & \begin{minipage}[H]{0.22\columnwidth}
					\centering
					\raisebox{-.3\height}{\includegraphics[width=0.55\linewidth]{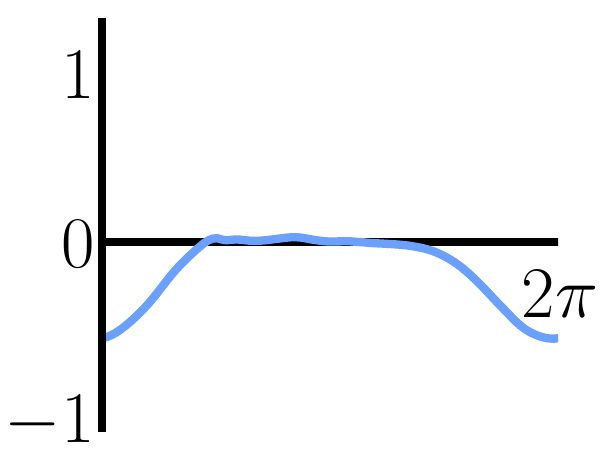}}
				\end{minipage} \\
				\hline 
				{Escherichia coli~\cite{garcia2004modeling}}   & \makecell[c]{$S$ in Eq.~\eqref{eq Escherichia coli}}
				
				& \makecell[c]{Concentration\\ of autoinducer }& \begin{minipage}[H]{0.22\columnwidth}
					\centering
					\raisebox{-.3\height}{\includegraphics[width=0.55\linewidth]{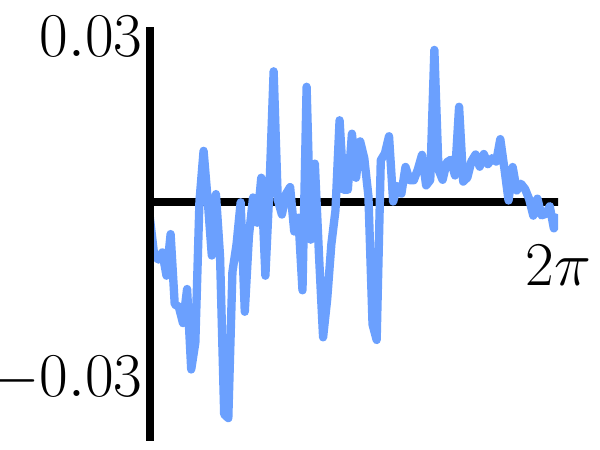}}
				\end{minipage} \\
				\hline 
				
				{Skeleton for Cdks~\cite{gerard2011skeleton}}   &  \makecell[c]{E2F in Eq.~\eqref{eq Skeleton}}
				
				& \makecell[c]{Transcription\\ factor E2F} & \begin{minipage}[H]{0.22\columnwidth}
					\centering
					\raisebox{-.3\height}{\includegraphics[width=0.55\linewidth]{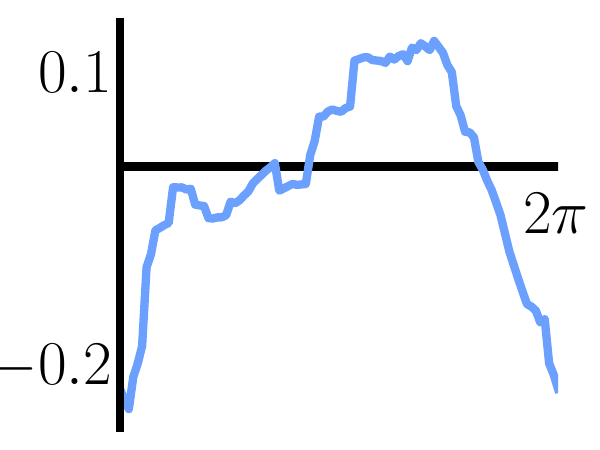}}
				\end{minipage} \\
				\hline  
				{\makecell[c]{Mammalian \\Circadian Clock~\cite{leloup2003toward}}}   &  \makecell[c]{$P_{CP}$}
				
				& \makecell[c]{ Nonphosphorylated\\ protein  PER \\ in the cytosol} & \begin{minipage}[H]{0.22\columnwidth}
					\centering
					\raisebox{-.3\height}{\includegraphics[width=0.55\linewidth]{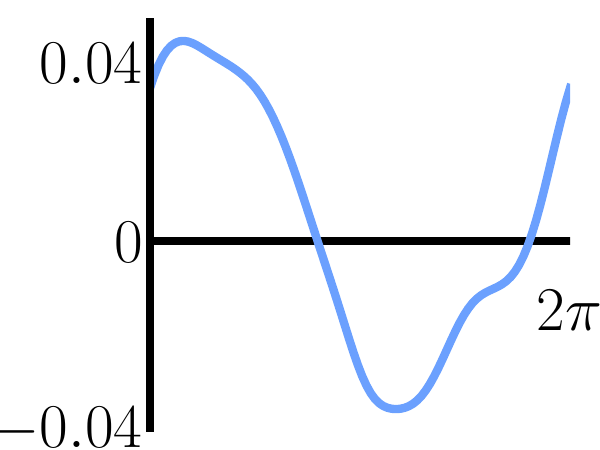}}
				\end{minipage} \\
				\hline  
				{Cdks for cell cycle~\cite{gerard2009temporal}}   & \makecell[c]{ATR} 
				
				&  \makecell[c]{ ATR kinase}  & \begin{minipage}[H]{0.22\columnwidth}
					\centering
					\raisebox{-.3\height}{\includegraphics[width=0.55\linewidth]{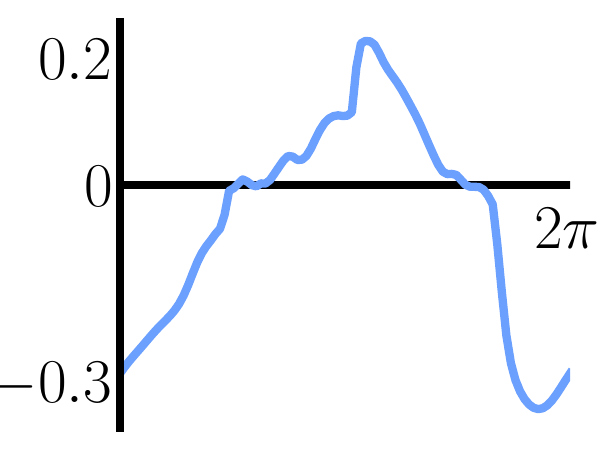}}
				\end{minipage} \\
				\hline  
				\bottomrule
			\end{tabular}
		}
	\end{table*}



\stopcontents[supp]



\bibliographystyle{apsrev4-2}     
\bibliography{main}

@Misc{methods,
  note = {Materials and methods are available as supplementary material},
}

@article{winfree1974patterns,
	title={Patterns of phase compromise in biological cycles},
	author={Winfree, A. T.},
	journal={J. Math. Biol.},
	volume={1},
	number={1},
	pages={73--93},
	year={1974},
	publisher={Springer}
}

@book{oksendal2013stochastic,
	title={Stochastic Differential Equations: An Introduction with Applications},
	author={Oksendal, B.},
	year={2013},
	publisher={Springer Science \& Business Media},
}

@inproceedings{dinh2016density,
	title={Density estimation using Real NVP},
	author={Dinh, L. and Sohl-Dickstein, J. and Bengio, S.},
	booktitle={International Conference on Learning Representations},
	year={2016},
	publisher = {OpenReview},
	address={San Juan, Puerto Rico}
}

@inproceedings{paszke2019pytorch,
	title={Pytorch: An imperative style, high-performance deep learning library},
	author={Paszke, A. and Gross, S. and Massa, F. and Lerer, A. and Bradbury, J. and Chanan, G. and Killeen, T. and Lin, Z. and Gimelshein, N. and Antiga, L. and others},
	booktitle={Advances in Neural Information Processing Systems},
	volume={32},
	year={2019},
	publisher = {Curran Associates, Inc.},
	location={Vancouver},
}

@article{brown2004phase,
	title={On the phase reduction and response dynamics of neural oscillator populations},
	author={Brown, E. and Moehlis, J. and Holmes, P.},
	journal={Neur. Comput.},
	volume={16},
	number={4},
	pages={673--715},
	year={2004},
	publisher={MIT Press},
}

@article{goldobin2010dynamics,
	title={Dynamics of limit-cycle oscillators subject to general noise},
	author={Goldobin, D. S. and Teramae, J. and Nakao, H. and Ermentrout, G. B.},
	journal={Phys. Rev. Lett.},
	volume={105},
	number={15},
	pages={154101},
	year={2010},
	publisher={APS},
}

@article{NOVICENKO20121090,
	title = {Phase reduction of weakly perturbed limit cycle oscillations in time-delay systems},
	journal = {Phys. D},
	volume = {241},
	number = {12},
	pages = {1090-1098},
	year = {2012},
	issn = {0167-2789},
	author = {Novičenko, V. and Pyragas, K.},
}

@article{guckenheimer1975isochrons,
	title={Isochrons and phaseless sets},
	author={Guckenheimer, J.},
	journal={J. Math. Biol.},
	volume={1},
	number={3},
	pages={259--273},
	year={1975},
	publisher={Springer}
}

@book{winfree1980geometry,
	title={The Geometry of Biological Time},
	author={Winfree, A. T.},
	volume={2},
	year={1980},
	publisher={Springer}
}

@article{paszke2017automatic,
	title={Automatic differentiation in pytorch},
	author={Paszke, A. and Gross, S. and Chintala, S. and Chanan, G. and Yang, E. and DeVito, Z. and Lin, Z. and Desmaison, A. and Antiga, L. and Lerer, A.},
	journal={NeurIPS 2017 Workshop Autodiff Decision Program Chairs},
	year={2017}
}

@incollection{kuramoto1984chemical,
	title={Chemical turbulence},
	author={Kuramoto, Y.},
	booktitle={Chemical oscillations, waves, and turbulence},
	pages={111--140},
	year={1984},
	publisher={Springer},
address={Heidelberg}
}

@article{kuramoto2019concept,
	title={On the concept of dynamical reduction: the case of coupled oscillators},
	author={Kuramoto, Y. and Nakao, H.},
	journal={Philos. Trans. Royal Soc. A},
	volume={377},
	number={2160},
	pages={20190041},
	year={2019},
	publisher={The Royal Society Publishing}
}

@article{nakao2016phase,
	title={Phase reduction approach to synchronisation of nonlinear oscillators},
	author={Nakao, H.},
	journal={Contemp. Phys.},
	volume={57},
	number={2},
	pages={188--214},
	year={2016},
	publisher={Taylor \& Francis}
}

@article{rosenblum1996phase,
	title={Phase synchronization of chaotic oscillators},
	author={Rosenblum, M. G. and Pikovsky, A. S. and Kurths, J.},
	journal={Phys. Rev. Lett.},
	volume={76},
	number={11},
	pages={1804},
	year={1996},
	publisher={APS}
}

@incollection{rosenblum2001phase,
	title={Phase synchronization: from theory to data analysis},
	author={Rosenblum, M. and Pikovsky, A. and Kurths, J. and Sch{\"a}fer, C. and Tass, P. A.},
	booktitle={Handbook of biological physics},
	volume={4},
	pages={279--321},
	year={2001},
	publisher={Elsevier},
address={Amsterdam}

}

@book{strogatz2018nonlinear,
	title={Nonlinear dynamics and chaos: with applications to physics, biology, chemistry, and engineering},
	author={Strogatz, S. H.},
	year={2018},
	publisher={CRC press},
address={Florida}

}

@article{gray1994synchronous,
	title={Synchronous oscillations in neuronal systems: mechanisms and functions},
	author={Gray, C. M.},
	journal={J. Comput. Neurosci.},
	volume={1},
	number={1},
	pages={11--38},
	year={1994},
	publisher={Springer}
}

@article{varela2001brainweb,
	title={The brainweb: phase synchronization and large-scale integration},
	author={Varela, F. and Lachaux, J. and Rodriguez, E. and Martinerie, J.},
	journal={Nat. Rev. Neurosci.},
	volume={2},
	number={4},
	pages={229--239},
	year={2001},
	publisher={Nature Publishing Group}
}

@article{palva2005phase,
	title={Phase synchrony among neuronal oscillations in the human cortex},
	author={Palva, J. M. and Palva, S. and Kaila, K.},
	journal={J. Neurosci.},
	volume={25},
	number={15},
	pages={3962--3972},
	year={2005},
	publisher={Soc Neuroscience}
}

@article{sauseng2008does,
	title={What does phase information of oscillatory brain activity tell us about cognitive processes?},
	author={Sauseng, P. and Klimesch, W.},
	journal={Neurosci. Biobehav. Rev.},
	volume={32},
	number={5},
	pages={1001--1013},
	year={2008},
	publisher={Elsevier}
}

@article{acebron2005kuramoto,
	title={The Kuramoto model: A simple paradigm for synchronization phenomena},
	author={Acebr{\'o}n, J. A. and Bonilla, L. L. and Vicente, C. J. P. and Ritort, F. and Spigler, R.},
	journal={Rev. Mod. Phys.},
	volume={77},
	number={1},
	pages={137},
	year={2005},
	publisher={APS}
}

@article{fitzhugh1969mathematical,
	title={Mathematical models of excitation and propagation in nerve},
	author={FitzHugh, R.},
	journal={Biol. Eng.},
	volume={9},
	pages={1--85},
	year={1969},
	publisher={McCraw Hill}
}

@article{garcia2004modeling,
	title={Modeling a synthetic multicellular clock: repressilators coupled by quorum sensing},
	author={Garcia-Ojalvo, J. and Elowitz, M. B. and Strogatz, S. H.},
	journal={Proceedings of the National Academy of Sciences},
	volume={101},
	number={30},
	pages={10955--10960},
	year={2004},
	publisher={National Acad Sciences}
}

@article{ermentrout1996type,
	title={Type I membranes, phase resetting curves, and synchrony},
	author={Ermentrout, B.},
	journal={Neur. Comput.},
	volume={8},
	number={5},
	pages={979--1001},
	year={1996},
	publisher={MIT Press One Rogers Street, Cambridge, MA 02142-1209, USA journals-info~…}
}

@article{ermentrout2003simulating,
	title={Simulating, analyzing, and animating dynamical systems: a guide to XPPAUT for researchers and students},
	author={Ermentrout, B. and Mahajan, A.},
	journal={Appl. Mech. Rev.},
	volume={56},
	number={4},
	pages={B53--B53},
	year={2003}
}

@article{pietras2019network,
	title={Network dynamics of coupled oscillators and phase reduction techniques},
	author={Pietras, B. and Daffertshofer, A.},
	journal={Phys, Rep.},
	volume={819},
	pages={1--105},
	year={2019},
	publisher={Elsevier}
}

@article{iten2020discovering,
	title={Discovering physical concepts with neural networks},
	author={Iten, R. and Metger, T. and Wilming, H. and Del Rio, L. and Renner, R.},
	journal={Phys. Rev. Lett.},
	volume={124},
	number={1},
	pages={010508},
	year={2020},
	publisher={APS}
}

@article{wang2023scientific,
	title={Scientific discovery in the age of artificial intelligence},
	author={Wang, H. and Fu, T. and Du, Y. and Gao, W. and Huang, K. and Liu, Z. and Chandak, P. and Liu, S. and Van Katwyk, P. and Deac, A. and others},
	journal={Nature},
	volume={620},
	number={7972},
	pages={47--60},
	year={2023},
	publisher={Nature Publishing Group UK London}
}

@article{davies2021advancing,
	title={Advancing mathematics by guiding human intuition with AI},
	author={Davies, A. and Veli{\v{c}}kovi{\'c}, P. and Buesing, L. and Blackwell, S. and Zheng, D. and Toma{\v{s}}ev, N. and Tanburn, R. and Battaglia, P. and Blundell, C. and Juh{\'a}sz, A. and others},
	journal={Nature},
	volume={600},
	number={7887},
	pages={70--74},
	year={2021},
	publisher={Nature Publishing Group}
}

@article{pecora1998master,
	title={Master stability functions for synchronized coupled systems},
	author={Pecora, L. M. and Carroll, T. L.},
	journal={Phys. Rev. Lett.},
	volume={80},
	number={10},
	pages={2109},
	year={1998},
	publisher={APS}
}

@book{hamilton2020time,
	title={Time series analysis},
	author={Hamilton, J. D.},
	year={2020},
	publisher={Princeton university press},
address={Princeton}
}

@article{nakao2007noise,
	title={Noise-Induced Synchronization and Clustering in Ensembles of Uncoupled Limit-Cycle Oscillators},
	author={Nakao, H. and Arai, K. and Kawamura, Y.},
	journal={Phys. Rev. Lett.},
	volume={98},
	number={18},
	pages={184101},
	year={2007}
}

@article{li2007effects,
	title={Effects of chemical synapses on the enhancement of signal propagation in coupled neurons near the canard regime},
	author={Li, X. and Wang, J. and Hu, W.},
	journal={Phys. Rev. E},
	volume={76},
	number={4},
	pages={041902},
	year={2007},
	publisher={APS}
}

@article{mcmillen2002synchronizing,
	title={Synchronizing genetic relaxation oscillators by intercell signaling},
	author={McMillen, D. and Kopell, N. and Hasty, J. and Collins, J.},
	journal={Proc. Natl. Acad. Sci.},
	volume={99},
	number={2},
	pages={679--684},
	year={2002},
	publisher={National Acad Sciences}
}

@article{zhong2023modulating,
	title={Modulating biological rhythms: A noncomputational strategy harnessing nonlinearity and decoupling frequency and amplitude},
	author={Zhong, Z. and Lin, W. and Qin, B.},
	journal={Phys. Rev. Lett.},
	volume={131},
	number={13},
	pages={138401},
	year={2023},
	publisher={APS}
}

@article{monga2019optimal,
	title={Optimal phase control of biological oscillators using augmented phase reduction},
	author={Monga, Bharat and Moehlis, Jeff},
	journal={Biol. Cybern.},
	volume={113},
	number={1},
	pages={161--178},
	year={2019},
	publisher={Springer}
}

@article{monga2019phase,
	title={Phase reduction and phase-based optimal control for biological systems: a tutorial},
	author={Monga, B. and Wilson, D. and Matchen, T. and Moehlis, J.},
	journal={Biol. Cybern.},
	volume={113},
	number={1},
	pages={11--46},
	year={2019},
	publisher={Springer}
}

@article{wang2023desynchronizing,
	title={Desynchronizing oscillators coupled in multi-cluster networks through adaptively controlling partial networks},
	author={Wang, K. and Yang, L. and Zhou, S. and Lin, W.},
	journal={Chaos},
	volume={33},
	number={9},
	year={2023},
	publisher={AIP Publishing}
}

@article{kim2012mechanism,
	title={A mechanism for robust circadian timekeeping via stoichiometric balance},
	author={Kim, J. and Forger, D. B.},
	journal={Mol. Syst. Biol.},
	volume={8},
	number={1},
	pages={630},
	year={2012},
	publisher={John Wiley \& Sons, Ltd Chichester, UK}
}

@article{zhang2024learning,
	title={Learning Hamiltonian neural Koopman operator and simultaneously sustaining and discovering conservation laws},
	author={Zhang, J. and Zhu, Q. and Lin, W.},
	journal={Phys. Rev. Res.},
	volume={6},
	number={1},
	pages={L012031},
	year={2024},
	publisher={APS}
}

@book{glass1988clocks,
	title={From clocks to chaos: The rhythms of life},
	author={Glass, L. and Mackey, M. C.},
	year={1988},
	publisher={Princeton University Press},
address={Princeton}
}

@article{taylor2009dynamical,
	title={Dynamical quorum sensing and synchronization in large populations of chemical oscillators},
	author={Taylor, A. F. and Tinsley, M. R. and Wang, F. and Huang, Z. and Showalter, K.},
	journal={Science},
	volume={323},
	number={5914},
	pages={614--617},
	year={2009},
	publisher={American Association for the Advancement of Science}
}

@article{zhang2024machine,
	title={Machine-learning-coined noise induces energy-saving synchrony},
	author={Zhang, J. and Yang, L. and Zhu, Q. and Grebogi, C. and Lin, W.},
	journal={Phys. Rev. E},
	volume={110},
	number={1},
	pages={L012203},
	year={2024},
	publisher={APS}
}

@article{mitotic,
	title={A minimal cascade model for the mitotic oscillator involving cyclin and cdc2 kinase.},
	author={Goldbeter, A.},
	journal={Proc. Natl. Acad. Sci.},
	volume={88},
	number={20},
	pages={9107--9111},
	year={1991},
	publisher={National Acad Sciences}
}

@article{morris,
	title={Voltage oscillations in the barnacle giant muscle fiber},
	author={Morris, C. and Lecar, H.},
	journal={Biophys. J.},
	volume={35},
	number={1},
	pages={193--213},
	year={1981},
	publisher={Elsevier}
}

@incollection{evans2013mechanical,
	title={Mechanical astronomy: A route to the ancient discovery of epicycles and eccentrics},
	author={Evans, J. and Cannan, C. C.},
	booktitle={From Alexandria, through Baghdad: Surveys and studies in the ancient Greek and medieval Islamic mathematical sciences in honor of JL Berggren},
	pages={145--174},
	year={2013},
	publisher={Springer},
address={Heidelberg}
}

@book{toomer1998ptolemy,
	title={Ptolemy’s Almagest},
	author={Toomer, G. J.},
	year={1998},
	publisher={Princeton University Press},
address={Princeton}
}

@article{gonze2002robustness,
	title={Robustness of circadian rhythms with respect to molecular noise},
	author={Gonze, D. and Halloy, J. and Goldbeter, A.},
	journal={Proc. Natl. Acad. Sci.},
	volume={99},
	number={2},
	pages={673--678},
	year={2002},
	publisher={National Acad Sciences}
}

@article{yan2019robust,
	title={Robust synchronization of the cell cycle and the circadian clock through bidirectional coupling},
	author={Yan, J. and Goldbeter, A.},
	journal={J. Roy. Soc. Interface},
	volume={16},
	number={158},
	pages={20190376},
	year={2019},
	publisher={The Royal Society}
}

@article{altinok2007cell,
	title={A cell cycle automaton model for probing circadian patterns of anticancer drug delivery},
	author={Altinok, A. and L{\'e}vi, F. and Goldbeter, A.},
	journal={Adv. Drug Deliv. Rev},
	volume={59},
	number={9-10},
	pages={1036--1053},
	year={2007},
	publisher={Elsevier}
}

@article{levi2008implications,
	title={Implications of circadian clocks for the rhythmic delivery of cancer therapeutics},
	author={L{\'e}vi, F. and Altinok, A. and Clairambault, J. and Goldbeter, A.},
	journal={Philos. Trans. Roy. Soc. A: Math. Phys. Eng. Sci.},
	volume={366},
	number={1880},
	pages={3575--3598},
	year={2008},
	publisher={The Royal Society London}
}

@article{goldbeter2022multi,
	title={Multi-synchronization and other patterns of multi-rhythmicity in oscillatory biological systems},
	author={Goldbeter, A. and Yan, J.},
	journal={Interface Focus},
	volume={12},
	number={3},
	pages={20210089},
	year={2022},
	publisher={The Royal Society}
}

@article{leloup2003toward,
	title={Toward a detailed computational model for the mammalian circadian clock},
	author={Leloup, J. and Goldbeter, A.},
	journal={Proc. Natl. Acad. Sci.},
	volume={100},
	number={12},
	pages={7051--7056},
	year={2003},
	publisher={National Acad Sciences}
}

@article{zhang2022neural,
	title={Neural stochastic control},
	author={Zhang, J. and Zhu, Q. and Lin, W.},
	journal={Advances in Neural Information Processing Systems},
	volume={35},
	pages={9098--9110},
	year={2022}
}

@article{volos2015memristor,
	title={Memristor: A New Concept in Synchronization of Coupled Neuromorphic Circuits.},
	author={Volos, C. K. and Kyprianidis, I. M. and Stouboulos, I. N. and Tlelo-Cuautle, E. and Vaidyanathan, S.},
	journal={J. Eng. Sci. Technol. Rev.},
	volume={8},
	number={2},
	year={2015}
}

@article{korotkov2019dynamics,
	title={The dynamics of ensemble of neuron-like elements with excitatory couplings},
	author={Korotkov, A. G. and Kazakov, A. O. and Levanova, T. A. and Osipov, G. V.},
	journal={Commun. Nonlinear Sci. Numer. Simul.},
	volume={71},
	pages={38--49},
	year={2019},
	publisher={Elsevier}
}

@book{kloeden1992stochastic,
	title={Stochastic differential equations},
	author={Kloeden, P. E. and Platen, E. and Platen, E.},
	year={1992},
	publisher={Springer}
}

@article{rubin2004high,
	title={High frequency stimulation of the subthalamic nucleus eliminates pathological thalamic rhythmicity in a computational model},
	author={Rubin, J. E. and Terman, D.},
	journal={J. Comput. Neurosci.},
	volume={16},
	pages={211--235},
	year={2004},
	publisher={Springer}
}

@article{sel1968self,
	title={Self-Oscillations in Glycolysis 1. A Simple Kinetic Model},
	author={Selkov, E. E.},
	journal={Eur. J. Biochem.},
	volume={4},
	number={1},
	pages={79--86},
	year={1968},
	publisher={Wiley Online Library}
}

@article{wilson1972excitatory,
	title={Excitatory and inhibitory interactions in localized populations of model neurons},
	author={Wilson, H. R. and Cowan, J. D.},
	journal={Biophys. J.},
	volume={12},
	number={1},
	pages={1--24},
	year={1972},
	publisher={Elsevier}
}

@article{gerard2009temporal,
	title={Temporal self-organization of the cyclin/Cdk network driving the mammalian cell cycle},
	author={G{\'e}rard, Claude and Goldbeter, Albert},
	journal={Proc. Natl. Acad. Sci.},
	volume={106},
	number={51},
	pages={21643--21648},
	year={2009},
	publisher={National Acad Sciences}
}

@article{gerard2011skeleton,
	title={A skeleton model for the network of cyclin-dependent kinases driving the mammalian cell cycle},
	author={G{\'e}rard, C. and Goldbeter, A.},
	journal={Interface Focus},
	volume={1},
	number={1},
	pages={24--35},
	year={2011},
	publisher={The Royal Society}
}

@article{chen2018neural,
	title={Neural ordinary differential equations},
	author={Chen, R. T. Q. and Rubanova, Y. and Bettencourt, J. and Duvenaud, D. K.},
	journal={Advances in neural information processing systems},
	volume={31},
	year={2018}
}

@article{pathak2018model,
	title={Model-free prediction of large spatiotemporally chaotic systems from data: A reservoir computing approach},
	author={Pathak, J. and Hunt, B. and Girvan, M. and Lu, Z. and Ott, E.},
	journal={Phys. Rev. Lett.},
	volume={120},
	number={2},
	pages={024102},
	year={2018},
	publisher={APS}
}

@article{brunton2016discovering,
	title={Discovering governing equations from data by sparse identification of nonlinear dynamical systems},
	author={Brunton, S. L. and Proctor, J. L. and Kutz, J. N.},
	journal={Proc. Natl. Acad. Sci. },
	volume={113},
	number={15},
	pages={3932--3937},
	year={2016},
	publisher={National Acad Sciences}
}

@article{berry1984quantal,
	title={Quantal phase factors accompanying adiabatic changes},
	author={Berry, M. V.},
	journal={Proc. Roy. Soc. A. Math. Phys. Sci.},
	volume={392},
	number={1802},
	pages={45--57},
	year={1984},
	publisher={The Royal Society London}
}

@article{simon1983holonomy,
	title={Holonomy, the quantum adiabatic theorem, and Berry's phase},
	author={Simon, B.},
	journal={Phys. Rev. Lett.},
	volume={51},
	number={24},
	pages={2167},
	year={1983},
	publisher={APS}
}

@article{aharonov1959significance,
	title={Significance of electromagnetic potentials in the quantum theory},
	author={Aharonov, Y. and Bohm, D.},
	journal={Phys. Rev.},
	volume={115},
	number={3},
	pages={485},
	year={1959},
	publisher={APS}
}

@article{hannay1985angle,
	title={Angle variable holonomy in adiabatic excursion of an integrable Hamiltonian},
	author={Hannay, J. H.},
	journal={J. Phys. A: Math. Gener.},
	volume={18},
	number={2},
	pages={221},
	year={1985},
	publisher={IOP Publishing}
}

@inproceedings{pancharatnam1956generalized,
	title={Generalized theory of interference, and its applications: Part I. Coherent pencils},
	author={Pancharatnam, S.},
	booktitle={Proc. Ind. Acad. Sci. A},
	volume={44},
	pages={247--262},
	year={1956},
	organization={Springer}
}

@book{wilczek1989geometric,
	title={Geometric phases in physics},
	author={Wilczek, F. and Shapere, A.},
	volume={5},
	year={1989},
	publisher={World Scientific},
address={Singapore}
}

@article{kepler1991geometric,
	title={Geometric phase shifts under adiabatic parameter changes in classical dissipative systems},
	author={Kepler, T. B. and Kagan, M. L.},
	journal={Phys. Rev. Lett.},
	volume={66},
	number={7},
	pages={847},
	year={1991},
	publisher={APS}
}

@article{dolcemascolo2020effective,
  title={Effective low-dimensional dynamics of a mean-field coupled network of slow-fast spiking lasers},
  author={Dolcemascolo, Axel and Miazek, Alexandre and Veltz, Romain and Marino, Francesco and Barland, St{\'e}phane},
  journal={Phys. Rev. E},
  volume={101},
  number={5},
  pages={052208},
  year={2020},
  publisher={APS}
}

@article{bomze1983lotka,
  title={Lotka-Volterra equation and replicator dynamics: a two-dimensional classification},
  author={Bomze, Immanuel M},
  journal={Biol. Cybern.},
  volume={48},
  number={3},
  pages={201--211},
  year={1983},
  publisher={Springer}
}

@article{hsieh1996phase,
  title={Phase-locked loop techniques. A survey},
  author={Hsieh, Guan-Chyun and Hung, James C},
  journal={IEEE Trans. Ind. Electron.},
  volume={43},
  number={6},
  pages={609--615},
  year={1996},
  publisher={IEEE}
}

@article{hellen2017electronic,
  title={Electronic circuit analog of synthetic genetic networks: Revisited},
  author={Hellen, Edward H and Kurths, J{\"u}rgen and Dana, Syamal K},
  journal={Eur. Phys. J. Spec. Top.},
  volume={226},
  pages={1811--1828},
  year={2017},
  publisher={Springer}
}

@article{brophy2014principles,
  title={Principles of genetic circuit design},
  author={Brophy, Jennifer AN and Voigt, Christopher A},
  journal={Nat. Meth.},
  volume={11},
  number={5},
  pages={508--520},
  year={2014},
  publisher={Nature Publishing Group US New York}
}

@article{holtz2010engineering,
  title={Engineering static and dynamic control of synthetic pathways},
  author={Holtz, William J and Keasling, Jay D},
  journal={Cell},
  volume={140},
  number={1},
  pages={19--23},
  year={2010},
  publisher={Elsevier}
}

@article{hellen2011electronic,
	title={An electronic analog of synthetic genetic networks},
	author={Hellen, Edward H and Volkov, Evgenii and Kurths, Jurgen and Dana, Syamal Kumar},
	journal={PLoS One},
	volume={6},
	number={8},
	pages={e23286},
	year={2011},
	publisher={Public Library of Science San Francisco, USA}
}

@article{teo2020merging,
	title={The merging of biological and electronic circuits},
	author={Teo, Jonathan JY and Sarpeshkar, Rahul},
	journal={Iscience},
	volume={23},
	number={11},
	year={2020},
	publisher={Elsevier}
}

@article{hodgkin1952quantitative,
  title={A quantitative description of membrane current and its application to conduction and excitation in nerve},
  author={Hodgkin, Alan L and Huxley, Andrew F},
  journal={The Journal of physiology},
  volume={117},
  number={4},
  pages={500},
  year={1952}
}

@article{hodgkin1952propagation,
  title={Propagation of electrical signals along giant nerve fibres},
  author={Hodgkin, Alan Lloyd and Huxley, Andrew Fielding},
  journal={Proceedings of the Royal Society of London. Series B-Biological Sciences},
  volume={140},
  number={899},
  pages={177--183},
  year={1952},
  publisher={The Royal Society London}
}

@article{chen2017weak,
  title={Weak synchronization and large-scale collective oscillation in dense bacterial suspensions},
  author={Chen, Chong and Liu, Song and Shi, Xia-qing and Chat{\'e}, Hugues and Wu, Yilin},
  journal={Nature},
  volume={542},
  number={7640},
  pages={210--214},
  year={2017},
  publisher={Nature Publishing Group UK London}
}

@article{yeung1999time,
  title={Time delay in the Kuramoto model of coupled oscillators},
  author={Yeung, MK Stephen and Strogatz, Steven H},
  journal={Physical review letters},
  volume={82},
  number={3},
  pages={648},
  year={1999},
  publisher={APS}
}

@article{kiyohara2005novel,
  title={A novel mutation in kaiC affects resetting of the cyanobacterial circadian clock},
  author={Kiyohara, Yota B and Katayama, Mitsunori and Kondo, Takao},
  journal={Journal of bacteriology},
  volume={187},
  number={8},
  pages={2559--2564},
  year={2005},
  publisher={American Society for Microbiology}
}

@article{wilson2018greater,
	title={Greater accuracy and broadened applicability of phase reduction using isostable coordinates},
	author={Wilson, Dan and Ermentrout, Bard},
	journal={Journal of mathematical biology},
	volume={76},
	number={1},
	pages={37--66},
	year={2018},
	publisher={Springer}
}

@article{wilson2019augmented,
	title={Augmented phase reduction of (not so) weakly perturbed coupled oscillators},
	author={Wilson, Dan and Ermentrout, Bard},
	journal={SIAM Review},
	volume={61},
	number={2},
	pages={277--315},
	year={2019},
	publisher={SIAM}
}

@article{shirasaka2017phase,
	title={Phase reduction theory for hybrid nonlinear oscillators},
	author={Shirasaka, Sho and Kurebayashi, Wataru and Nakao, Hiroya},
	journal={Physical Review E},
	volume={95},
	number={1},
	pages={012212},
	year={2017},
	publisher={APS}
}

@article{botvinick2025invariant,
	title={Invariant measures in time-delay coordinates for unique dynamical system identification},
	author={Botvinick-Greenhouse, Jonah and Martin, Robert and Yang, Yunan},
	journal={Physical Review Letters},
	volume={135},
	number={16},
	pages={167202},
	year={2025},
	publisher={APS}
}

@article{hiruta2025autoencoder,
	title={Autoencoder for limit cycle in Kolmogorov flow},
	author={Hiruta, Yoshiki and Ishimoto, Kenta},
	journal={Journal of the Physical Society of Japan},
	volume={94},
	number={6},
	pages={064401},
	year={2025},
	publisher={The Physical Society of Japan}
}

@article{mauroy2012use,
	title={On the use of Fourier averages to compute the global isochrons of (quasi) periodic dynamics},
	author={Mauroy, Alexandre and Mezi{\'c}, Igor},
	journal={Chaos: An Interdisciplinary Journal of Nonlinear Science},
	volume={22},
	number={3},
	year={2012},
	publisher={AIP Publishing}
}

@article{mauroy2014global,
	title={Global isochrons and phase sensitivity of bursting neurons},
	author={Mauroy, Alexandre and Rhoads, Blane and Moehlis, Jeff and Mezic, Igor},
	journal={SIAM Journal on Applied Dynamical Systems},
	volume={13},
	number={1},
	pages={306--338},
	year={2014},
	publisher={SIAM}
}

@article{yawata2024phase,
	title={Phase autoencoder for limit-cycle oscillators},
	author={Yawata, Koichiro and Fukami, Kai and Taira, Kunihiko and Nakao, Hiroya},
	journal={Chaos: An Interdisciplinary Journal of Nonlinear Science},
	volume={34},
	number={6},
	year={2024},
	publisher={AIP Publishing}
}
\end{document}